%% file: main.tex
\documentclass[10pt]{article}
\usepackage{times}

\DeclareMathVersion{varnormal}
\DeclareMathVersion{varbold}

\SetSymbolFont{letters}{normal}{OML}{txmi}{m}{it}
\SetSymbolFont{letters}{bold}{OML}{txmi}{bx}{it}
\SetSymbolFont{letters}{varnormal}{OML}{mdugm}{m}{it}
\SetSymbolFont{letters}{varbold}{OML}{mdugm}{b}{it}
\SetSymbolFont{operators}{normal}{OT1}{txr}{m}{n}
\SetSymbolFont{operators}{bold}{OT1}{txr}{bx}{n}
\SetSymbolFont{operators}{varnormal}{OT1}{mdugm}{m}{n}
\SetSymbolFont{operators}{varbold}{OT1}{mdugm}{b}{n}
\SetSymbolFont{symbols}{normal}{OMS}{txsy}{m}{n}
\SetSymbolFont{symbols}{bold}{OMS}{txsy}{bx}{n}
\SetSymbolFont{symbols}{varnormal}{OMS}{mdugm}{m}{n}
\SetSymbolFont{symbols}{varbold}{OMS}{mdugm}{b}{n}
\SetSymbolFont{largesymbols}{normal}{OMX}{txex}{m}{n}
\SetSymbolFont{largesymbols}{bold}{OMX}{txex}{bx}{n}
\SetSymbolFont{largesymbols}{varnormal}{OMX}{mdugm}{m}{n}
\SetSymbolFont{largesymbols}{varbold}{OMX}{mdugm}{b}{n}

\SetMathAlphabet{\mathrm}{normal}{OT1}{txr}{m}{n}
\SetMathAlphabet{\mathrm}{bold}{OT1}{txr}{bx}{n}
\SetMathAlphabet{\mathrm}{varnormal}{OT1}{mdugm}{m}{n}
\SetMathAlphabet{\mathrm}{varbold}{OT1}{mdugm}{b}{n}

\SetMathAlphabet{\mathit}{normal}{OT1}{txr}{m}{it}
\SetMathAlphabet{\mathit}{bold}{OT1}{txr}{bx}{it}
\SetMathAlphabet{\mathit}{varnormal}{OT1}{mdugm}{m}{it}
\SetMathAlphabet{\mathit}{varbold}{OT1}{mdugm}{b}{it}

\usepackage[leqno]{amsmath} 
\usepackage{amsthm}\theoremstyle{plain}\newtheorem{remark}{Remark}
\usepackage{mathrsfs}
\usepackage{amssymb}
\usepackage{empheq}
\usepackage{subcaption}
\usepackage[paper=a4paper,top=36pt,left=100pt,right=100pt,bottom=50pt,foot=25pt]{geometry}
\usepackage{titlesec}
\usepackage{lineno}
\usepackage{hyperref}
\usepackage{nicefrac}
\usepackage{float}
\usepackage{subfloat}
\usepackage[dvipsnames]{xcolor}
\usepackage[numbers,compress]{natbib}

\titleformat{\section}{\normalfont\bfseries}{\thesection}{1em}{}
\titleformat{\subsection}{\normalfont}{\thesubsection}{1em}{}

\renewcommand\appendix{\par
  \setcounter{section}{0}
  \setcounter{subsection}{0}
  \setcounter{figure}{0}
  \setcounter{table}{0}
  \setcounter{equation}{0}
  \renewcommand\thesection{Appendix \Alph{section}}
  \renewcommand\thefigure{\Alph{section}\arabic{figure}}
  \renewcommand\thetable{\Alph{section}\arabic{table}}
  \renewcommand\theequation{\Alph{section}\arabic{equation}}  
}

\newcommand\supp{\par
  \setcounter{section}{0}
  \setcounter{subsection}{0}
  \setcounter{figure}{0}
  \setcounter{table}{0}
  \setcounter{equation}{0}  
  \renewcommand\thesection{S\arabic{section}}
  \renewcommand\thefigure{S-\arabic{figure}}
  \renewcommand\thetable{S-\arabic{table}}
  \renewcommand\theequation{S-\arabic{equation}}  
}

\title{\Large 
Strongly imposing the free surface boundary condition for 
wave equations with 
finite difference operators}
\author{Longfei Gao\thanks{Email address: longfei.gao@austin.utexas.edu}
\\ {\it \small Oden Institute for Computational Engineering and Sciences,} 
\\ {\it \small The University of Texas at Austin, Austin, TX 78712, USA}
}
\date{}

\begin{document}


\maketitle

\begin{abstract}

Acoustic and elastic wave equations are routinely used in geophysical and engineering studies to simulate the propagation of waves, with a broad range of applications, including seismology, near surface characterization, non-destructive structural evaluation, etc.
Finite difference methods remain popular choices for these simulations due to their simplicity and efficiency.
In particular, the family of finite difference methods based on the summation-by-parts operators and the simultaneous-approximation-terms technique have been proposed for these simulations, which offers great flexibility in addressing boundary and interface conditions.
For the applications mentioned above, surface of the earth 
is usually associated with the free surface boundary condition. 
In this study, we demonstrate that the weakly imposed free surface boundary condition through the simultaneous-approximation-terms technique can have issue when the source terms, which introduces abrupt disturbances to the wave field, are placed too close to the surface. 
%
In response, we propose to build the free surface boundary condition into the summation-by-parts finite difference operators and hence strongly and automatically impose the free surface boundary condition to address this issue.
The procedure is very simple for acoustic wave equation, requiring resetting a few rows and columns in the existing 
difference operators only. 
For the elastic wave equation, the procedure is more involved and requires special design of the grid layout and summation-by-parts operators that satisfy additional requirements, as revealed by the discrete energy analysis.
In both cases, the energy conserving property is preserved.
Numerical examples 
are presented to demonstrate the effectiveness of the proposed approach.

\end{abstract}

\providecommand{\keywords}[1]{{\small \textbf{Key words.} #1}}
\vspace{0.25em}

\providecommand{\MSCcodes}[1]{{\small \textbf{MSC codes.} #1}}
\MSCcodes{65M06, 35Q86, 86-08}

\section{Introduction}\label{sec_introduction}

The phenomena of wave propagation has broad applications in geoscience and civil engineering such as earthquake studies, site characterization, non-destructive structural testing, etc.; see \cite{aki2002quantitative, tromp2020seismic, kallivokas2013site, drinkwater2006ultrasonic} for a few examples.
Among the many discretization techniques that have been proposed for wave simulations, finite difference methods are popular choices for these applications due to their simplicity and efficiency. 
In particular, a family of finite difference methods based on the summation-by-parts (SBP) operators, which dates back to \cite{kreiss1974finite}, has found its success in wave simulations; see \cite{nilsson2007stable, appelo2009stable, sjogreen2012fourth, wang2017convergence, wang2019fourth} for a few examples.
In recent years, several studies (see, e.g., \cite{o2017energy, gao2019sbp, gao2019combining, gao2022nonuniform, o2022high}) have applied these methods to wave simulations on staggered grids, which is a popular choice of grid layout in seismic studies.

With the SBP operators, a common approach to impose boundary conditions is by appending penalty terms (i.e., weakly), sometimes referred to as the simultaneous-approximation-terms (SAT) technique; see \cite{carpenter1994time}.
One of the advantages of the SAT technique is its flexibility, since it enables the separation of concern between designing the discrete operators that approximate differential operators and addressing the the boundary conditions, which allows a modular procedure when approaching a discretization task and the reusability of designed operators.

In this study, we consider a particular simulation setting pertinent to the practical applications mentioned above, namely, a point source is placed near or on the boundary associated with the free surface boundary condition. 
In seismic studies for example, this mimics the survey setting where controlled
seismic sources%
\footnotemark
\footnotetext
{
In fact, this setting is also relevant to seismic studies involving natural (i.e., earthquake) sources, which happens deep beneath the surface.
This is because in commonly used imaging techniques such as reverse time migration or full waveform inversion, the signals or data misfits at the sensor locations, which are placed on or near the earth surface, are used as the source terms in a back propagation or adjoint simulation that follows the forward propagation driven by the earthquake source to eventually form the image or update of the subterranean media. 
In this second simulation, the source location (i.e., the sensor location) is on or near the surface.
}
(e.g., those generated by airgun in marine environment or explosive in land environment) are used to inject energy to the subterranean media that drive the wave propagation. 
Signals generated by these sources are collected at sensor locations to infer the subterranean media. 
Similar application settings appear in other applications such as site characterization and non-destructive testing.

In this setting, we demonstrate that the weakly imposed boundary condition can have severe violation of the free surface boundary conditions and lead to inaccuracies in the simulation results.
%
%
To address this issue, we propose to build the free surface boundary condition into the SBP operators and hence strongly impose the boundary condition.
Both acoustic and elastic waves are considered.
For the acoustic case, the proposed procedure is very simple, involving only resetting a few rows and columns in the SBP operators to impose the free surface boundary condition strongly.
In the elastic case, the procedure is more involved and requires special design of the grid layout and SBP operators that satisfy additional requirement, as revealed by discrete energy analysis.
Numerical examples are presented to demonstrate the advantages of the strong approach over the weak approach for both acoustic and elastic waves, including improved accuracy and relaxed time step restriction.

In the remainder, we start with the 1D case to explain the procedures for both approaches and compare their outcomes in section \ref{section_1D_case}. 
We then extend the study to the 2D acoustic case in section \ref{section_2D_acoustic_case} and to the 2D elastic case in section \ref{section_2D_elastic_case}. 
Finally, we conclude in section \ref{section_conclusions}.

\section{The 1D case}\label{section_1D_case}
\subsection{Background}\label{subsection_1D_background}
For clarity, we start our discussion with the 1D case, where both the acoustic and elastic wave equations reduce to the following system:
\begin{linenomath}
\begin{subequations}
\label{1D_wave_equation}
\begin{empheq}[left=\empheqlbrace]{alignat = 2}
\displaystyle \enskip \rho \frac{\partial v}{\partial t} \enskip &= \enskip \displaystyle \! \frac{\partial \sigma}{\partial x} \ + \ s^v \, ; \label{1D_wave_equation_a} \\
\displaystyle \enskip \frac{1}{\rho c^2} \frac{\partial \sigma}{\partial t} \enskip &= \enskip \displaystyle \! \frac{\partial v}{\partial x}  \ + \ s^\sigma \, , \label{1D_wave_equation_b}
\end{empheq}
\end{subequations}
\end{linenomath}
where $\rho$ and $c$ are given physical parameters, standing for density and wave-speed, respectively; $s^v$ and $s^\sigma$ are source terms that drive the wave propagation; $\sigma$ and $v$ are the sought solution variables.
When interpreted as the acoustic wave system, $\sigma$ represents the negative of pressure;
when interpreted as the elastic wave system, $\sigma$ represents the stress. 
In both cases, $v$ represents the velocity. 

In seismic survey involving artificially generated sources, 
$s^v$ and $s^\sigma$ are often abstracted as point sources. 
In this study, we consider the simple case where the point source location coincides with a grid point. 
The term $s^v$ that applies on velocity is often referred to as the directional source, which can be generated by, e.g., striking a sledgehammer on a metal plate on the ground; 
the term $s^\sigma$ that applies on pressure or stress is often referred to as the compressional source, which can be generated using, e.g., airgun or explosive. 
%
%
%

Supposing that wave system \eqref{1D_wave_equation} is defined over an interval $(x_L,x_R)$, its associated physical energy can be expressed as:
\begin{equation}
\label{1D_continuous_energy}
\mathscr E \ = \ \frac{1}{2} \int_{x_L}^{x_R} \rho v^2 dx \ + \ \frac{1}{2} \int_{x_L}^{x_R} \frac{1}{\rho c^2} \sigma^2 dx \, ,
\end{equation}
where the two terms represent the kinetic and potential energy in the system, respectively.
Omitting the source terms in \eqref{1D_wave_equation} for now,
taking the time derivative on both sides of \eqref{1D_continuous_energy} and substituting the equations from \eqref{1D_wave_equation}, we arrive at:
\begin{equation}
\label{1D_continuous_energy_derivative}
\frac{d \mathscr E}{d t} \ = \ - \ \sigma(x_L) \cdot v(x_L) \ + \ \sigma(x_R) \cdot v(x_R) \, .
\end{equation} 
In other words, in the absence of {\it source terms}, time derivative of the physical energy depends on {\it boundary data} only. If free surface boundary condition is associated with both boundaries, i.e., $\sigma(x_L) = \sigma(x_R) = 0$, we have $\frac{d \mathscr E}{d t} = 0$, i.e., the physical energy is conserved.

\begin{figure}[H]
\captionsetup{width=1\textwidth, font=small,labelfont=small}
\centering\includegraphics[scale=0.0675]{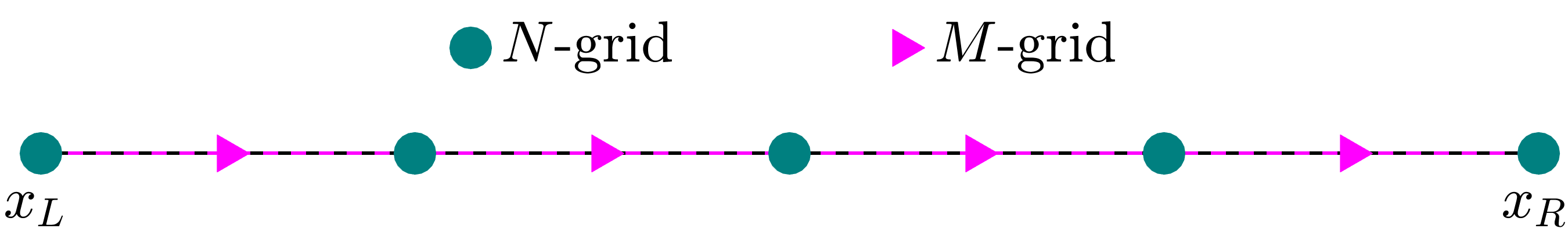}
\caption{Illustration of the grid layout underlying the discretization presented below.
}
\label{1D_grid_layout_extrapolating}
\end{figure}

The continuous wave system \eqref{1D_wave_equation} can be discretized in space using the SBP-SAT approach.
As an example, one such discretization on staggered grids has been presented in \cite{gao2019sbp}. 
The underlying grid layout is illustrated in Figure \ref{1D_grid_layout_extrapolating}, where the two sub-grids are referred to as the $N$-grid 
and $M$-grid, 
respectively.
For brevity, we only sketch the outline of the SBP-SAT discretization below. For a more detailed discussion, the readers are referred to \cite{gao2019sbp}.

Omitting the source terms and boundary condition, \eqref{1D_wave_equation} can be discretized in space using the SBP operators as follows:

\begin{linenomath}
\begin{subequations}
\label{1D_semi_discretized_extrapolating_wo_SATs}
\begin{empheq}[left=\empheqlbrace]{alignat = 2}
\displaystyle \enskip \mathcal A^M \boldsymbol \rho^M \frac{d V}{d t} 
\enskip &= \enskip \displaystyle 
\mathcal A^M \mathcal D^N \Sigma \, ; \label{1D_semi_discretized_extrapolating_wo_SATs_a} \\[0.25ex]
\displaystyle \enskip \mathcal A^N \boldsymbol \beta^N \frac{d \Sigma}{d t} 
\enskip &= \enskip \displaystyle 
\mathcal A^N \mathcal D^M V \, ,
\label{1D_semi_discretized_extrapolating_wo_SATs_b}
\end{empheq}
\end{subequations}
\end{linenomath}
where $V$ and $\Sigma$ are the discrete solution vectors, occupying the $M$-grid and $N$-grid, respectively, 
$\boldsymbol \rho^M$ and $\boldsymbol \beta^N$ are the discretization of physical parameters $\rho$ and $\nicefrac{1}{\rho c^2}$ on the $M$-grid and $N$-grid, respectively.
Matrices $\mathcal D^N$ and $\mathcal D^M$ are the difference operators; $\mathcal A^M$ and $\mathcal A^N$ are the norm matrices; together, they constitute the SBP operators.

The discrete energy associated with \eqref{1D_semi_discretized_extrapolating_wo_SATs} is defined as 
\begin{equation}
\label{1D_discrete_energy}
E \ = \ \frac{1}{2} V^T \left(\mathcal A^M \boldsymbol \rho^M \right) V 
  \ + \ \frac{1}{2} \Sigma^T \left(\mathcal A^N \boldsymbol \beta^N \right) \Sigma \, ,
\end{equation}
where the norm matrices $\mathcal A^M$ and $\mathcal A^N$ play the roles of the integrals from \eqref{1D_continuous_energy}.
Taking the time derivative on both sides of \eqref{1D_discrete_energy} and substituting the equations from \eqref{1D_semi_discretized_extrapolating_wo_SATs}, we arrive at:
\begin{equation}
\label{1D_discrete_energy_derivative}
\frac{d E}{d t} \ = \ \Sigma^T \left[ \mathcal A^N \mathcal D^M \ + \ \left( \mathcal A^M \mathcal D^N \right)^T \right] V \, .
\end{equation} 
To simplify the notation, we define
\begin{equation}
\label{1D_Q_definition}
Q \ = \ \mathcal A^N \mathcal D^M \ + \ \left( \mathcal A^M \mathcal D^N \right)^T .
\end{equation}
When designing the SBP operators, $Q$ is asked to satisfy the following property:
\begin{equation}
\label{1D_Q_property}
Q \ = \ - \ \mathcal E^L \left(\mathcal P^L\right)^T + \ \mathcal E^R \left(\mathcal P^R\right)^T ,
\end{equation}
where $\mathcal E^L$ and $\mathcal E^R$ are the canonical basis vectors whose first and last entries are 1, respectively, and therefore, select the first and last entries of $\Sigma$ when applied on it;
$\mathcal P^L$ and $\mathcal P^R$ are projection operators that can be applied on vector $V$ to provide approximations to $v(x_L)$ and $v(x_R)$, respectively.

We now can rewrite \eqref{1D_discrete_energy_derivative} as
\begin{equation}
\label{1D_discrete_energy_derivative_boundary}
\frac{d E}{d t} \ = \ - \ \left( \mathcal E^L \right)^T \! \Sigma \cdot \left( \mathcal P^L \right)^T \! V 
                    \ + \ \left( \mathcal E^R \right)^T \! \Sigma \cdot \left( \mathcal P^R \right)^T \! V \, .
\end{equation} 
Noticing 
that $\left( \mathcal E^L \right)^T \! \Sigma$ and $\left( \mathcal P^L \right)^T \! V$ are approximations to $\sigma(x_L)$ and $v(x_L)$ and 
that $\left( \mathcal E^R \right)^T \! \Sigma$ and $\left( \mathcal P^R \right)^T \! V$ are approximations to $\sigma(x_R)$ and $v(x_R)$, \eqref{1D_discrete_energy_derivative_boundary} is the discrete equivalent to \eqref{1D_continuous_energy_derivative}.

Moreover, the free surface boundary condition 
can be imposed weakly by appending two trailing terms i.e, SATs, to \eqref{1D_semi_discretized_extrapolating_wo_SATs_a}, leading to the following system:
\begin{linenomath}
\begin{subequations}
\label{1D_semi_discretized_extrapolating_with_SATs}
\begin{empheq}[left=\empheqlbrace]{alignat = 2}
\displaystyle \enskip \mathcal A^M \boldsymbol \rho^M \frac{d V}{d t} 
\enskip &= \enskip \displaystyle 
\mathcal A^M \mathcal D^N \Sigma 
\ + \ \mathcal P^L \left[ \left( \mathcal E^L \right)^T \! \Sigma - \bf 0 \right] 
\ - \ \mathcal P^R \left[ \left( \mathcal E^R \right)^T \! \Sigma - \bf 0 \right] \, ; \label{1D_semi_discretized_extrapolating_with_SATs_a} \\[0.25ex]
\displaystyle \enskip \mathcal A^N \boldsymbol \beta^N \frac{d \Sigma}{d t} 
\enskip &= \enskip \displaystyle 
\mathcal A^N \mathcal D^M V \, .
\label{1D_semi_discretized_extrapolating_with_SATs_b}
\end{empheq}
\end{subequations}
\end{linenomath}
One can easily verify that $\tfrac{d E}{d t} = 0$ for this system, i.e., the discrete energy is conserved, just as the continuous energy $\mathscr E$ is conserved under the free surface boundary condition; see \eqref{1D_continuous_energy_derivative}.

\begin{remark}\label{remark_flexibility}
For the SBP-SAT approach described above, the design of SBP operators and the design of SATs can be detached from each other. 
When designing the SBP operators, one only needs to be careful with leaving some suitable residual terms on the boundary; see \eqref{1D_discrete_energy_derivative_boundary}.
The boundary conditions are addressed later via the SATs and are not part of the consideration when designing SBP operators. 
For a different type of boundary condition, one only needs to change the SATs accordingly and can reuse the same SBP operators.
This flexibility offers great advantage in practice. 
%
%
\end{remark}

Below we give a concrete example of the SBP operators (originally presented in \cite{gao2019sbp}) in \eqref{SBP_matrices_1D_extrapolating} for the grid points near the left boundary (see Figure \ref{1D_grid_layout_extrapolating}).
These operators correspond to the case of unit grid spacing, i.e., $\Delta x=1$. 
For general cases, $\mathcal D^N$ and $\mathcal D^M$ need to be scaled by $\tfrac{1}{\Delta x}$ whereas $\mathcal A^N$ and $\mathcal A^M$ need to be scaled by $\Delta x$.
The bottom right corners of $\mathcal A^M$ and $\mathcal A^N$ are the mirror opposite of their top left counterparts; the bottom right corners of $\mathcal D^N$ and $\mathcal D^M$ are the negative mirror opposite of their top left counterparts.
They are omitted here to conserve space.%
\footnotemark
\footnotetext
{
We note here that the set of operators presented in \eqref{SBP_matrices_1D_extrapolating} is not the unique choice. 
Other choices exist; see \cite{o2017energy} for another example. 
However, numerical experiments reveal that the operators presented in \cite{o2017energy} can lead to severe restriction on time step size. 
Specifically, using the same 1D test with homogeneous medium that leads to \eqref{1D_extrapolating_CFL_limit_weak} and \eqref{1D_extrapolating_CFL_limit_strong}, it is revealed that the maximally allowed time step size allowed by the operators presented in \cite{o2017energy} is less than half ($\approx 43.3\%$) of that associated with the interior stencil; cf. \eqref{1D_extrapolating_CFL_limit_weak} and \eqref{1D_extrapolating_CFL_limit_strong}.
For this reason, we omit this alternative choice in this study.
}

\begin{subequations}
\label{SBP_matrices_1D_extrapolating}
\footnotesize
\begin{align}
\label{SBP_matrices_1D_extrapolating_A}
&
\mathcal A^M = \left[ 
\arraycolsep=3.2pt
\begin{array}{r r r r r r r r r r r r r r r r r} 
\nicefrac{13}{12} &                 & \multicolumn{1}{r|}{} \\
		          & \nicefrac{7}{8} & \multicolumn{1}{r|}{} \\
			      &                 & \multicolumn{1}{r|}{ \nicefrac{25}{24} }  \\ \cline{1-3} \\[-1.75ex]
				  &  			    &    			& 		   1 &			  & \\
				  & 			    &    			& 			 & 			1 & \\[-0.75ex]
				  &  			    &    			& 			 &   		  & \ddots
\end{array} 
\right],
&
&
\mathcal A^N = \left[ 
\arraycolsep=3.2pt
\begin{array}{r r r r r r r r r r r r r r r r r} 
\nicefrac{7}{18} & 				   & 		   		& \multicolumn{1}{r|}{} \\
		         & \nicefrac{9}{8} &     			& \multicolumn{1}{r|}{} \\
  			     &	               & 1  			& \multicolumn{1}{r|}{} \\
			     &    			   &    			& \multicolumn{1}{r|}{ \nicefrac{71}{72} }    \\ \cline{1-4} \\[-1.75ex]
				 &   			   &    			&    			&  1  			&    		  \\
				 &    			   &  				&    			&    			&  1  		& \\[-0.75ex]
				 &  			   &    			&   		    &    			&      		& \ddots
\end{array} 
\right];
&
\allowdisplaybreaks
\\
\label{SBP_matrices_1D_extrapolating_D}
&
\mathcal D^N = \left[ 
\arraycolsep=3.2pt
\begin{array}{r r r r r r r r r r r r r r r r r} 
-\nicefrac{79}{78} &  \nicefrac{27}{26} & -\nicefrac{1}{26}  &  \nicefrac{1}{78}  & \multicolumn{1}{r|}{0}      \\
 \nicefrac{2}{21}  & -\nicefrac{9}{7}   &  \nicefrac{9}{7}   & -\nicefrac{2}{21}  & \multicolumn{1}{r|}{0}      \\
 \nicefrac{1}{75}  &  0			        & -\nicefrac{27}{25} &  \nicefrac{83}{75} & \multicolumn{1}{r|}{ -\nicefrac{1}{25} } & \\ \cline{1-5}
				   &   					&  \nicefrac{1}{24}  & -\nicefrac{9}{8}   & \nicefrac{9}{8}   						 & -\nicefrac{1}{24} 		  & \\
				   &     				&    				 &  \nicefrac{1}{24}  & -\nicefrac{9}{8} 	  					 &  \nicefrac{9}{8}  		  & -\nicefrac{1}{24} \\[-0.5ex]
				   & 					& 					 & 					  & \multicolumn{1}{l}{\ddots} 				 & \multicolumn{1}{l}{\ddots} & \multicolumn{1}{l}{\ddots} & \multicolumn{1}{l}{\ddots}  
\end{array} 
\right], 
&
&
\mathcal D^M = \left[ 
\arraycolsep=3.2pt
\begin{array}{r r r r r r r r r r r r r r r r r} 
-2       		  &  3    			  & -1    			   &  0    			    & \multicolumn{1}{r|}{0}   \\
-1       		  &  1    			  &  0    			   &  0    			    & \multicolumn{1}{r|}{0}   \\
 \nicefrac{1}{24} & -\nicefrac{9}{8}  &  \nicefrac{9}{8}   & -\nicefrac{1}{24}  & \multicolumn{1}{r|}{0} & \\
-\nicefrac{1}{71} &  \nicefrac{6}{71} & -\nicefrac{83}{71} &  \nicefrac{81}{71} & \multicolumn{1}{r|}{ -\nicefrac{3}{71} } & \\ \cline{1-5}
				  &    				  &  \nicefrac{1}{24}  & -\nicefrac{9}{8}   &  \nicefrac{9}{8}   		 			   & -\nicefrac{1}{24} 		    & \\
				  &     			  &    				   &  \nicefrac{1}{24}  & -\nicefrac{9}{8}   		 			   &  \nicefrac{9}{8}  		    & -\nicefrac{1}{24} \\[-0.5ex]
				  & 				  & 				   & 					& \multicolumn{1}{l}{\ddots} 			   & \multicolumn{1}{l}{\ddots} & \multicolumn{1}{l}{\ddots} & \multicolumn{1}{l}{\ddots}  
\end{array} 
\right].
&
\end{align}
\end{subequations}
With the above matrices, $Q$ takes the following form:
\begin{equation}
\label{Q_1D_matrix_weak}
Q =
\footnotesize
\left[ 
\arraycolsep=3.2pt
\begin{array}{r r r r r r r r r r r r r r r r}
-\nicefrac{15}{8} & \nicefrac{5}{4} & -\nicefrac{3}{8} \\
&    &    &    &    &    &    &    &    &    &    &    \\
&    &    &    &    &    &    &    &    &    &    &    \\
&    &    &    &    &    &    &    &    &    &    &    \\
&    &    &    &    &    &    &    &    &    &    & \nicefrac{3}{8} & -\nicefrac{5}{4} & \nicefrac{15}{8}
\end{array}
\right],
\end{equation}
i.e., only the first and last rows are nonzero and correspond to $-\left(\mathcal P^L\right)^T$ and $\left(\mathcal P^R\right)^T$, respectively.

\subsection{Numerical experiments}\label{subsection_1D_numerical_experiments}
In the following, we conduct a series of numerical experiments to illustrate the behavior of the above discretization in the presence of source terms near the boundary.
For clarity, homogeneous medium with unit density and wave-speed (i.e., $\rho=1~\text{kg}/\text{m}^3$ and $c=1~\text{m}/\text{s}$) is considered.
The free surface boundary condition is associated with both boundaries (i.e., $\sigma=0$ at $x_L$ and $x_R$). 

A compressional point source is used to drive the wave propagation. 
Temporal profile of the point source is specified as the Ricker wavelet with central frequency 5~Hz and time delay 0.25~s; see \cite[p.~684]{gao2019combining} for more detail about the source specification.
The maximal frequency in the source content is considered as 12.5~Hz, which leads to a minimal wavelength of 0.08~m.
The length of the interval $(x_L,x_R)$ is specified as 20 of such minimal wavelengths (i.e., $1.6$~m).
%
%

The stencil used for these experiments is the standard fourth-order staggered grid stencil in the interior (i.e., 
$
\nicefrac{ [\nicefrac{1}{24}, \enskip -\nicefrac{9}{8}, \enskip \nicefrac{9}{8}, \enskip -\nicefrac{1}{24}] }{ \Delta x }
$;
see, e.g., \cite{levander1988fourth,fornberg1999spatial})
with adaptations near the boundaries, as illustrated in \eqref{SBP_matrices_1D_extrapolating}, to satisfy the SBP property.
%
The staggered leapfrog scheme is used for time integration; see, e.g., \cite{ghrist2000staggered} for more detail.
%
%
Denoting the Courant number as $C = \nicefrac{\Delta t}{\Delta x}$ (the unit wave-speed $c$ has been omitted), the CFL restriction for the aforementioned interior stencil and time integration scheme on unbounded domain stipulates 
\begin{equation}
\label{Interior_CFL_limit}
C \leq C_{\max} = \nicefrac{6}{7}.
\end{equation}
Interested readers may consult \cite{levander1988fourth} or \ref{section_Courant_number_6_over_7} for more detail. 

\begin{figure}[H]
\captionsetup{width=1\textwidth, font=small,labelfont=small}
\centering\includegraphics[scale=0.0725]{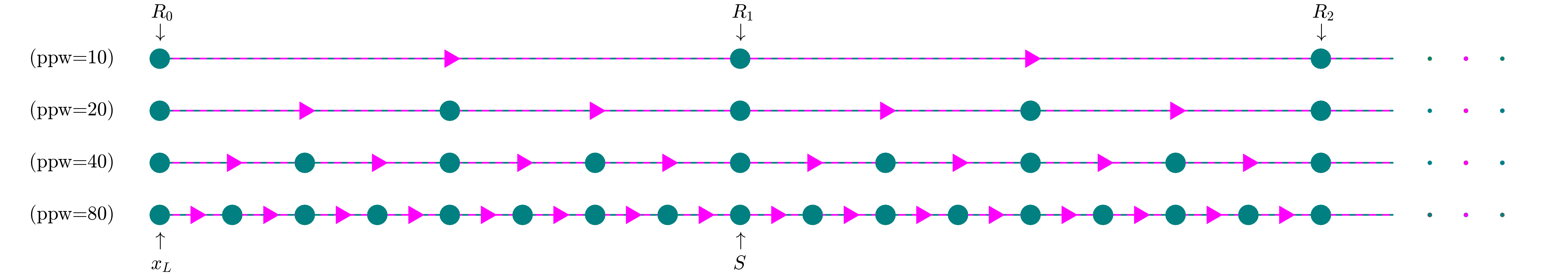}
\caption{Grid layouts and source receiver locations for the 1D experiments.
%
%
The four grid lines, from top to the bottom, correspond to the left ends of the staggered grids with 10, 20, 40, and 80 ppw, respectively. 
Source and receiver locations are placed at fixed physical locations for all four simulations, and indicated by the arrows and letters ($S$ for source; $R$ for receiver). 
} 
\label{Figure_grid_layout_source_receiver}
\end{figure}

Four experiments have been conducted, which use 10, 20, 40, and 80 grid points per minimal wavelength (ppw), 
respectively. 
The grid layouts and source receiver locations in these experiments are illustrated in Figure \ref{Figure_grid_layout_source_receiver}.
Specifically, the point sources in these experiments are placed at a fixed physical location, which are 1, 2, 4, and 8 grid point(s) away from the left boundary $x_L$, respectively.
The time step length is specified as 2.5e-4~s for all four simulations, and the simulation time is specified as 6~s, which amounts to 24000 time steps.

\begin{figure}[H]
\captionsetup{width=1\textwidth, font=small,labelfont=small}
\centering
\begin{subfigure}[b]{1\textwidth}
\captionsetup{width=1\textwidth, font=small,labelfont=small}
\centering\includegraphics[scale=0.15]{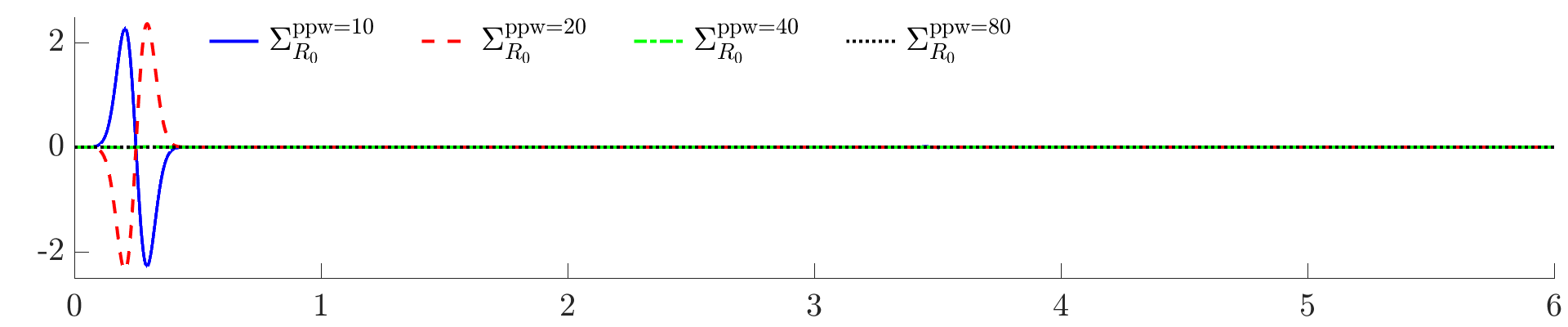}
\caption{Time histories of $\Sigma$ at $R_0$.} 
\label{P_0_6_s_extrapolating}
\end{subfigure}\hfill
\\[2ex]
\begin{subfigure}[b]{1\textwidth}
\captionsetup{width=1\textwidth, font=small,labelfont=small}
\centering\includegraphics[scale=0.15]{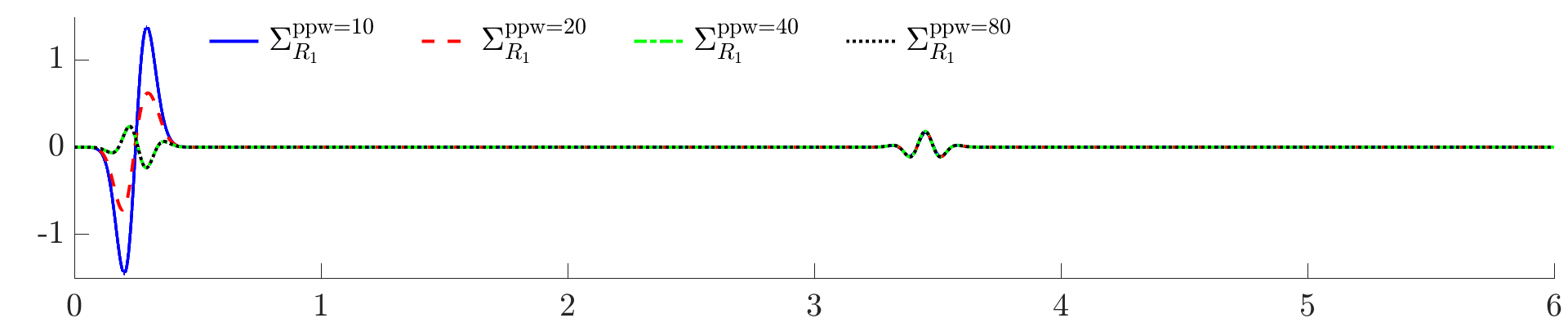}
\caption{Time histories of $\Sigma$ at $R_1$.}
\label{P_1_6_s_extrapolating}
\end{subfigure}\hfill
\\[2ex]
\begin{subfigure}[b]{1\textwidth}
\captionsetup{width=1\textwidth, font=small,labelfont=small}
\centering\includegraphics[scale=0.15]{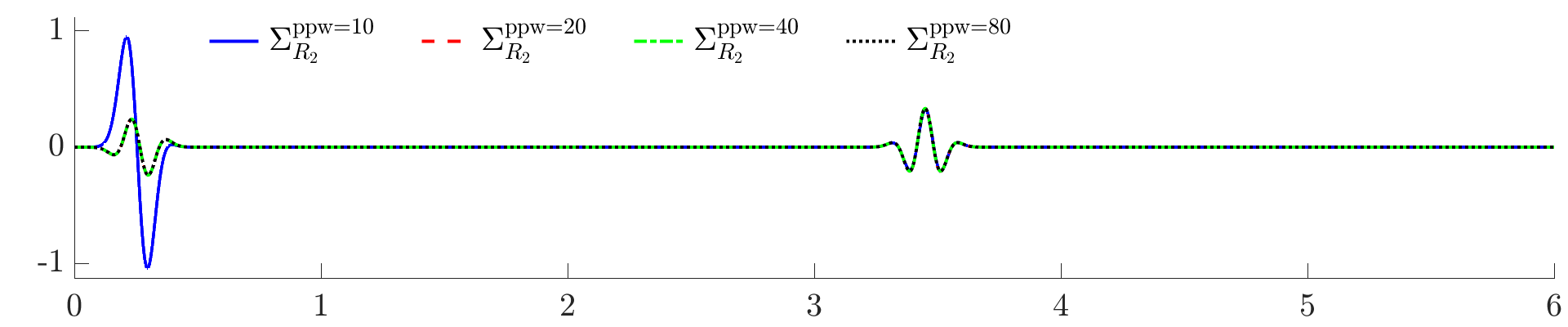}
\caption{Time histories of $\Sigma$ at $R_2$.}
\label{P_2_6_s_extrapolating}
\end{subfigure}\hfill
\caption{Time histories of the solution variable $\Sigma$ at the three receiver locations depicted in Figure \ref{Figure_grid_layout_source_receiver} for the four experiments where the distances, in terms of number of grid points, between the source and the free surface boundary are 1, 2, 4, and 8, respectively. 
}
\label{P_extrapolating_weak}
\end{figure}

Time histories of the solution variable $\Sigma$ are recorded at three fixed physical locations ($R_0$, $R_1$, and $R_2$ in Figure \ref{Figure_grid_layout_source_receiver}) and displayed in Figure \ref{P_extrapolating_weak}. 
We note here that if the free surface boundary condition is strictly satisfied, $\Sigma$ at $R_0$ (i.e., the left boundary) should remain zero for the entire simulation duration for all four experiments. 
From Figure \ref{P_extrapolating_weak}, we observe that when the point source is placed too close to the free surface (in terms of the number of grid points; see Figure \ref{Figure_grid_layout_source_receiver}), which corresponds to the smaller ppw cases, the weakly imposed free surface boundary condition is severely violated, leading to significant inaccuracy in the recorded time histories.
\footnotemark
\footnotetext
{
Two additional experiments are presented in the {\it Supplementary Material} (\ref*{supp_section_finer_grid} and \ref*{supp_section_interior_source}) to demonstrate that such violation is not due to the lack of grid resolution and that such violation subsides when the source is placed far away from the free surface.
}

Moreover, numerical experiments reveal that the CFL restriction for the above discretization setting (encompassing the difference operators, weakly imposed boundary terms, and time integration scheme) is
\begin{equation}
\label{1D_extrapolating_CFL_limit_weak}
C \leq C_{\max} \approx 0.6355,
\end{equation}
which is approximately 74.1\% of that associated with the interior stencil on an unbounded domain; see \eqref{Interior_CFL_limit}.
%

\subsection{Strong imposition of the free surface boundary condition}\label{subsection_1D_strong_imposition}

To address the severe violation of the free surface boundary condition illustrated above, we propose to impose the free surface boundary condition strongly by incorporating it into the difference operators (cf. Remark \ref{remark_flexibility}).
While this can be achieved by starting anew the process of designing the SBP operators, it can also be achieved by adapting the existing operators originally designed with weakly imposed boundary conditions in mind. 
We find this later approach easier to comprehend and illuminate the connections and contrasts between the weak and strong approaches better and will present it in the following.

Taking the SBP operators illustrated in \eqref{SBP_matrices_1D_extrapolating} and its associated grid layout (see Figure \ref{1D_grid_layout_extrapolating}) as an example, 
the first entry of the solution vector $\Sigma$, denoted as $\Sigma(1)$ and corresponding to the left boundary $x_L$, is updated by the first row of $\mathcal D^M$; see \eqref{1D_semi_discretized_extrapolating_wo_SATs}.
From the free surface boundary condition at $x_L$, i.e., $\sigma(x_L) = 0$, we know that $\Sigma(1)$ is supposed to remain zero during the simulation. 
We can simply set the first row of $\mathcal D^M$ to zero to achieve this. 
On the other hand, since $\Sigma(1)$ will be multiplied by the first column of $\mathcal D^N$ to update $V$, we can also set the first column of $\mathcal D^N$ to zero to achieve the same effect.%

In doing so, information of the free surface boundary condition is included in the simulation, and therefore, the SATs in \eqref{1D_semi_discretized_extrapolating_with_SATs} are no longer needed. 
For all other grid points, their stencils remain (effectively) unchanged and the accuracy of their approximations inherits from the original design of $\mathcal D^M$ and $\mathcal D^N$.
Finally, if we set both the first row of $\mathcal D^M$ and the first column of $\mathcal D^N$ to zero, it is easy to verify that the first row of $Q$ is also zero; see \eqref{1D_Q_definition} and \eqref{Q_1D_matrix_weak}. 
Performing the corresponding operations for the right boundary $x_R$ (i.e., setting the last row of $\mathcal D^M$ and the last column of $\mathcal D^N$ to zero), we have that 
\begin{equation}
\label{Q_1D_property_strong}
Q = \boldsymbol 0.
\end{equation}
Recalling the discrete energy analysis result from \eqref{1D_discrete_energy_derivative}, we have that $\frac{d E}{d t} = 0$ because of \eqref{Q_1D_property_strong}, i.e., the semi-discretized system \eqref{1D_semi_discretized_extrapolating_wo_SATs}, 
%
%
with operators modified as described above,
is also energy-conserving.

\begin{figure}[H]
\captionsetup{width=1\textwidth, font=small,labelfont=small}
\centering
\begin{subfigure}[b]{1\textwidth}
\captionsetup{width=1\textwidth, font=small,labelfont=small}
\centering\includegraphics[scale=0.15]{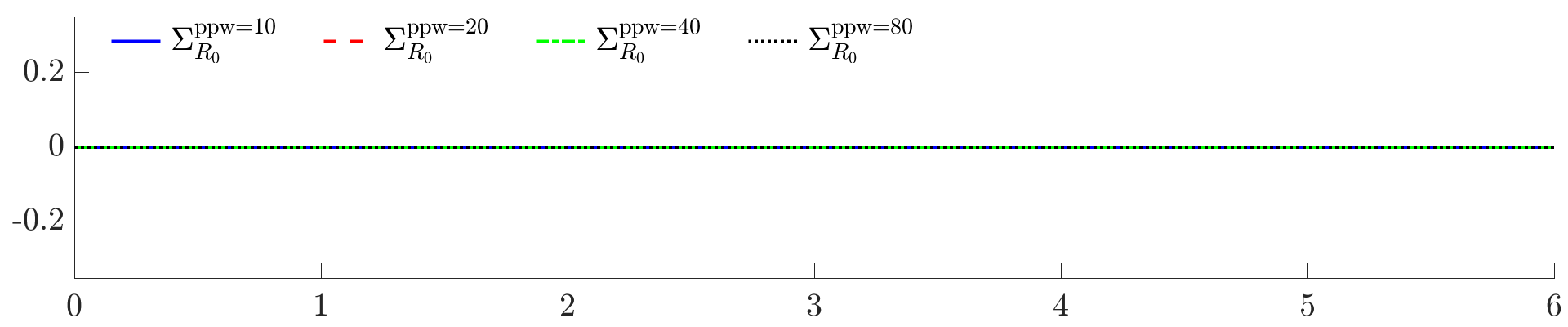}
\caption{Time history of $\Sigma$ at $R_0$.}
\label{P_0_6_s_extrapolating_strong}
\end{subfigure}\hfill
\\[2ex]
\begin{subfigure}[b]{1\textwidth}
\captionsetup{width=1\textwidth, font=small,labelfont=small}
\centering\includegraphics[scale=0.15]{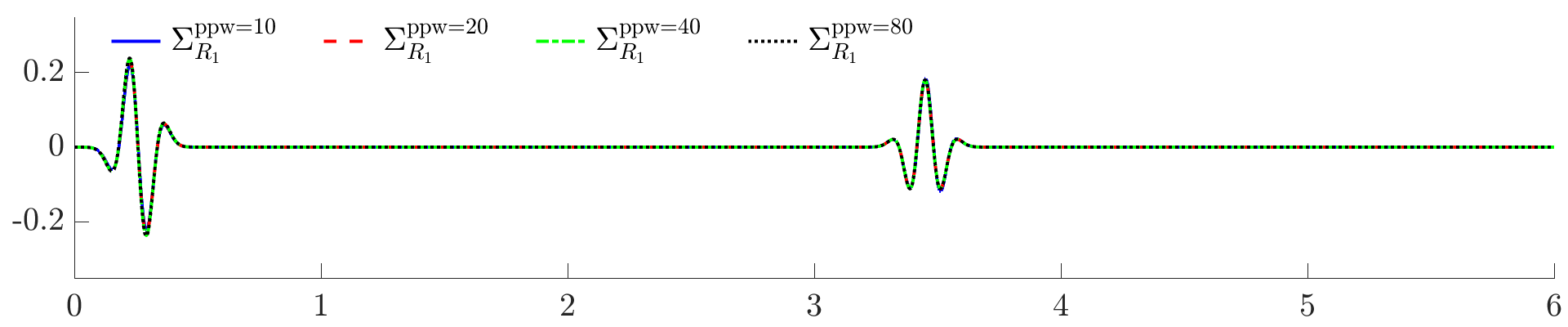}
\caption{Time history of $\Sigma$ at $R_1$.}
\label{P_1_6_s_extrapolating_strong}
\end{subfigure}\hfill
\\[2ex]
\begin{subfigure}[b]{1\textwidth}
\captionsetup{width=1\textwidth, font=small,labelfont=small}
\centering\includegraphics[scale=0.15]{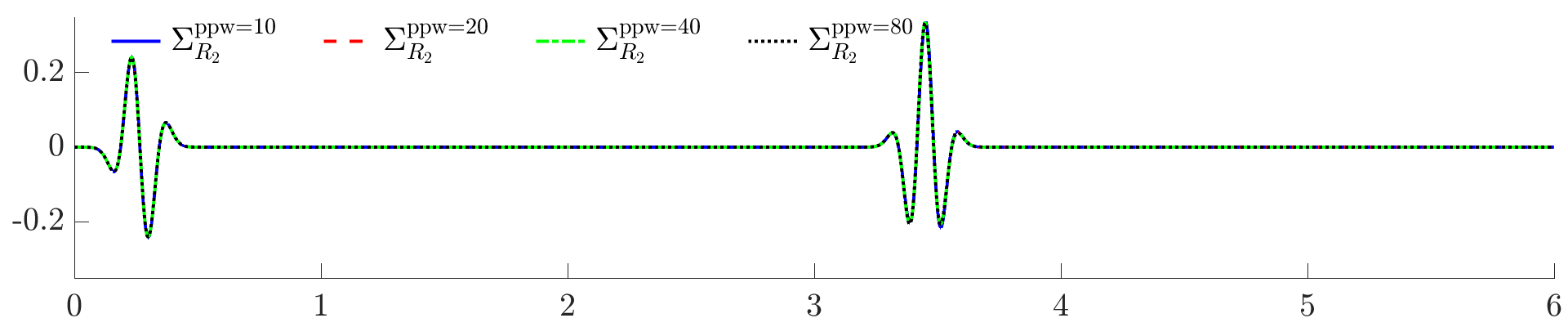}
\caption{Time history of $\Sigma$ at $R_2$.}
\label{P_2_6_s_extrapolating_strong}
\end{subfigure}\hfill
\caption{Time histories of the solution variable $\Sigma$ at the three receiver locations depicted in Figure \ref{Figure_grid_layout_source_receiver}.
The experimental setup is identical to that for Figure \ref{P_extrapolating_weak}.
The free surface boundary condition is imposed strongly.
}
\label{P_extrapolating_strong}
\end{figure}

Conducting the same experiments as those in section \ref{subsection_1D_numerical_experiments} with the free surface boundary condition imposed strongly, the recorded signals are displayed in Figure \ref{P_extrapolating_strong}, from where we observe that the violation of the free surface boundary condition disappears and that results from the different simulations agree well with each other (cf. Figure \ref{P_extrapolating_weak}).
%

Moreover, numerical experiments reveal that the CFL restriction for the above discretization with strongly imposed boundary conditions is 
\begin{equation}
\label{1D_extrapolating_CFL_limit_strong}
C \leq C_{\max} = \nicefrac{6}{7},
\end{equation}
i.e., the same as that allowed by the interior stencil; see \eqref{Interior_CFL_limit} and compare to \eqref{1D_extrapolating_CFL_limit_weak}. 
%

For a cross validation of the simulation results from different discretizations, including the two discretizations described above (i.e., the weak and strong cases, respectively) and a finite element discretization, the readers are referred to the {\it Supplementary Material} \ref*{supp_cross_comparison}. 
For the spectral radii information associated with the semi-discretized systems, which provide confirmation to the above claims on CFL restrictions, the readers are referred to the {\it Supplementary Material} \ref*{supp_spectral_radii}.

\section{The 2D acoustic case}\label{section_2D_acoustic_case}
Extension of the above procedure to the 2D acoustic wave equation is straightforward.
We can simply set the rows of the 2D equivalence of $\mathcal D^M$ that correspond to grid points on the free surface to zero and set the columns of the 2D equivalence of $\mathcal D^N$ that correspond to grid points on the free surface to zero to strongly impose the free surface boundary condition, provided that the pressure grid points are placed on the free surface.

The acoustic equation and semi-discretized systems are included in the {\it Supplementary Material} \ref*{supp_acoustic_background} and omitted here to conserve space. A numerical example is presented below to demonstrate the effectiveness of the strongly imposed free surface boundary condition, in comparison to the weak approach.
Specifically, homogeneous medium with unit density and wave-speed is considered, as in the 1D case. The source specification is also the same as in the 1D case, which leads to a minimal wavelength of 0.08 m. 
The simulation domain is 6 by 6 of such minimal wavelengths (i.e., 0.48 m on both directions).
The free surface boundary condition is associated with all boundaries.  

\begin{figure}[H]
\captionsetup{width=1\textwidth, font=small,labelfont=small}
\centering\includegraphics[scale=0.0825]{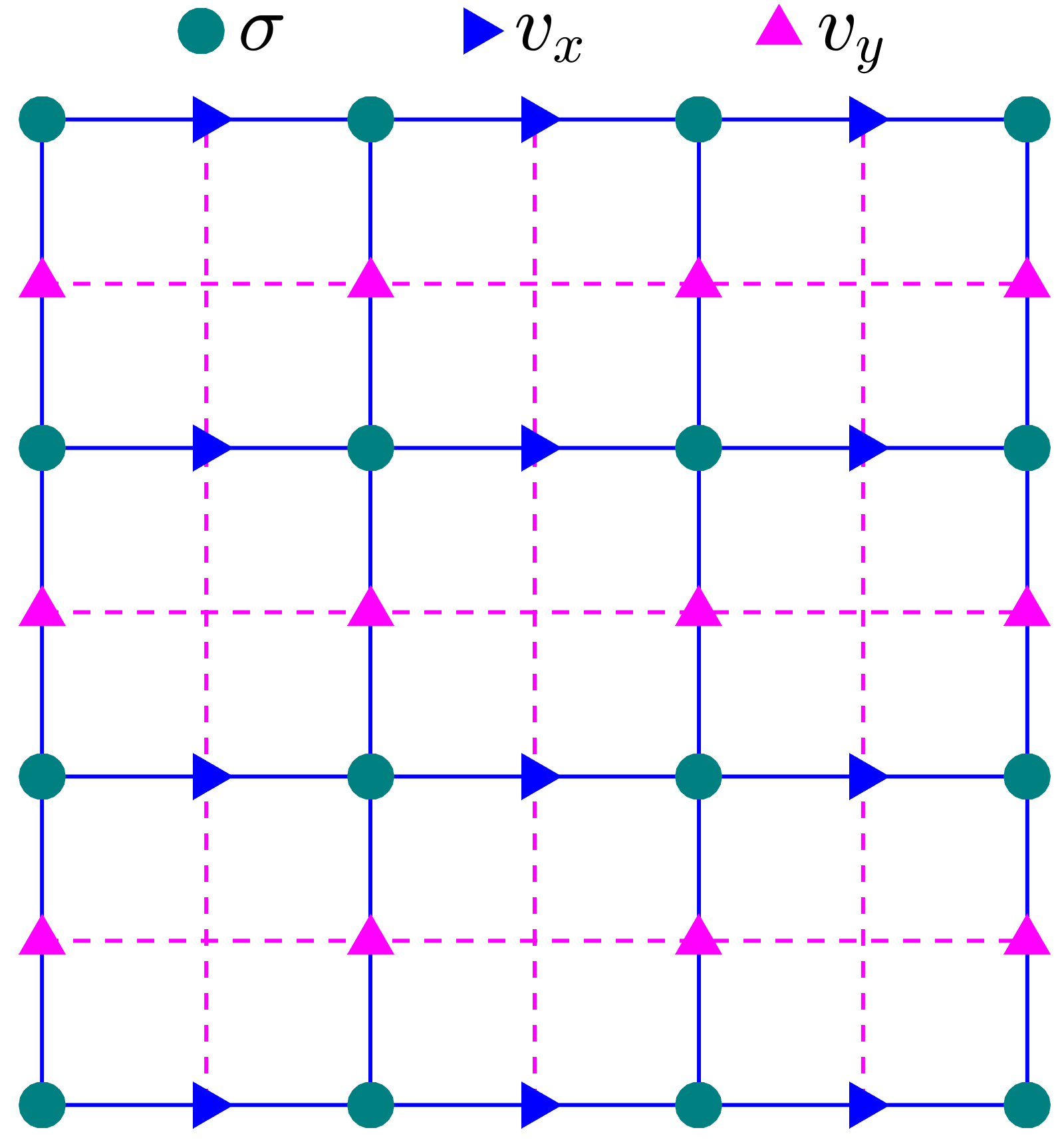}
\caption{Grid layout for the 2D acoustic wave equation.} 
\label{Figure_grid_layout_acoustic_full_FSBC}
\end{figure}

Three experiments are conducted for the strong and weak case each, which use 10, 30, and 50 grid points per minimal wavelength (ppw). The grid layout is illustrated in Figure \ref{Figure_grid_layout_acoustic_full_FSBC}.
The point source and receivers are placed at fixed physical locations at one third and two thirds of the simulation domain horizontally (i.e., 0.16 m and 0.32 m).
The source is placed at one grid point below the surface for the case where ppw is 10.
Three receivers are placed at the surface, at 3 and 6 grid points below the surface for the case where ppw is 10.
The time step length is specified as 5e-4 for all simulations, and the simulation time is specified as 3 s.

\begin{figure}[H]
\captionsetup{width=1\textwidth, font=small,labelfont=small}
\centering
\begin{subfigure}[b]{1\textwidth}
\captionsetup{width=1\textwidth, font=small,labelfont=small}
\centering\includegraphics[scale=0.15]{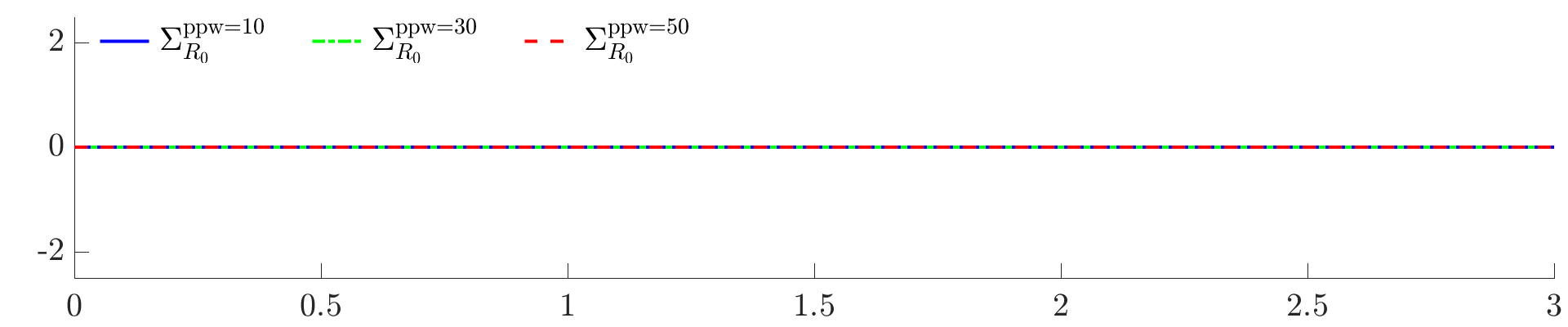}
\caption{Time history of $\Sigma$ at the surface.}
\label{Acoustic_full_FSBC_P_0_strong}
\end{subfigure}\hfill
\\[2ex]
\begin{subfigure}[b]{1\textwidth}
\captionsetup{width=1\textwidth, font=small,labelfont=small}
\centering\includegraphics[scale=0.15]{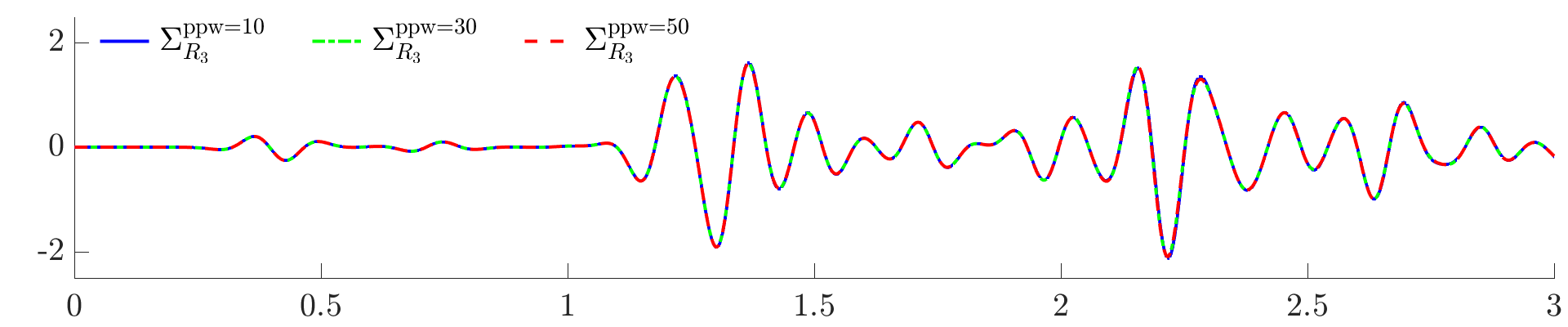}
\caption{Time history of $\Sigma$ at 3 grid points below the surface for the case ppw = 10.}
\label{Acoustic_full_FSBC_P_1_strong}
\end{subfigure}\hfill
\\[2ex]
\begin{subfigure}[b]{1\textwidth}
\captionsetup{width=1\textwidth, font=small,labelfont=small}
\centering\includegraphics[scale=0.15]{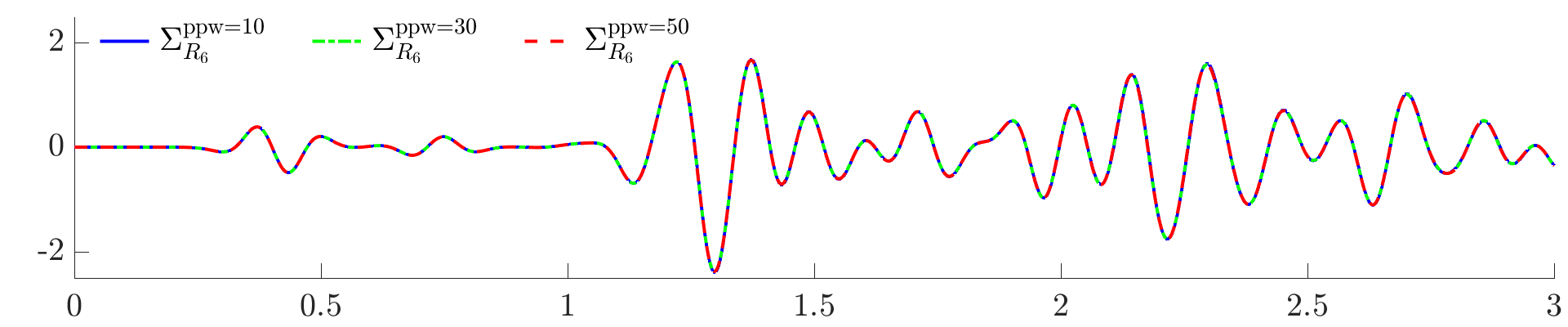}
\caption{Time history of $\Sigma$ at 6 grid points below the surface for the case ppw = 10.}
\label{Acoustic_full_FSBC_P_2_strong}
\end{subfigure}\hfill
\caption{Time histories of the 2D acoustic experiments. The free surface boundary condition is imposed strongly.
}
\label{Acoustic_full_FSBC_strong}
\end{figure}

\begin{figure}[H]
\captionsetup{width=1\textwidth, font=small,labelfont=small}
\centering
\begin{subfigure}[b]{1\textwidth}
\captionsetup{width=1\textwidth, font=small,labelfont=small}
\centering\includegraphics[scale=0.15]{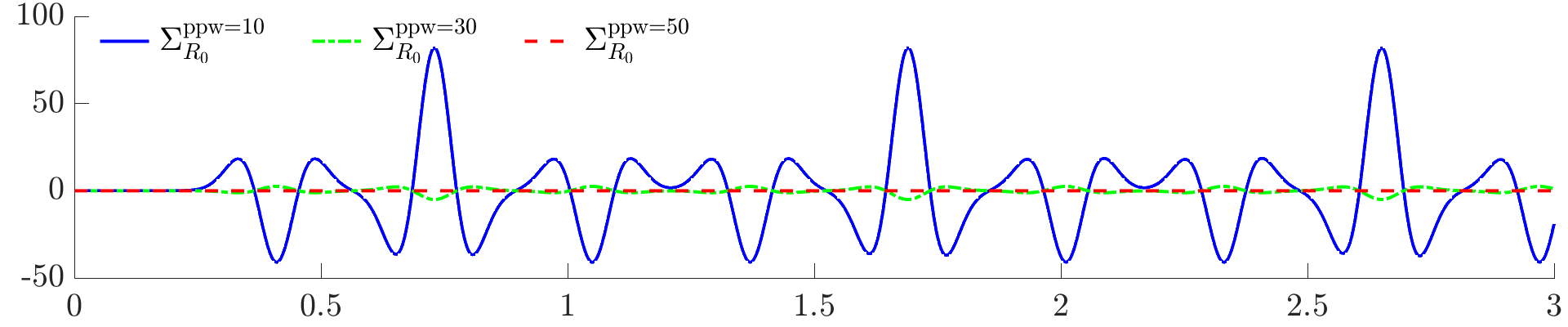}
\caption{Time history of $\Sigma$ at the surface.}
\label{Acoustic_full_FSBC_P_0_weak}
\end{subfigure}\hfill
\\[2ex]
\begin{subfigure}[b]{1\textwidth}
\captionsetup{width=1\textwidth, font=small,labelfont=small}
\centering\includegraphics[scale=0.15]{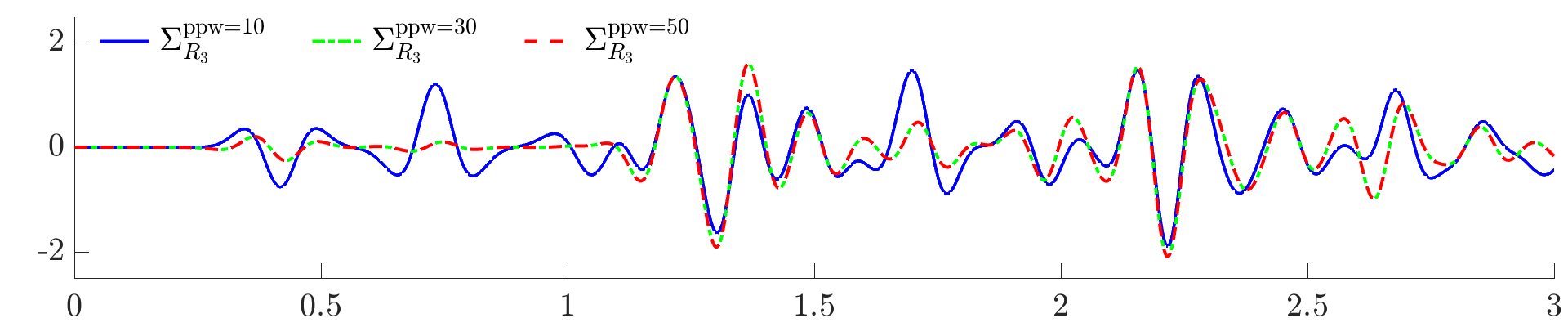}
\caption{Time history of $\Sigma$ at 3 grid points below the surface for the case ppw = 10.}
\label{Acoustic_full_FSBC_P_1_weak}
\end{subfigure}\hfill
\\[2ex]
\begin{subfigure}[b]{1\textwidth}
\captionsetup{width=1\textwidth, font=small,labelfont=small}
\centering\includegraphics[scale=0.15]{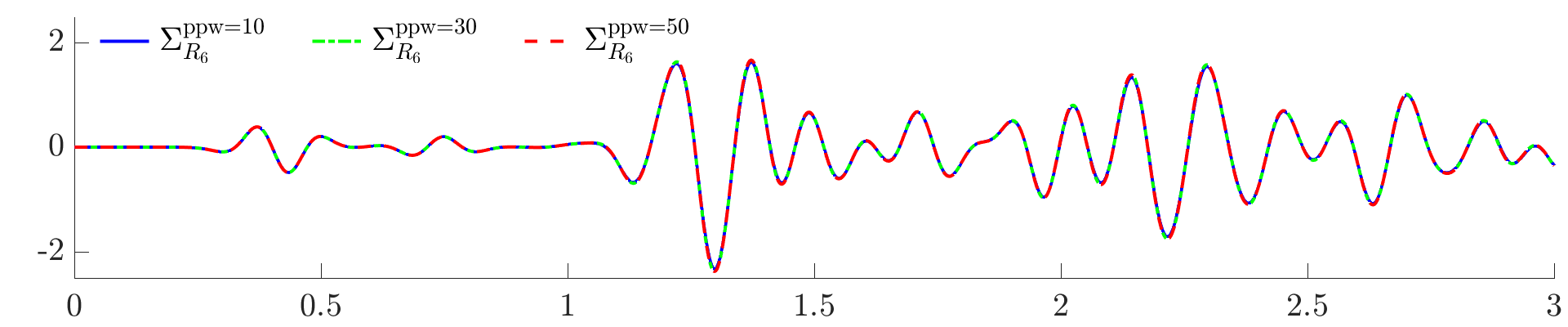}
\caption{Time history of $\Sigma$ at 6 grid points below the surface for the case ppw = 10.}
\label{Acoustic_full_FSBC_P_2_weak}
\end{subfigure}\hfill
\caption{Time histories of the 2D acoustic experiments. The free surface boundary condition is imposed weakly.
}
\label{Acoustic_full_FSBC_weak}
\end{figure}

Time histories of the solution variable $\Sigma$ at the three receiver locations are presented in Figures \ref{Acoustic_full_FSBC_strong} and \ref{Acoustic_full_FSBC_weak}, for the strong and weak case, respectively. 
As a reminder, variable $\sigma$ and its discretization $\Sigma$ represents the negative of pressure for the acoustic case. 
From these figures, we observe that the strongly imposed boundary condition delivers satisfactory results in all three simulations. 
We note here that because of the free surface boundary condition, $\Sigma$ is supposed to be zero at the surface; see Figure \ref{Acoustic_full_FSBC_P_0_strong}.
On the other hand, when the point source is placed too close to the surface (in terms of the number of grid points in between), the weakly imposed free surface boundary condition is severely violated, leading to inaccuracies in the simulated results. 
Such inaccuracy diminishes as the receiver depth increases.
Time histories of the solution variable $V_x$ and $V_y$ are included in the {\it Supplementary Material} \ref*{supp_additional_figures_2D_acoustic}.

Moreover, numerical experiments reveal that the CFL restriction for the strong and weak cases are $\nicefrac{6}{7}$ and 0.6355, respectively, the same as in the 1D case. 
In other words, for the strong case, there is no additional penalty on the time step size compared to that from the interior stencil; see \eqref{Interior_CFL_limit}, while for the weak case, there is a penalty of around 26\% on the time step size allowed.

\section{The 2D elastic case}\label{section_2D_elastic_case}
Unlike the acoustic case, extension of the above procedure on strongly imposing the free surface boundary condition to the 2D elastic case is less straightforward and requires intricate design in the grid layout and the operators, which will be the focus of this section.

\subsection{Background} \label{section_background_elastic}
The 2D isotropic elastic wave system is presented in the following:
\begin{linenomath}
\begin{subequations}
\label{2D_elastic_wave_system}
\begin{empheq}[left=\empheqlbrace]{alignat = 2}
\displaystyle \rho \frac{\partial v_x}{\partial t} \enskip &= \enskip \displaystyle \frac{\partial \sigma_{xx}}{\partial x} + \frac{\partial \sigma_{xy}}{\partial y}  \, ; 
\label{Elastic_wave_equation_2D_analysis_Vx} \\
\displaystyle \rho \frac{\partial v_y}{\partial t} \enskip &= \enskip \displaystyle \frac{\partial \sigma_{xy}}{\partial x} + \frac{\partial \sigma_{yy}}{\partial y}  \, ; 
\label{Elastic_wave_equation_2D_analysis_Vy} \\
\displaystyle s_{xxkl} \frac{\partial \sigma_{kl}}{\partial t} \enskip &= \enskip \displaystyle \frac{\partial v_x}{\partial x}; 
\label{Elastic_wave_equation_2D_analysis_Sxx} \\
\displaystyle s_{xykl} \frac{\partial \sigma_{kl}}{\partial t} \enskip &= \enskip \displaystyle \frac{1}{2} \left( \frac{\partial v_y}{\partial x} + \frac{\partial v_x}{\partial y} \right) \, ; 
\label{Elastic_wave_equation_2D_analysis_Sxy} \\
\displaystyle s_{yykl} \frac{\partial \sigma_{kl}}{\partial t} \enskip &= \enskip \displaystyle \frac{\partial v_y}{\partial y},
\label{Elastic_wave_equation_2D_analysis_Syy}
\end{empheq}
\end{subequations}
\end{linenomath}
where velocity components ($v_x$ and $v_y$) and stress components ($\sigma_{xx}$, $\sigma_{xy}$, and $\sigma_{yy}$) are the sought solution variables; density ($\rho$) and compliance tensor ($s_{ijkl}$) are given physical parameters.
%
%
%
In \eqref{2D_elastic_wave_system}, the Einstein summation convention applies to the subscript indices $k$ and $l$, which go through $x$ and $y$.
%

We note here that in \eqref{2D_elastic_wave_system}, the constitution relation is expressed via the compliance tensor, which is convenient for analysis and derivation. 
Its equivalent form, where the constitution relation is expressed via the stiffness tensor, is convenient for implementation. 
The analysis and derivation based on \eqref{2D_elastic_wave_system} can be easily translated to its equivalent system. 
To conserve space, we omit their link here. Interest readers may referred to \cite[Appendix A]{gao2020explicit} for more detail.

\begin{figure}[H]
\captionsetup{width=1\textwidth, font=small,labelfont=small}
\centering\includegraphics[scale=0.0875]{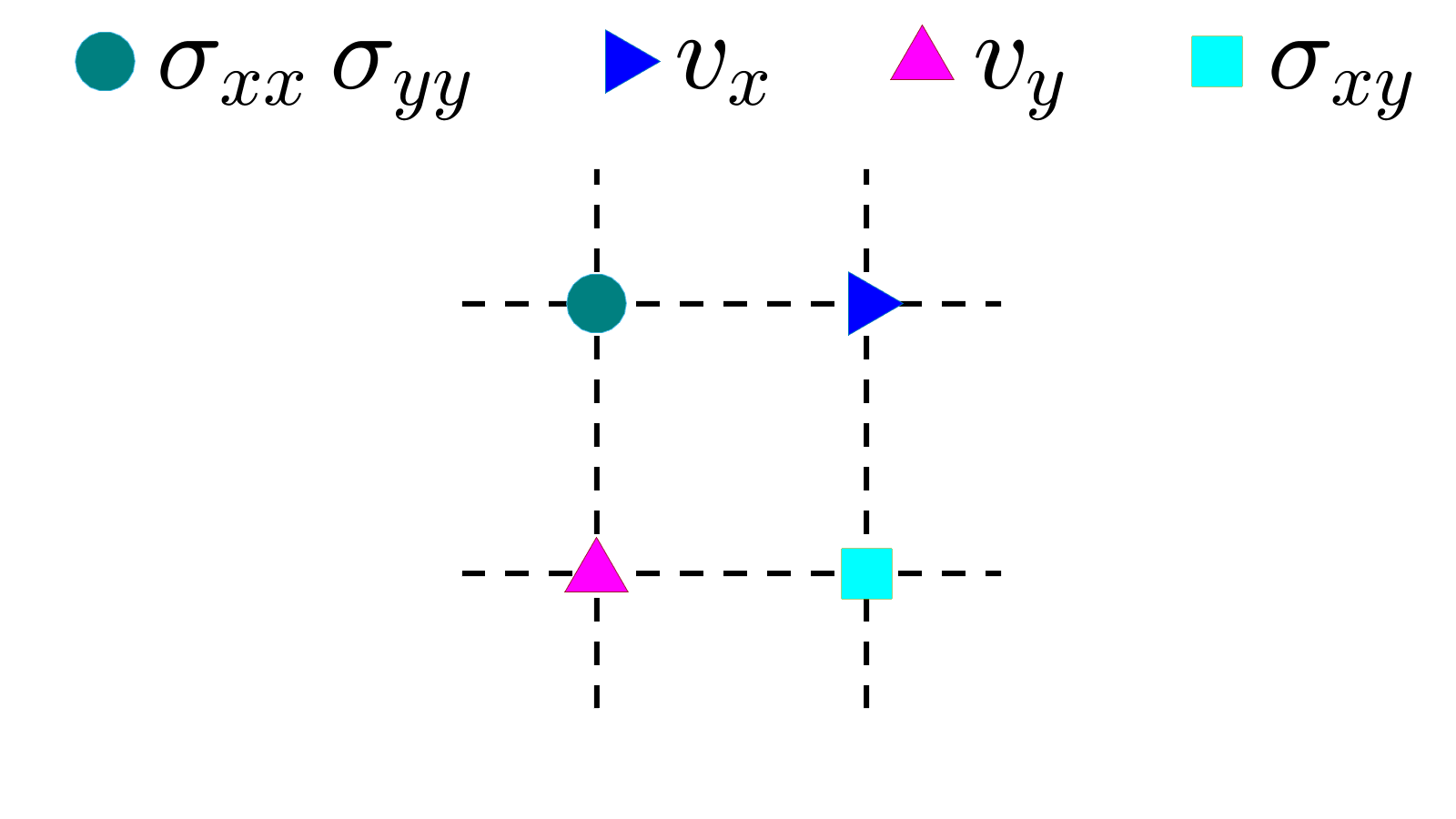}
\vspace{-1.em}
\caption{Standard grid layout for the 2D elastic wave equation without considering boundaries.} 
\label{Figure_grid_layout_elastic_without_boundary}
\end{figure}

As illustrated in Figure \ref{Figure_grid_layout_elastic_without_boundary}, when discretized on staggered grids using the standard grid layout, the normal stress components $\sigma_{xx}$ and $\sigma_{yy}$ are placed on the same location, 
the horizontal velocity $v_x$ and vertical velocity $v_y$ are shifted from the normal stress components on $x$- and $y$- directions, respectively, 
the shear stress component $\sigma_{xy}$ is shifted from the normal stress components on both directions.

Now considering the top boundary for example, the free surface boundary condition stipulates $\sigma_{xy} = \sigma_{yy} = 0$.
Strongly imposing these boundary conditions is non-trivial in the standard staggered grid setting for two reasons.
First, due to the grid staggering, the two stress components involved in the free surface boundary condition (i.e., $\sigma_{xy}$ and $\sigma_{yy}$) are not on the surface simultaneously.
Second, among the two normal stress components that share the same grid (i.e., $\sigma_{xx}$ and $\sigma_{yy}$) and hence share the same operators, one of them, namely, $\sigma_{xx}$, is not constrained by the free surface boundary condition. Strongly imposing the condition $\sigma_{yy} = 0$ by modifying the operators as described in section \ref{subsection_1D_strong_imposition} will lead to erroneous update to $\sigma_{xx}$.

\subsection{Grid layout} \label{section_grid_layout}
\begin{figure}[H]
\captionsetup{width=1\textwidth, font=small,labelfont=small}
\centering\includegraphics[scale=0.0875]{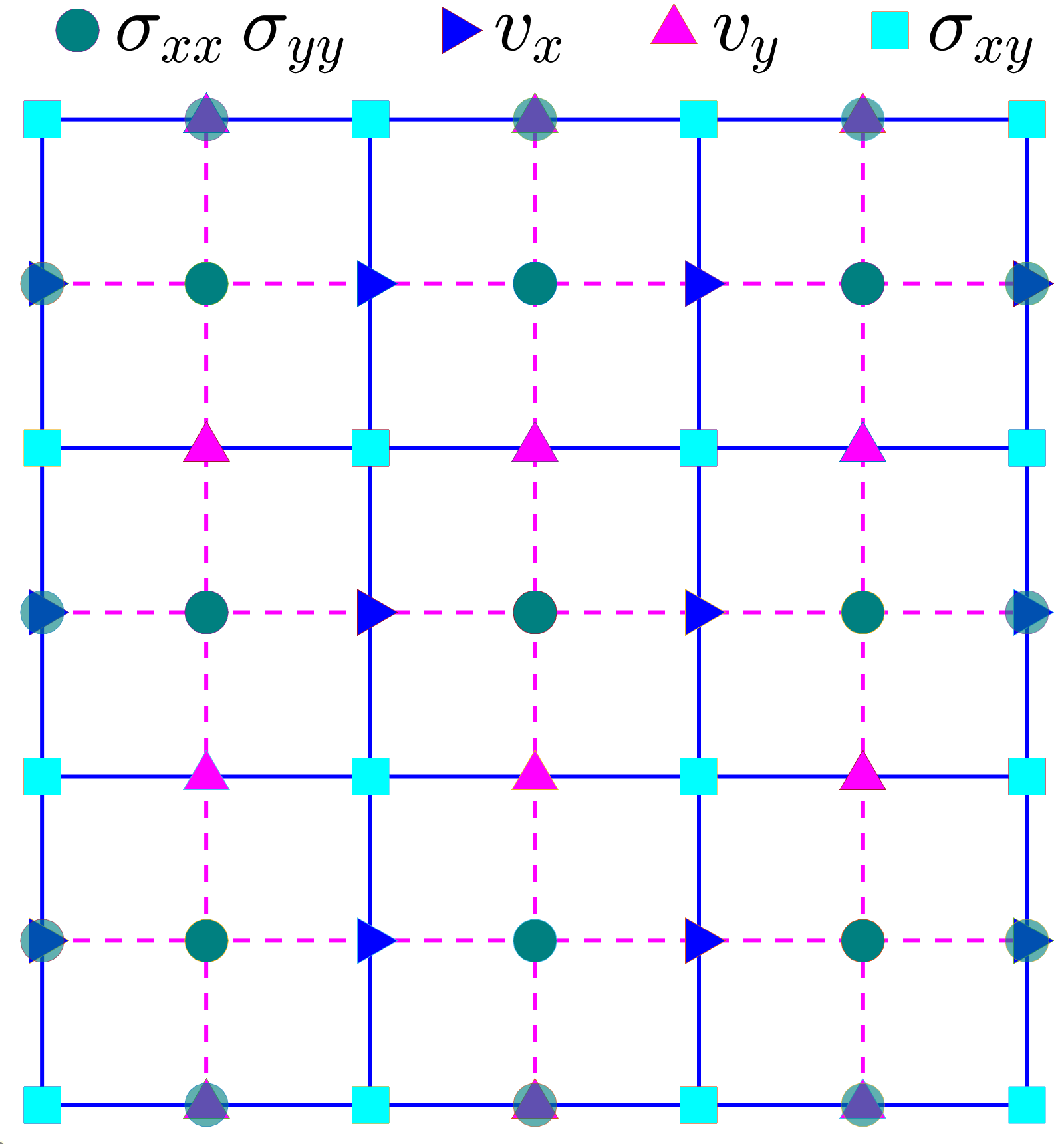}
\caption{Grid layout for the 2D elastic wave equation suitable for strongly imposing the free surface boundary condition.} 
\label{Elastic_2D_grid_layout_full_FSBC_shaded}
\end{figure}

In order to strongly impose the free surface boundary condition, we have to modify the grid layout slightly at the boundary. 
We find the grid layout illustrated in Figure \ref{Elastic_2D_grid_layout_full_FSBC_shaded} a suitable choice.
In Figure \ref{Elastic_2D_grid_layout_full_FSBC_shaded}, the shear stress component $\sigma_{xy}$ is aligned with the $N$-grid points on both directions. 
Again taking the top boundary for example, we have inserted an extra grid point occupied by the normal stress components $\sigma_{xx}$ and $\sigma_{yy}$ on the boundary.
Now both $\sigma_{xy}$ and $\sigma_{yy}$ are present on the top surface. This addresses the first hindrance mentioned above.
For the purpose of imposing the free surface boundary condition strongly, modifying the operators as described in section \ref{subsection_1D_strong_imposition} will still lead to erroneous updates to $\sigma_{xx}$ on the top surface. 
However, noticing that $\sigma_{xx}$ is only used to update $v_x$; see \eqref{2D_elastic_wave_system}, which is not present on the top surface. 
This means that the erroneous updates to $\sigma_{xx}$ will not have an effect on the rest of the simulation; $\sigma_{xx}$ on the top surface is essentially isolated from the solution variables at the other grid points.
This addresses the second hindrance mentioned above.

In fact, these extra grid points can be considered as only {\it conceptually} there, i.e., they can be omitted in {\it practical} implementations. 
This is because, on one hand, $\sigma_{xx}$ on the top surface is never used. 
On the other hand, $\sigma_{yy}$ is zero on the top surface and hence always has zero contribution to the variables it is used to update. This is equivalent to not taking its contribution by truncating the difference operators (instead of setting the rows and columns to zero as described in section \ref{subsection_1D_strong_imposition}).

\begin{remark}
\label{remark_truncate}
The omission mentioned above lead to some inconsistency between how the grid columns for $\sigma_{xy}$ and $v_x$ and the grid columns for $\sigma_{xx}$, $\sigma_{yy}$, and $v_y$ are treated.
For the convenience of discussion, we refer to the former as $\mathbb N$-grid columns and the later as $\mathbb M$-grid columns henceforward.
Namely, on $\mathbb N$-grid columns, the operators are reset to zero but not truncated; on $\mathbb M$-grid columns, the operators are truncated and grid points on the surface are omitted.

To reflect this choice, these extra grid points are shaded in Figure \ref{Elastic_2D_grid_layout_full_FSBC_shaded}.
This may cause slight confusion but brings the benefit that the grid points and hence the sizes of the operators on these two grid columns are the same, which greatly simplifies the analysis and implementation. 
As we will see later, by design, the operators associated with these two grid columns can be made to be {\normalfont effectively} the same, which is convenient for practical implementation. 
\end{remark}

\subsection{SBP operators} \label{section_SBP_operators_elastic}
As alluded to in Remark \ref{remark_truncate}, we still need to design the suitable SBP operators for this grid layout with the strongly imposed free surface boundary condition in mind.
We note here that although $\sigma_{yy}$ on the top surface is omitted in numerical implementation, it still needs to be formally included in the stencil. Otherwise, the effect of the boundary condition $\sigma_{yy} = 0$ will not be imposed.

To start, we conduct discrete energy analysis for the elastic wave system \eqref{2D_elastic_wave_system} to gain insights on what are the requirements for these SBP operators.
The semi-discretized system is presented below:

\begin{subequations}
\label{Semi_discretized_elastic_wave_equation_2D_analysis}
\begin{empheq}[left=\empheqlbrace]{alignat = 2}
\displaystyle \mathcal A^{V_x} \boldsymbol{\rho}^{V_x} \frac{d V_x}{d t} &\enskip = \enskip \displaystyle \mathcal A^{V_x} \mathcal D^{\Sigma_{xx}}_x \Sigma_{xx} + \mathcal A^{V_x} \mathcal D^{\Sigma_{xy}}_y \Sigma_{xy} \, ; 
\label{Discretized_elastic_wave_equation_2D_analysis_Vx} \\
\displaystyle \mathcal A^{V_y} \boldsymbol{\rho}^{V_y} \frac{d V_y}{d t} &\enskip = \enskip \displaystyle \mathcal A^{V_y} \mathcal D^{\Sigma_{xy}}_x \Sigma_{xy} + \mathcal A^{V_y} \mathcal D^{\Sigma_{yy}}_y \Sigma_{yy} \, ;
\label{Discretized_elastic_wave_equation_2D_analysis_Vy} \\
\displaystyle \mathcal A^{\Sigma_{xx}} S^{\Sigma_{kl}}_{xxkl} \frac{d \Sigma_{kl}}{d t} &\enskip = \enskip \displaystyle \mathcal A^{\Sigma_{xx}} \mathcal D^{V_x}_x V_x \, ; 
\label{Discretized_elastic_wave_equation_2D_analysis_Sxx} \\
\displaystyle \mathcal A^{\Sigma_{xy}} S^{\Sigma_{kl}}_{xykl} \frac{d \Sigma_{kl}}{d t} &\enskip = \enskip \displaystyle \frac{1}{2} \mathcal A^{\Sigma_{xy}} \left( \mathcal D^{V_y}_x V_y + \mathcal D^{V_x}_y V_x \right) \, ; 
\label{Discretized_elastic_wave_equation_2D_analysis_Sxy} \\
\displaystyle \mathcal A^{\Sigma_{yy}} S^{\Sigma_{kl}}_{yykl} \frac{d \Sigma_{kl}}{d t} &\enskip = \enskip \displaystyle \mathcal A^{\Sigma_{yy}} \mathcal D^{V_y}_y V_y \, ,
\label{Discretized_elastic_wave_equation_2D_analysis_Syy}
\end{empheq}
\end{subequations}
In \eqref{Semi_discretized_elastic_wave_equation_2D_analysis}, superscripts such as $^{V_x}$ indicate the grids with which the underlying quantities or operators are associated; the Einstein summation convention still applies, but not to those appearing in these superscripts. 
%

The 2D SBP operators appearing in \eqref{Semi_discretized_elastic_wave_equation_2D_analysis} are constructed from their 1D counterparts via tensor product as
\begin{subequations}
\label{Elastic_2D_norm_matrices}
\begin{alignat}{2}
\mathcal A^{\Sigma_{xy}}                            = \overset{\mathbb N}{\mathcal A^N_x} \otimes \overset{\mathbb N}{\mathcal A^N_y} , \quad
&
\mathcal A^{V_x}                                    = \overset{\mathbb M}{\mathcal A^N_x} \otimes \overset{\mathbb N}{\mathcal A^M_y} , \quad
&
\mathcal A^{V_y}                                    = \overset{\mathbb N}{\mathcal A^M_x} \otimes \overset{\mathbb M}{\mathcal A^N_y} , \quad
&
\mathcal A^{\Sigma_{xx}} = \mathcal A^{\Sigma_{yy}} = \overset{\mathbb M}{\mathcal A^M_x} \otimes \overset{\mathbb N}{\mathcal A^M_y} ,
\end{alignat}
\end{subequations}
and
\begin{subequations}
\label{Elastic_2D_difference_operators}
\begin{alignat}{4}
& 
\mathcal D_{x}^{\Sigma_{xy}} = \overset{\mathbb N}{\mathcal D^N_x} \otimes \mathcal I^N_y , \quad 
& &
\mathcal D_{x}^{V_x}         = \overset{\mathbb M}{\mathcal D^N_x} \otimes \mathcal I^M_y , \quad 
& &
\mathcal D_{x}^{V_y}         = \overset{\mathbb N}{\mathcal D^M_x} \otimes \mathcal I^N_y , \quad 
& &
\mathcal D_{x}^{\Sigma_{xx}} = \overset{\mathbb M}{\mathcal D^M_x} \otimes \mathcal I^M_y , \\
& 
\mathcal D_{y}^{\Sigma_{xy}} = \mathcal I^N_x \otimes \overset{\mathbb N}{\mathcal D^N_y} , \quad 
& &
\mathcal D_{y}^{V_x}         = \mathcal I^N_x \otimes \overset{\mathbb N}{\mathcal D^M_y} , \quad 
& &
\mathcal D_{y}^{V_y}         = \mathcal I^M_x \otimes \overset{\mathbb M}{\mathcal D^N_y} , \quad 
& &
\mathcal D_{y}^{\Sigma_{yy}} = \mathcal I^M_x \otimes \overset{\mathbb M}{\mathcal D^M_y} .
\end{alignat}
\end{subequations}
Since the $\mathbb N$ and $\mathbb M$ grid rows and columns have different grid layouts, they are associated with different sets of 1D SBP operators. 
In \eqref{Elastic_2D_norm_matrices} and \eqref{Elastic_2D_difference_operators}, the overhead symbols $\mathbb N$ and $\mathbb M$ are attached to distinguish.

The discrete energy associated with \eqref{Semi_discretized_elastic_wave_equation_2D_analysis} is
\begin{equation}
\label{Discrete_energy_2D_elastic_wave_system}
E \ = \ \tfrac{1}{2} V_i^T \left( \mathcal A^{V_i}  \boldsymbol{\rho}^{V_i} \right) V_i^{\phantom{T}} 
\! + \, 
\tfrac{1}{2} \Sigma_{ij}^T \left( \mathcal A^{\Sigma_{ij}}  S^{\Sigma_{kl}}_{ijkl} \right) \Sigma_{kl}^{\phantom{T}}\ ,
\end{equation}
Taking the time derivative on both sides of \eqref{Discrete_energy_2D_elastic_wave_system} and substituting the equations from \eqref{Semi_discretized_elastic_wave_equation_2D_analysis}, we arrive at 
\begin{linenomath}
\begin{subequations}
\label{Energy_analysis_derivation_2D_elastic}
\begin{alignat}{4}
\frac{d E}{d t} 
\ = \ \ &
V_x^T \left[ 
\mathcal A^{V_x} \mathcal D^{\Sigma_{xx}}_x 
+
\left(\mathcal A^{\Sigma_{xx}} \mathcal D^{V_x}_x\right)^T 
\right] \Sigma_{xx}
\label{Energy_analysis_derivation_2D_elastic_Vx_Sxx}\\
%
\ + \ \ &
V_y^T \left[ 
\mathcal A^{V_y} \mathcal D^{\Sigma_{xy}}_x 
+
\left( \mathcal A^{\Sigma_{xy}} \mathcal D^{V_y}_x \right)^T 
\right] \Sigma_{xy}
\label{Energy_analysis_derivation_2D_elastic_Sxy_Vy}\\
%
\ + \ \ &
V_y^T \left[
\mathcal A^{V_y} \mathcal D^{\Sigma_{yy}}_y 
+
\left( \mathcal A^{\Sigma_{yy}} \mathcal D^{V_y}_y \right)^T
\right] \Sigma_{yy}
\label{Energy_analysis_derivation_2D_elastic_Sxy_Vx}\\
%
\ + \ \ &
V_x^T \left[ 
\mathcal A^{V_x} \mathcal D^{\Sigma_{xy}}_y
+
\left( \mathcal A^{\Sigma_{xy}} \mathcal D^{V_x}_y \right)^T
\right] \Sigma_{xy} \, .
\label{Energy_analysis_derivation_2D_elastic_Vy_Syy}
\end{alignat}
\end{subequations}
\end{linenomath}
Substituting the definitions from \eqref{Elastic_2D_norm_matrices} and \eqref{Elastic_2D_difference_operators} and invoking the tensor product properties, we arrive at
\begin{linenomath}
\begin{subequations}
\label{Energy_analysis_derivation_2D_elastic_1D_operators}
\begin{alignat}{4}
\frac{d E}{d t} 
\ = \ \ &
V_x^T \bigg[ 
\Big( \overset{\mathbb M}{\mathcal A^N_x} \overset{\mathbb M}{\mathcal D^M_x} \Big) 
\otimes \overset{\mathbb N}{\mathcal A^M_y} 
+ 
\Big( \overset{\mathbb M}{\mathcal A^M_x} \overset{\mathbb M}{\mathcal D^N_x} \Big)^T 
\otimes \overset{\mathbb M}{\mathcal A^M_y}
\bigg] \Sigma_{xx}
\label{Energy_analysis_derivation_2D_elastic_Vx_Sxx}\\[-0.5ex]
%
\ + \ \ &
\Sigma_{xy}^T \bigg[ 
\Big( \overset{\mathbb N}{\mathcal A^N_x} \overset{\mathbb N}{\mathcal D^M_x} \Big) 
\otimes \overset{\mathbb N}{\mathcal A^N_y} 
+ 
\Big( \overset{\mathbb N}{\mathcal A^M_x} \overset{\mathbb N}{\mathcal D^N_x} \Big)^T 
\otimes \overset{\mathbb M}{\mathcal A^N_y}
\bigg] V_y
\label{Energy_analysis_derivation_2D_elastic_Sxy_Vy}\\[-0.5ex]
%
\ + \ \ &
V_y^T \bigg[ 
\overset{\mathbb N}{\mathcal A^M_x} \otimes 
\Big( \overset{\mathbb M}{\mathcal A^N_y} \overset{\mathbb M}{\mathcal D^M_y} \Big) 
+ 
\overset{\mathbb M}{\mathcal A^M_x} \otimes 
\Big( \overset{\mathbb M}{\mathcal A^M_y} \overset{\mathbb M}{\mathcal D^N_y} \Big)^T 
\bigg] \Sigma_{yy}
\label{Energy_analysis_derivation_2D_elastic_Vy_Syy}\\[-0.5ex] 
%
\ + \ \ &
\Sigma_{xy}^T \bigg[ 
\overset{\mathbb N}{\mathcal A^N_x} \otimes 
\Big( \overset{\mathbb N}{\mathcal A^N_y} \overset{\mathbb N}{\mathcal D^M_y} \Big) 
+
\overset{\mathbb M}{\mathcal A^N_x} \otimes 
\Big( \overset{\mathbb N}{\mathcal A^M_y} \overset{\mathbb N}{\mathcal D^N_y} \Big)^T 
\bigg] V_x
\label{Energy_analysis_derivation_2D_elastic_Sxy_Vx}
.
\end{alignat}
\end{subequations}
\end{linenomath}
If we have the following properties
\begin{linenomath}
\begin{subequations}
\label{requirement_equal_norm_matrices}
\begin{alignat}{3}
\overset{\mathbb N}{\mathcal A^N_x} &= \overset{\mathbb M}{\mathcal A^N_x} &= \mathcal A^N_x \\
\overset{\mathbb N}{\mathcal A^N_y} &= \overset{\mathbb M}{\mathcal A^N_y} &= \mathcal A^N_y  
\end{alignat}
\end{subequations}
\end{linenomath}
and
\begin{linenomath} 
\begin{subequations}
\label{requirement_skew_symmetry}
\begin{alignat}{3}
\Big( \overset{\mathbb N}{\mathcal A^N_x} \overset{\mathbb N}{\mathcal D^M_x} \Big) 
+
\Big( \overset{\mathbb N}{\mathcal A^M_x} \overset{\mathbb N}{\mathcal D^N_x} \Big)^T &= \boldsymbol 0 
\\
\Big( \overset{\mathbb M}{\mathcal A^N_x} \overset{\mathbb M}{\mathcal D^M_x} \Big) 
+
\Big( \overset{\mathbb M}{\mathcal A^M_x} \overset{\mathbb M}{\mathcal D^N_x} \Big)^T &= \boldsymbol 0 
\\
\Big( \overset{\mathbb N}{\mathcal A^N_y} \overset{\mathbb N}{\mathcal D^M_y} \Big) 
+
\Big( \overset{\mathbb N}{\mathcal A^M_y} \overset{\mathbb N}{\mathcal D^N_y} \Big)^T &= \boldsymbol 0 
\\
\Big( \overset{\mathbb M}{\mathcal A^N_y} \overset{\mathbb M}{\mathcal D^M_y} \Big) 
+
\Big( \overset{\mathbb M}{\mathcal A^M_y} \overset{\mathbb M}{\mathcal D^N_y} \Big)^T &= \boldsymbol 0
\end{alignat}
\end{subequations}
\end{linenomath}
in the 1D operators, we have that $\frac{d E}{d t} = 0$ following \eqref{Energy_analysis_derivation_2D_elastic_1D_operators} and the property of tensor product, i.e., the discrete energy is conserved, which mimics the behavior of the continuous elastic wave system.

We note her that while \eqref{requirement_skew_symmetry} is expected from common practice, \eqref{requirement_equal_norm_matrices} is the additional requirement related to the special grid layout, namely, the two sets ($\mathbb N$ and $\mathbb M$) of 1D SBP operators need to share the same norm matrices, which implies that the design of these two sets of operators need to be interwined.

Such operators do exist. 
Moreover, it turns out that we can ask for an additional requirement on these operators, namely,
\begin{linenomath}
\begin{subequations}
\label{requirement_equal_difference_operators}
\begin{alignat}{3}
\overset{\mathbb N}{\mathcal D^N_x} &= \overset{\mathbb M}{\mathcal D^N_x} &= \mathcal D^N_x \\
\overset{\mathbb N}{\mathcal D^N_y} &= \overset{\mathbb M}{\mathcal D^N_y} &= \mathcal D^N_y \, ,  
\end{alignat}
\end{subequations}
\end{linenomath}
so that the implementation can be greatly simplified.
These operators are presented below%
\footnotemark
\footnotetext
{
Convergence studies of these operators are included in the {\it Supplementary Material} \ref*{supp_convergence_tests}. 
}
\begin{linenomath}
\begin{subequations}
\label{SBP_matrices_1D_intertwined}
\footnotesize
\begin{align}
\label{SBP_matrices_1D_intertwined_A}
&
\mathcal A^M = \left[ 
\arraycolsep=3.2pt
\begin{array}{r r r r r r r r r r r r r r r r r} 
\nicefrac{13}{12} &                 & \multicolumn{1}{r|}{} \\
		          & \nicefrac{7}{8} & \multicolumn{1}{r|}{} \\
			      &                 & \multicolumn{1}{r|}{ \nicefrac{25}{24} }  \\ \cline{1-3} \\[-1.75ex]
				  &  			    &    			& 		   1 &			  & \\
				  & 			    &    			& 			 & 			1 & \\[-0.75ex]
				  &  			    &    			& 			 &   		  & \ddots
\end{array} 
\right],
&
&
\mathcal A^N = \left[ 
\arraycolsep=3.2pt
\begin{array}{r r r r r r r r r r r r r r r r r} 
\nicefrac{7}{18} & 				   & 		   		& \multicolumn{1}{r|}{} \\
		         & \nicefrac{9}{8} &     			& \multicolumn{1}{r|}{} \\
  			     &	               & 1  			& \multicolumn{1}{r|}{} \\
			     &    			   &    			& \multicolumn{1}{r|}{ \nicefrac{71}{72} }    \\ \cline{1-4} \\[-1.75ex]
				 &   			   &    			&    			&  1  			&    		  \\
				 &    			   &  				&    			&    			&  1  		& \\[-0.75ex]
				 &  			   &    			&   		    &    			&      		& \ddots
\end{array} 
\right];
&
\allowdisplaybreaks
\\
\label{SBP_matrices_1D_intertwined_D}
&
\mathcal D^N = \left[ 
\arraycolsep=3.2pt
\begin{array}{r r r r r r r r r r r r r r r r r} 
-\nicefrac{79}{78} &  \nicefrac{27}{26} & -\nicefrac{1}{26}  &  \nicefrac{1}{78}  & \multicolumn{1}{r|}{0}      \\
 \nicefrac{2}{21}  & -\nicefrac{9}{7}   &  \nicefrac{9}{7}   & -\nicefrac{2}{21}  & \multicolumn{1}{r|}{0}      \\
 \nicefrac{1}{75}  &  0			        & -\nicefrac{27}{25} &  \nicefrac{83}{75} & \multicolumn{1}{r|}{ -\nicefrac{1}{25} } & \\ \cline{1-5}
				   &   					&  \nicefrac{1}{24}  & -\nicefrac{9}{8}   & \nicefrac{9}{8}   						 & -\nicefrac{1}{24} 		  & \\
				   &     				&    				 &  \nicefrac{1}{24}  & -\nicefrac{9}{8} 	  					 &  \nicefrac{9}{8}  		  & -\nicefrac{1}{24} \\[-0.5ex]
				   & 					& 					 & 					  & \multicolumn{1}{l}{\ddots} 				 & \multicolumn{1}{l}{\ddots} & \multicolumn{1}{l}{\ddots} & \multicolumn{1}{l}{\ddots}  
\end{array} 
\right], 
&
&
\mathcal D^M = \left[ 
\arraycolsep=3.2pt
\begin{array}{r r r r r r r r r r r r r r r r r} 
\nicefrac{79}{28} & -\nicefrac{3}{14} & -\nicefrac{1}{28}  &  0    			    & \multicolumn{1}{r|}{0}   \\
-1       		  &  1    			  &  0    			   &  0    			    & \multicolumn{1}{r|}{0}   \\
 \nicefrac{1}{24} & -\nicefrac{9}{8}  &  \nicefrac{9}{8}   & -\nicefrac{1}{24}  & \multicolumn{1}{r|}{0} & \\
-\nicefrac{1}{71} &  \nicefrac{6}{71} & -\nicefrac{83}{71} &  \nicefrac{81}{71} & \multicolumn{1}{r|}{ -\nicefrac{3}{71} } & \\ \cline{1-5}
				  &    				  &  \nicefrac{1}{24}  & -\nicefrac{9}{8}   &  \nicefrac{9}{8}   		 			   & -\nicefrac{1}{24} 		    & \\
				  &     			  &    				   &  \nicefrac{1}{24}  & -\nicefrac{9}{8}   		 			   &  \nicefrac{9}{8}  		    & -\nicefrac{1}{24} \\[-0.5ex]
				  & 				  & 				   & 					& \multicolumn{1}{l}{\ddots} 			   & \multicolumn{1}{l}{\ddots} & \multicolumn{1}{l}{\ddots} & \multicolumn{1}{l}{\ddots}  
\end{array} 
\right].
&
\end{align}
\end{subequations}
\end{linenomath}
Comparing with \eqref{SBP_matrices_1D_extrapolating}, the only difference is in the first row of $\mathcal D^M$, which requires some explanation.

When interpreted as the $\mathbb M$ operators, these operators can be used directly since they are already truncated; see Remark \ref{remark_truncate}.
Taking the top surface for example, the first row of $\mathcal D^M$ applies on $\sigma_{yy}$ to update $v_y$ on the surface. 
The full stencil before truncation is actually $\small[-\nicefrac{18}{7} \,\, \nicefrac{79}{28} \,\, -\nicefrac{3}{14} \,\, -\nicefrac{1}{28} \,\, 0 \,\, 0 \,\, \cdots ]$ with an extra $-\nicefrac{18}{7}$ on the left, which corresponds to $\sigma_{yy}$ on the top surface. It can be omitted because of the free surface boundary condition, i.e., $\sigma_{yy} = 0$, on the top surface. This stencil is second order accurate, same as those for the other grid points within the boundary region.
When interpreted as $\mathbb N$ operators, we need to reset the rows in $D^M$ and columns in $D^N$ that corresponding to the free surface, as described in section \ref{subsection_1D_strong_imposition}.
With this resetting in mind, the stencil mentioned above (i.e., the first row of $\mathcal D^M$) is not asked to provide accurate update for $\sigma_{xy}$ on the surface.

We note here that the above operators are not the unique choice. 
Specifically, the additional requirement \eqref{requirement_equal_difference_operators} is optional and included here for the convenience of implementation. 
One can potentially use the extra degrees of freedom for other properties that one may desire.

\subsection{Numerical example} \label{section_numerical_example_elastic}
\begin{figure}[H]
\captionsetup{width=1\textwidth, font=small,labelfont=small}
\centering
\begin{subfigure}[b]{1\textwidth}
\captionsetup{width=1\textwidth, font=small,labelfont=small}
\centering\includegraphics[scale=0.3]{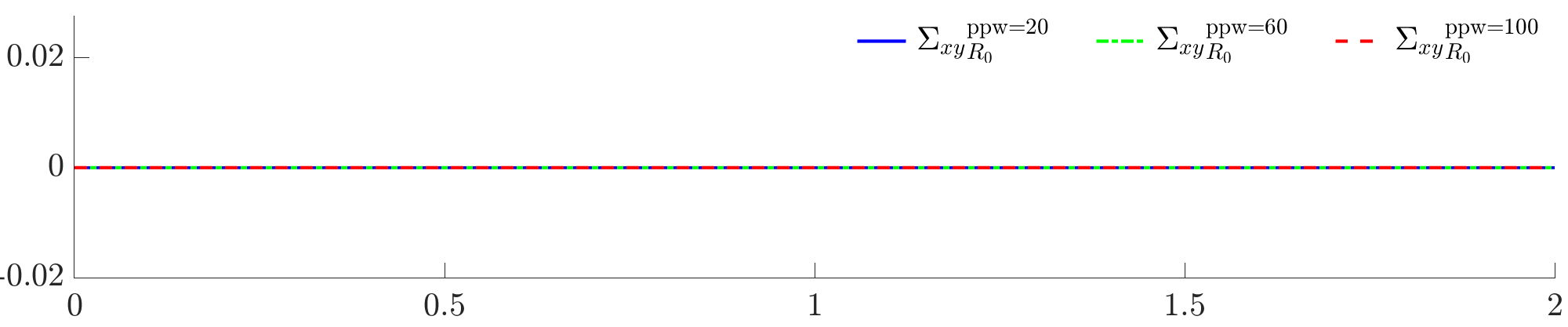}
\caption{Time history of $\Sigma_{xy}$ at the surface.}
\label{Elastic_full_FSBC_Sxy_0_strong}
\end{subfigure}\hfill
\\[2ex]
\begin{subfigure}[b]{1\textwidth}
\captionsetup{width=1\textwidth, font=small,labelfont=small}
\centering\includegraphics[scale=0.3]{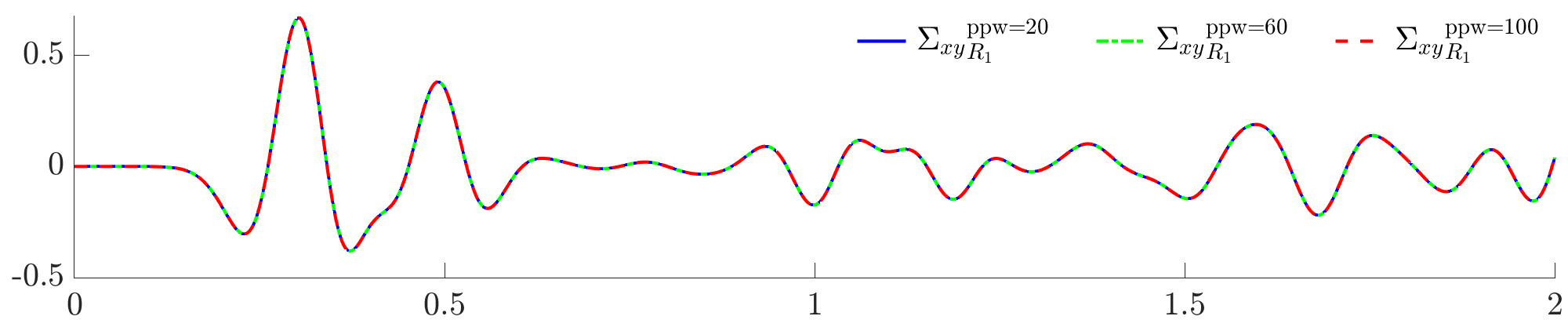}
\caption{Time history of $\Sigma_{xy}$ at 1 grid points below the surface for the case ppw = 20.}
\label{Elastic_full_FSBC_Sxy_1_strong}
\end{subfigure}\hfill
\\[2ex]
\begin{subfigure}[b]{1\textwidth}
\captionsetup{width=1\textwidth, font=small,labelfont=small}
\centering\includegraphics[scale=0.3]{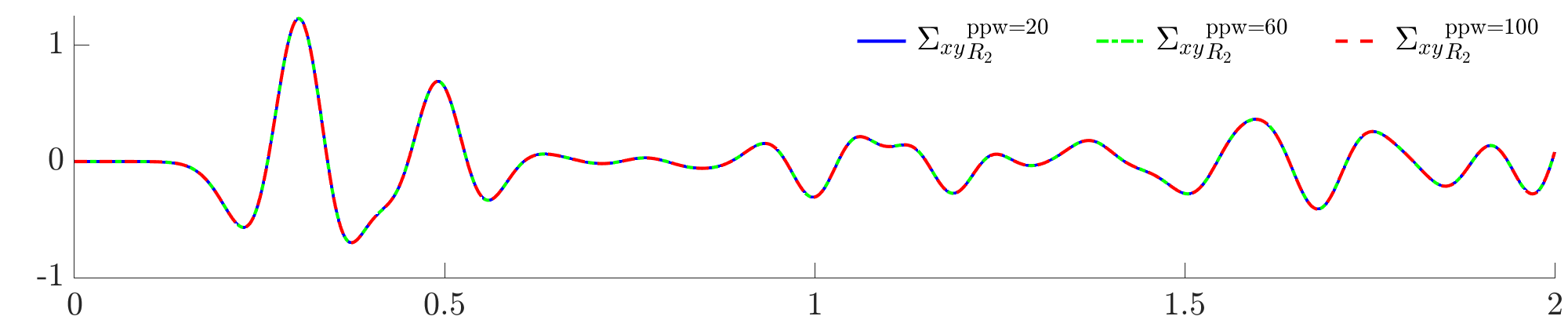}
\caption{Time history of $\Sigma_{xy}$ at 2 grid points below the surface for the case ppw = 20.}
\label{Elastic_full_FSBC_Sxy_2_strong}
\end{subfigure}\hfill
\caption{Time histories of the 2D elastic experiments. The free surface boundary condition is imposed strongly.
}
\label{Elastic_full_FSBC_Sxy_strong}
\end{figure}

\begin{figure}[H]
\captionsetup{width=1\textwidth, font=small,labelfont=small}
\centering
\begin{subfigure}[b]{1\textwidth}
\captionsetup{width=1\textwidth, font=small,labelfont=small}
\centering\includegraphics[scale=0.3]{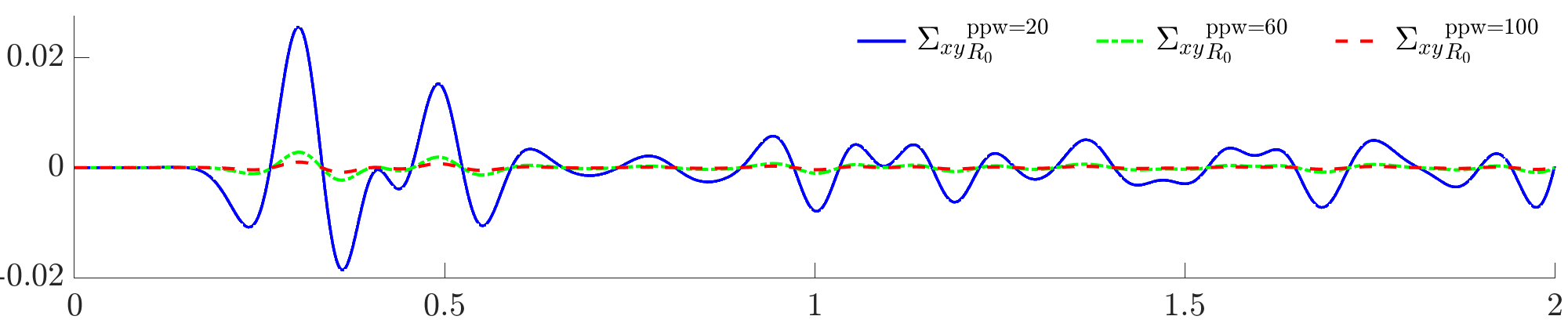}
\caption{Time history of $\Sigma_{xy}$ at the surface.}
\label{Elastic_full_FSBC_Sxy_0_weak}
\end{subfigure}\hfill
\\[2ex]
\begin{subfigure}[b]{1\textwidth}
\captionsetup{width=1\textwidth, font=small,labelfont=small}
\centering\includegraphics[scale=0.3]{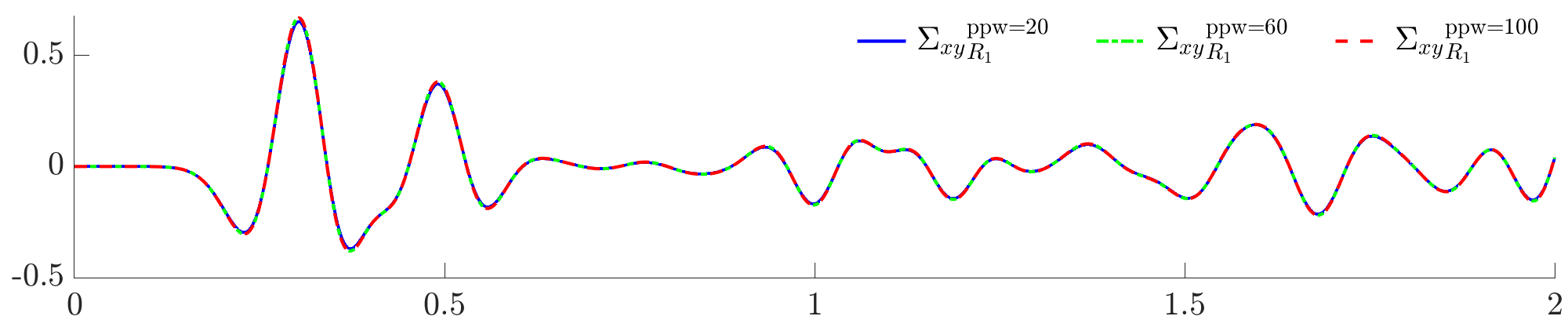}
\caption{Time history of $\Sigma_{xy}$ at 1 grid points below the surface for the case ppw = 20.}
\label{Elastic_full_FSBC_Sxy_1_weak}
\end{subfigure}\hfill
\\[2ex]
\begin{subfigure}[b]{1\textwidth}
\captionsetup{width=1\textwidth, font=small,labelfont=small}
\centering\includegraphics[scale=0.3]{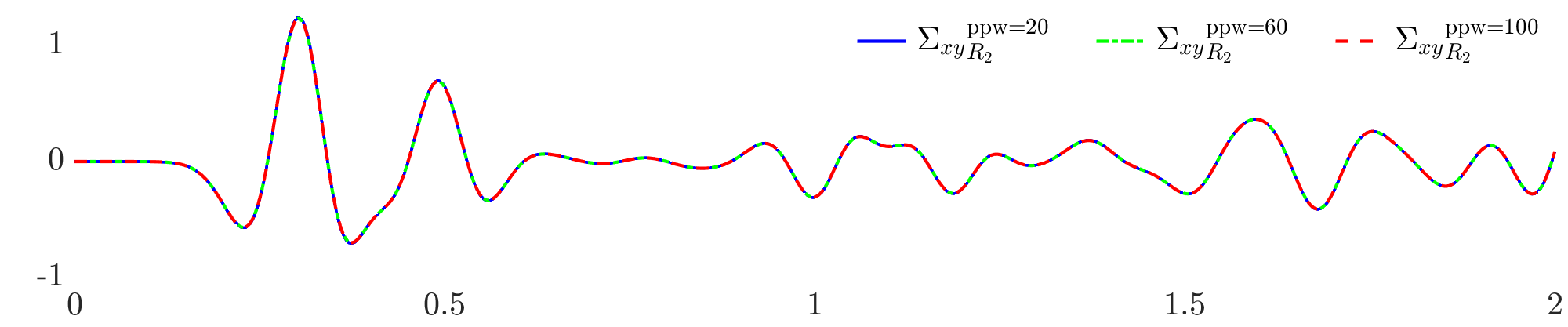}
\caption{Time history of $\Sigma_{xy}$ at 2 grid points below the surface for the case ppw = 20.}
\label{Elastic_full_FSBC_Sxy_2_weak}
\end{subfigure}\hfill
\caption{Time histories of the 2D elastic experiments. The free surface boundary condition is imposed weakly.
}
\label{Elastic_full_FSBC_Sxy_weak}
\end{figure}

Numerical experiments similar to that presented for the 2D acoustic case in section \ref{section_2D_acoustic_case} is conducted for the 2D elastic case.
The specification is mostly the same with a few exceptions. 
First, the grid points per minimal wavelength used here are increased to 20, 60, and 100. 
This follows the wisdom that elastic wave equation admits a surface wave mode that requires finer grid spacing to resolve than for the body waves; see \cite{kreiss2012boundary}. 
(However, experiments using 10, 30, 50 ppw reveal the same pattern.)
Moreover, the source is applied on $\sigma_{xx}$ and $\sigma_{yy}$ at one third of the simulation domain horizontally and half grid spacing below the surface vertically.
Time history of the solution variable $\Sigma_{xy}$ is recorded at two thirds of the simulation domain horizontally and at 0, 1 and 2 grid spacings below the surface vertically, which are shown in Figures \ref{Elastic_full_FSBC_Sxy_strong} and \ref{Elastic_full_FSBC_Sxy_weak}, for the strong and weak case, respectively. 
The time step length is specified as 2e-4 s for all simulations, and the simulation time is specified as 2 s.
To conserve space, time histories of the other solution components are omitted here and included in the {\it Supplementary Material} \ref*{supp_additional_figures_2D_elastic}.

From the above figures, we observe that both the strong and weak approaches deliver consistent and accurate results, with the strong approach offers only minor benefit at or near the free surface for the $\Sigma_{xy}$ component.
Their performance for the other solution components are even less distinguishable.

\begin{figure}[H]
\captionsetup{width=1\textwidth, font=small,labelfont=small}
\centering
\begin{subfigure}[b]{1\textwidth}
\captionsetup{width=1\textwidth, font=small,labelfont=small}
\centering\includegraphics[scale=0.3]{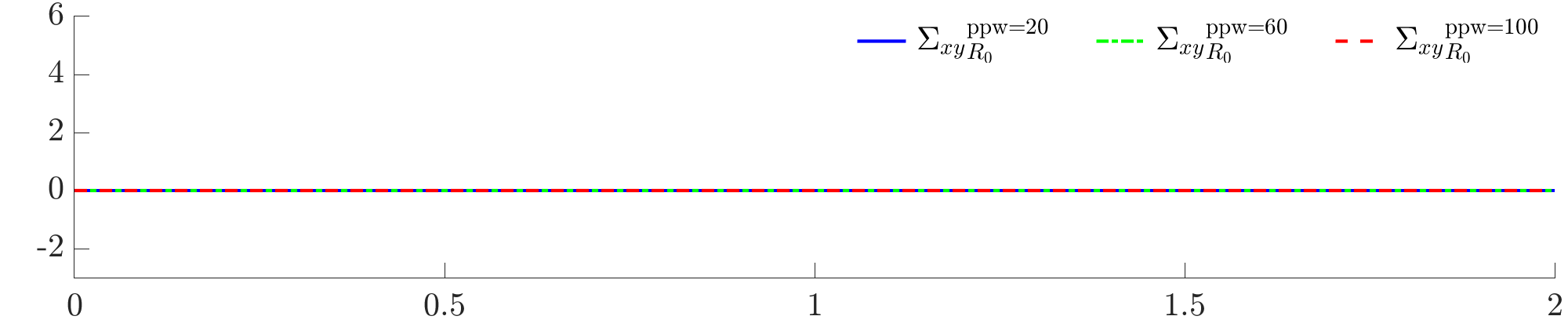}
\caption{Time history of $\Sigma_{xy}$ at the surface.}
\label{stress_test_Vy_elastic_full_FSBC_Sxy_0_strong}
\end{subfigure}\hfill
\\[2ex]
\begin{subfigure}[b]{1\textwidth}
\captionsetup{width=1\textwidth, font=small,labelfont=small}
\centering\includegraphics[scale=0.3]{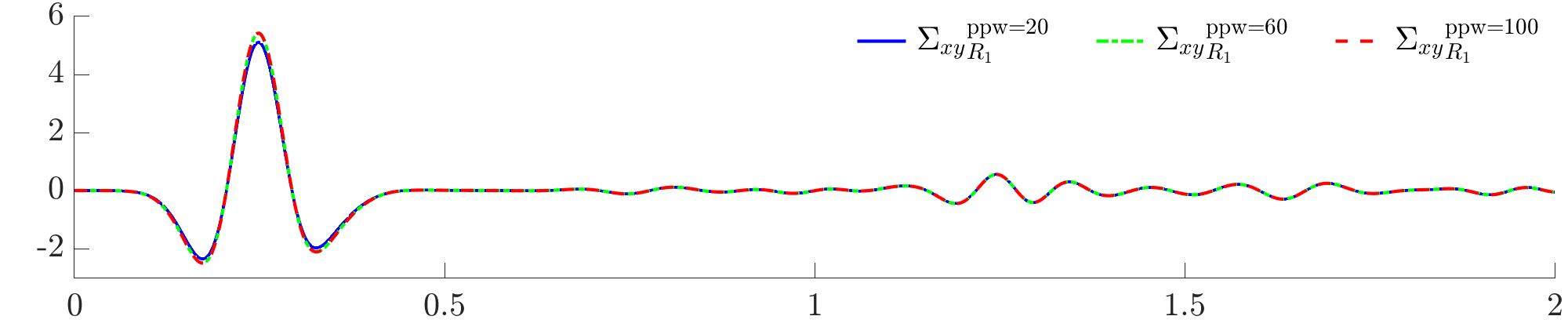}
\caption{Time history of $\Sigma_{xy}$ at 1 grid points below the surface for the case ppw = 20.}
\label{stress_test_Vy_elastic_full_FSBC_Sxy_1_strong}
\end{subfigure}\hfill
\\[2ex]
\begin{subfigure}[b]{1\textwidth}
\captionsetup{width=1\textwidth, font=small,labelfont=small}
\centering\includegraphics[scale=0.3]{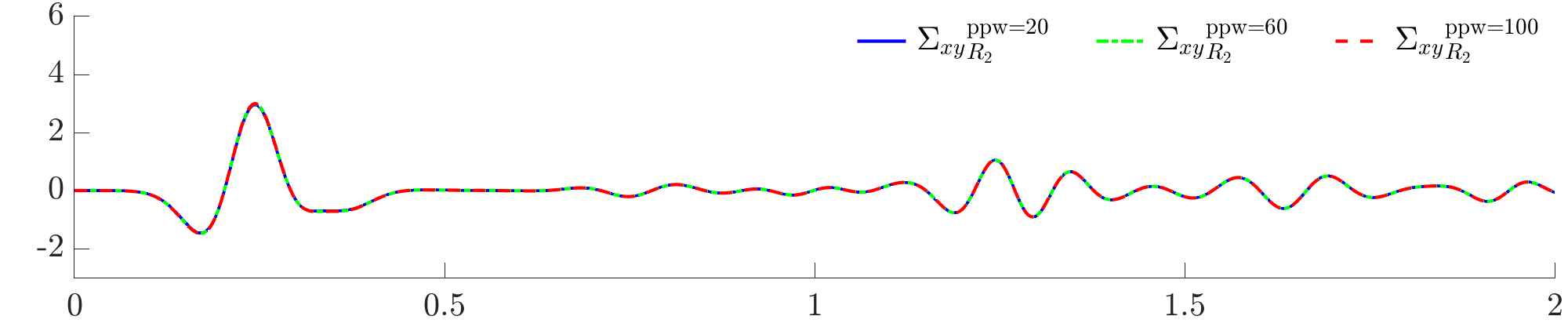}
\caption{Time history of $\Sigma_{xy}$ at 2 grid points below the surface for the case ppw = 20.}
\label{stress_test_Vy_elastic_full_FSBC_Sxy_2_strong}
\end{subfigure}\hfill
\caption{Time histories of the 2D elastic experiments. The free surface boundary condition is imposed strongly.
}
\label{stress_test_Vy_elastic_full_FSBC_Sxy_strong}
\end{figure}

\begin{figure}[H]
\captionsetup{width=1\textwidth, font=small,labelfont=small}
\centering
\begin{subfigure}[b]{1\textwidth}
\captionsetup{width=1\textwidth, font=small,labelfont=small}
\centering\includegraphics[scale=0.3]{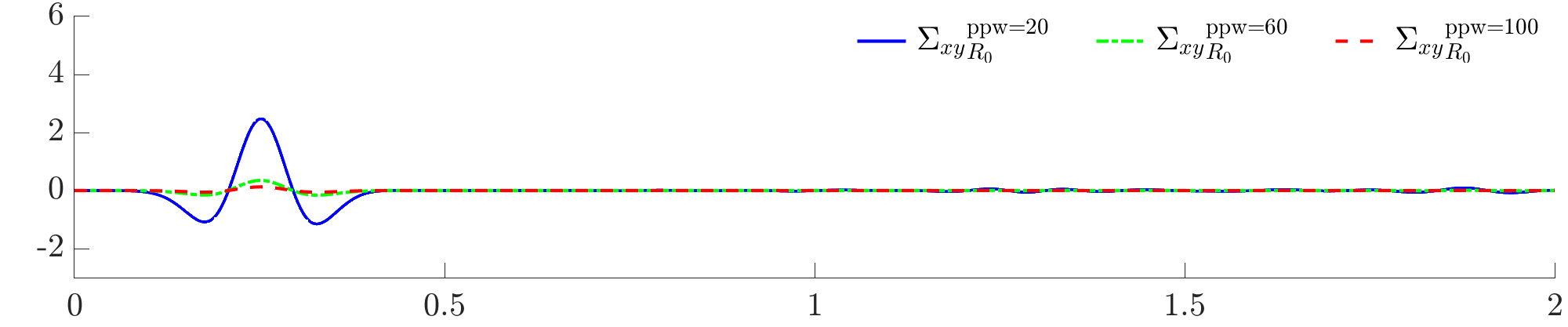}
\caption{Time history of $\Sigma_{xy}$ at the surface.}
\label{stress_test_Vy_elastic_full_FSBC_Sxy_0_weak}
\end{subfigure}\hfill
\\[2ex]
\begin{subfigure}[b]{1\textwidth}
\captionsetup{width=1\textwidth, font=small,labelfont=small}
\centering\includegraphics[scale=0.3]{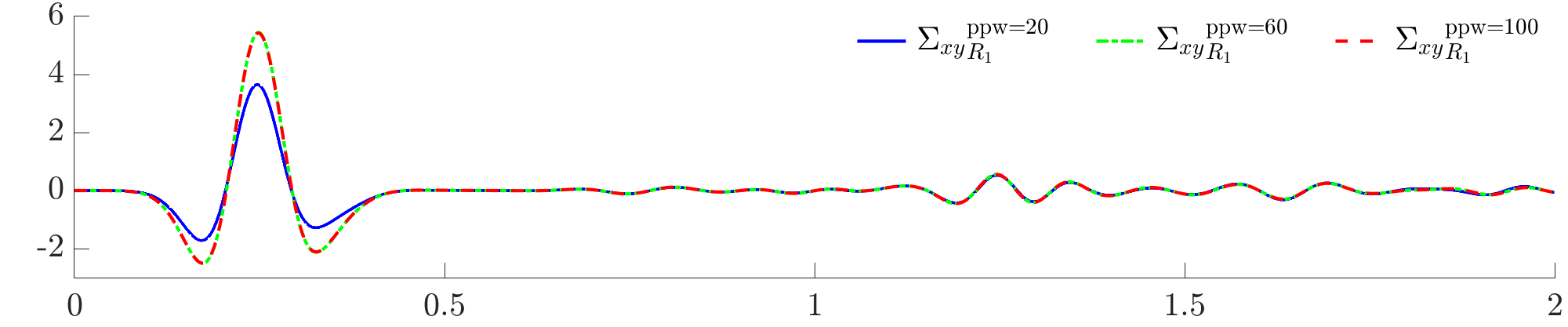}
\caption{Time history of $\Sigma_{xy}$ at 1 grid points below the surface for the case ppw = 20.}
\label{stress_test_Vy_elastic_full_FSBC_Sxy_1_weak}
\end{subfigure}\hfill
\\[2ex]
\begin{subfigure}[b]{1\textwidth}
\captionsetup{width=1\textwidth, font=small,labelfont=small}
\centering\includegraphics[scale=0.3]{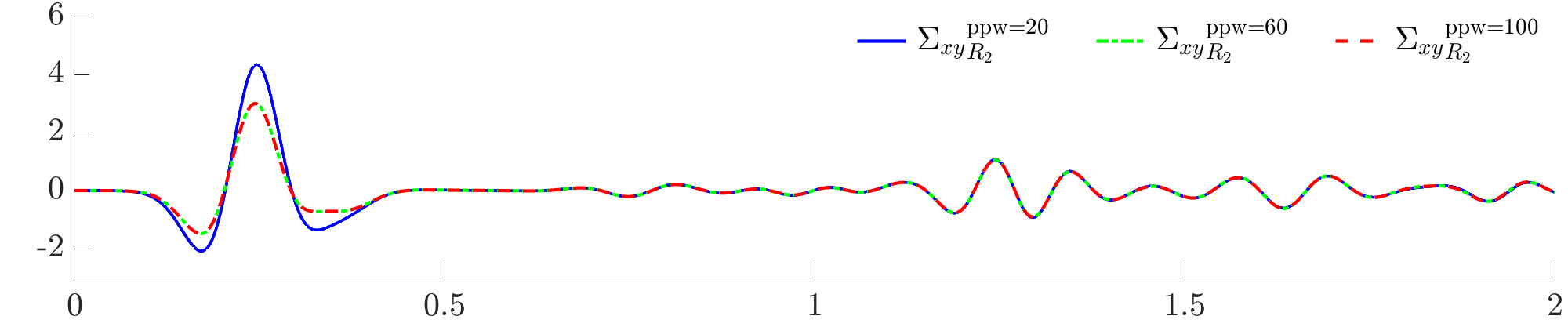}
\caption{Time history of $\Sigma_{xy}$ at 2 grid points below the surface for the case ppw = 20.}
\label{stress_test_Vy_elastic_full_FSBC_Sxy_2_weak}
\end{subfigure}\hfill
\caption{Time histories of the 2D elastic experiments. The free surface boundary condition is imposed weakly.
}
\label{stress_test_Vy_elastic_full_FSBC_Sxy_weak}
\end{figure}

A ``stress test'' of both approaches is conducted, where the source is applied on $V_y$ (i.e., a directional source) on the top surface and half grid spacing away from the left surface while the receivers are placed at 5 grid points (for the case ppw = 20) away from the source. 
In this case, the strong approach demonstrates noticeable improvement in the $\Sigma_{xy}$ component, as can be observed by comparing Figures \ref{stress_test_Vy_elastic_full_FSBC_Sxy_strong} and \ref{stress_test_Vy_elastic_full_FSBC_Sxy_weak}. 
There are little meaningful differences in the other recorded solution components, which are omitted here to conserve space but included in the {\it Supplementary Material} \ref*{supp_additional_figures_2D_elastic_stress_test_Vy}.

Another ``stress test'' is included in the {\it Supplementary Material} \ref*{supp_additional_figures_2D_elastic_stress_test_Sxy}, where the source is applied on $\Sigma_{xy}$ at the corner, more specifically, 1 grid spacing away from the top and left surfaces.
In this case, the strong approach performs markedly better than the weak approach for all five solution components.

Moreover, to accompany the above numerical experiments, convergence test results using a smooth manufactured solution for the elastic wave system \eqref{2D_elastic_wave_system} is included in {\it Supplementary Material} \ref{supp_2D_convergence_tests}, where the strong approach exhibits smaller error and higher convergence rates.

Finally, numerical experiments reveal that the CFL restriction for the strong and weak cases are $\nicefrac{6}{7}$ and 0.849, respectively. In other words, the weak case has a slight penalty on the time step size allowed while the strong case has none.

\begin{remark}
We note here that the results presented in this work is oriented towards a particular simulation setting (although appearing broadly in practice for various important wave-related applications), and should not be interpreted beyond its scope. 
Indeed, some authors have argued for the benefit of weakly imposed boundary conditions over the strongly imposed ones; see \cite{bazilevs2007weak} for an example in the fluid setting.
\end{remark}

\section{Conclusions} \label{section_conclusions}
We have considered the imposition of the free surface boundary condition for both the acoustic and elastic wave simulations in the framework of summation-by-parts finite difference discretizations.
Within this framework, the standard approach of imposing boundary conditions, i.e., weakly through appending simultaneous approximation terms, despite its flexibility, can have issues in certain practical situations, such as when a point source is placed near the surface, inducing abrupt changes in the wave field, which have been demonstrated in this work through numerical experiments.

As an alternative, we have demonstrated the possibility of imposing the free surface boundary condition strongly, i.e., by incorporating it into the summation-by-parts finite difference operators.
While easy to achieve for the acoustic case, this approach requires specially designed grid layout for the elastic case and operators that need to satisfy additional requirements, as revealed by the discrete energy analysis.
One possible set of such operators is presented in this work, which also takes into account the ease of implementation in its design.
Numerical experiments reveal that the strong approach can have significant advantages over the weak approach, in terms of both accuracy and time step size restrictions. 
Such advantages is particularly pronounced in the acoustic case.

The strong approach has its obvious disadvantage, namely, the boundary condition needs to be built into the operators and, therefore, renders the approach less flexible and less generally applicable.
However, in the particular case considered, i.e., free surface boundary condition for wave propagation, which has broad applications in geoscience and civil engineering, the additional effort is relatively small for the elastic case and trivial for the acoustic case.
For this reason, we recommend practitioners consider this option 
in their practical settings.

\appendix
\section{CFL restriction for the interior stencil and staggered leapfrog time integration scheme on unbounded domain}
\label{section_Courant_number_6_over_7}
\begin{figure}[H]
\vspace{-0.75em}
\captionsetup{width=1\textwidth, font=small,labelfont=small}
\centering\includegraphics[scale=0.075]{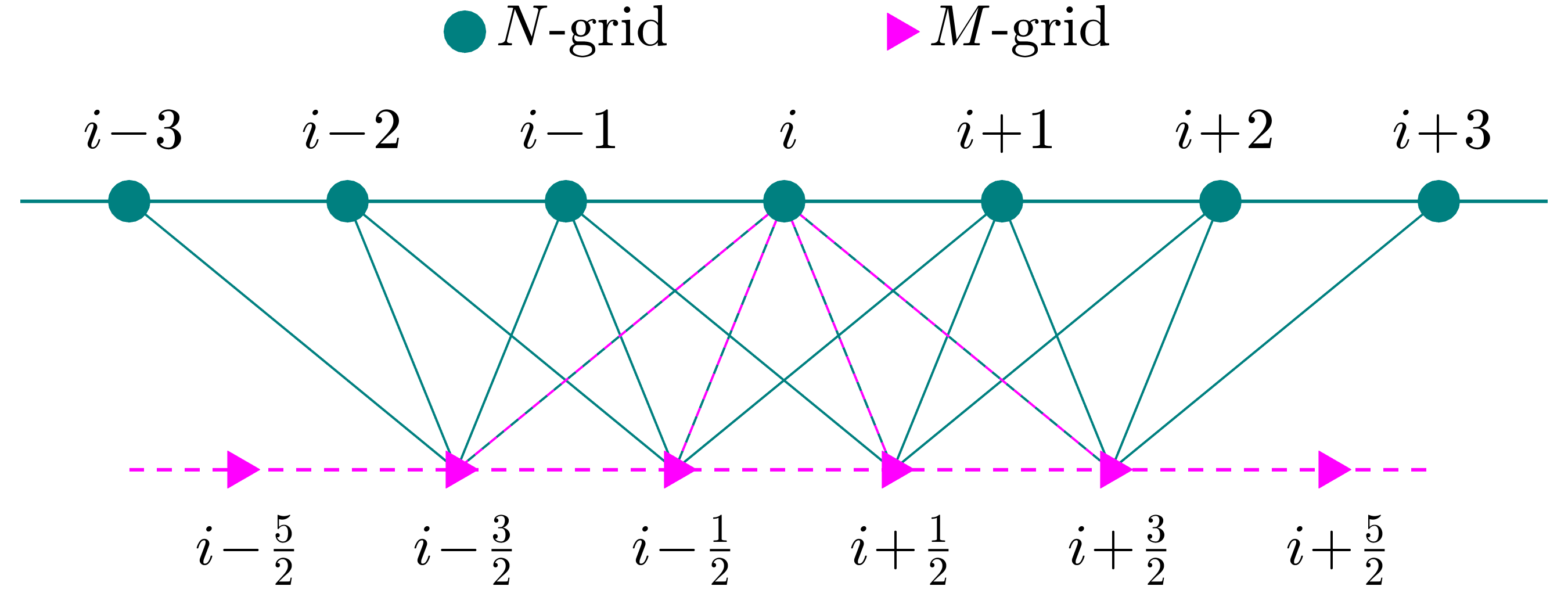}
\caption{
Illustration of the staggered spatial grids and the action of the stencil. 
%
}
\label{figure_Courant_number_grids_illustration}
\vspace{-0.75em}
\end{figure}

In this appendix, we consider the prototype 1D wave system \eqref{1D_wave_equation} with unit density and wave-speed on unbounded domain and derive the maximal Courant number associated with the standard fourth-order staggered grid spatial stencil 
$
\nicefrac{ [\nicefrac{1}{24}, \enskip -\nicefrac{9}{8}, \enskip \nicefrac{9}{8}, \enskip -\nicefrac{1}{24}] }{ \Delta x }
$
and the staggered leapfrog time integration scheme via von Neumann analysis.
%
%
Figure \ref{figure_Courant_number_grids_illustration} illustrates the staggered spatial grids and the action of the stencil.

In the following, we assume that variable $\sigma$ is discretized on the $N$-grid (which is associated with integer indices) and updated at integer time steps, while variable $v$ is discretized on the $M$-grid (which is associated with half-integer indices) and updated at half-integer time steps. 
Considering the approximation of $\sigma$ at the $i$th $N$-grid point, i.e., $\Sigma_i$, when updated from time step $n-1$ to $n$, the following formula is used:
\begin{equation}
\label{update_nm1_to_n}
\Sigma_{i}^{(n)} - \Sigma_{i}^{(n-1)} = \left( \tfrac{1}{24} V_{i-\frac{3}{2}}^{(n-\frac{1}{2})} - \tfrac{9}{8} V_{i-\frac{1}{2}}^{(n-\frac{1}{2})} + \tfrac{9}{8} V_{i+\frac{1}{2}}^{(n-\frac{1}{2})} - \tfrac{1}{24} V_{i+\frac{3}{2}}^{(n-\frac{1}{2})} \right) \tfrac{ \Delta t }{ \Delta x }.
\end{equation}
Similarly, when $\Sigma_i$ is updated from time step $n$ to $n+1$, the following formula is used:
\begin{equation}
\label{update_n_to_np1}
\Sigma_{i}^{(n+1)} - \Sigma_{i}^{(n)} = \left( \tfrac{1}{24} V_{i-\frac{3}{2}}^{(n+\frac{1}{2})} - \tfrac{9}{8} V_{i-\frac{1}{2}}^{(n+\frac{1}{2})} + \tfrac{9}{8} V_{i+\frac{1}{2}}^{(n+\frac{1}{2})} - \tfrac{1}{24} V_{i+\frac{3}{2}}^{(n+\frac{1}{2})} \right) \tfrac{ \Delta t }{ \Delta x }.
\end{equation}
Subtracting \eqref{update_nm1_to_n} from \eqref{update_n_to_np1} and substituting the update formulas for $V_{i-\frac{3}{2}}$, $V_{i-\frac{1}{2}}$, $V_{i+\frac{1}{2}}$, and $V_{i+\frac{3}{2}}$ from time step $n-\frac{1}{2}$ to $n+\frac{1}{2}$, we arrive at the following relation:
\begin{equation}
\label{update_relation_P}
\Sigma_i^{(n+1)} - 2 \Sigma_i^{(n)} + \Sigma_i^{(n-1)} = \left( 
\tfrac{1}{576} \Sigma_{i-3}^{(n)}
- \tfrac{3}{32} \Sigma_{i-2}^{(n)}          
+ \tfrac{87}{64} \Sigma_{i-1}^{(n)}
- \tfrac{365}{144} \Sigma_{i}^{(n)}
+ \tfrac{87}{64} \Sigma_{i+1}^{(n)}         
- \tfrac{3}{32} \Sigma_{i+2}^{(n)}        
+ \tfrac{1}{576} \Sigma_{i+3}^{(n)} \right) 
\left( \tfrac{ \Delta t }{ \Delta x } \right)^2.
\end{equation}

Considering a particular wave mode $p = e^{\alpha t} e^{\imath k x}$, where $\imath$ denotes the imaginary unit, we have, after simplification, the following relation:
\begin{equation}
\label{error_relation}
e^{\alpha \Delta t} - 2 + e^{- \alpha \Delta t} = 
\left(
 \tfrac{1}{72} \cos(k \Delta x)^3
-\tfrac{3}{8} \cos(k \Delta x)^2
+\tfrac{65}{24} \cos(k \Delta x)
-\tfrac{169}{72} \right) \left( \tfrac{ \Delta t }{ \Delta x } \right)^2.
\end{equation}
Using $\beta$ to denote $\left( \tfrac{1}{72} \cos(k \Delta x)^3 - \tfrac{3}{8} \cos(k \Delta x)^2 + \tfrac{65}{24} \cos(k \Delta x) - \tfrac{169}{72} \right)$ for brevity, we have that $\beta \in \left[-\tfrac{49}{9},0\right]$.
Further, denoting $\tfrac{ \Delta t }{ \Delta x }$ as $C$, we have $C>0$ for sensible discretizations. 
Finally, we denote $e^{\alpha \Delta t}$ as $G$, which is the growth factor between two neighboring time steps for 
the wave mode 
under consideration.
%
With these notations, we now have the following quadratic equation concerning $G$:
\begin{equation}
\label{quadratic_equation_G}
G^2 - \left(2 + \beta C^2 \right) G + 1 = 0\, ,
\end{equation}
%
which has two roots:
\begin{equation}
\label{roots_quadratic_equation_G}
G_1 = 1 + \tfrac{\beta C^2}{2} - \sqrt{ \left( 1 + \tfrac{\beta C^2}{2} \right)^2 - 1}
\text{\quad\enskip and \quad\enskip}
G_2 = 1 + \tfrac{\beta C^2}{2} + \sqrt{ \left( 1 + \tfrac{\beta C^2}{2} \right)^2 - 1}\, .
\end{equation}

Supposing that $\left( 1 + \tfrac{\beta C^2}{2} \right)^2 - 1 > 0$, we have that either $G_1\! < \!-1$, or $G_2\! > \!1$, i.e., this particular mode can grow out of bounds.
In the opposite case $\left( 1 + \tfrac{\beta C^2}{2} \right)^2 - 1 \leq 0$, we have that $\left| G_1 \right| = \left| G_2 \right| = 1$, 
i.e., this particular mode is stable.
For the overall simulation to be stable, $C$ needs to be chosen such that $\left( 1 + \tfrac{\beta C^2}{2} \right)^2 - 1 \leq 0$ holds for all $\beta \in \left[-\tfrac{49}{9},0\right]$.
This constraint is at its strictest when $\beta = -\tfrac{49}{9}$, which leads to $C \leq \nicefrac{6}{7}$, as stated in \eqref{Interior_CFL_limit}.

\addtolength{\bibsep}{-0.5em}

\renewcommand{\bibfont}{\normalfont\small}
\bibliographystyle{./bst_base/abbrv.bst}
\bibliography{refs}

\newpage
\supp
\input{supp}

%

\end{document}

%% file: supp.tex
\begin{center}
\Large
Supplementary Materials
\end{center}



\maketitle

\supp
\section{Additional 1D experiment with finer grid resolution}\label{supp_section_finer_grid}

\begin{figure}[H]
\captionsetup{width=1\textwidth, font=small,labelfont=small}
\centering
\begin{subfigure}[b]{1\textwidth}
\captionsetup{width=1\textwidth, font=small,labelfont=small}
\centering\includegraphics[scale=0.15]{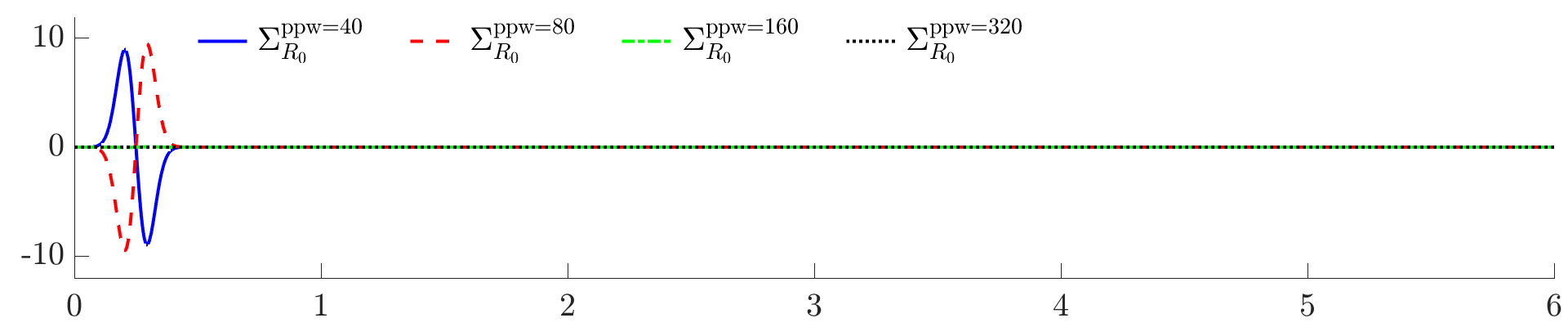}
\caption{Time histories of $\Sigma$ at $R_0$.} 
\label{supp_P_0_6_s_extrapolating_40_to_320}
\end{subfigure}\hfill
\\[2ex]
\begin{subfigure}[b]{1\textwidth}
\captionsetup{width=1\textwidth, font=small,labelfont=small}
\centering\includegraphics[scale=0.15]{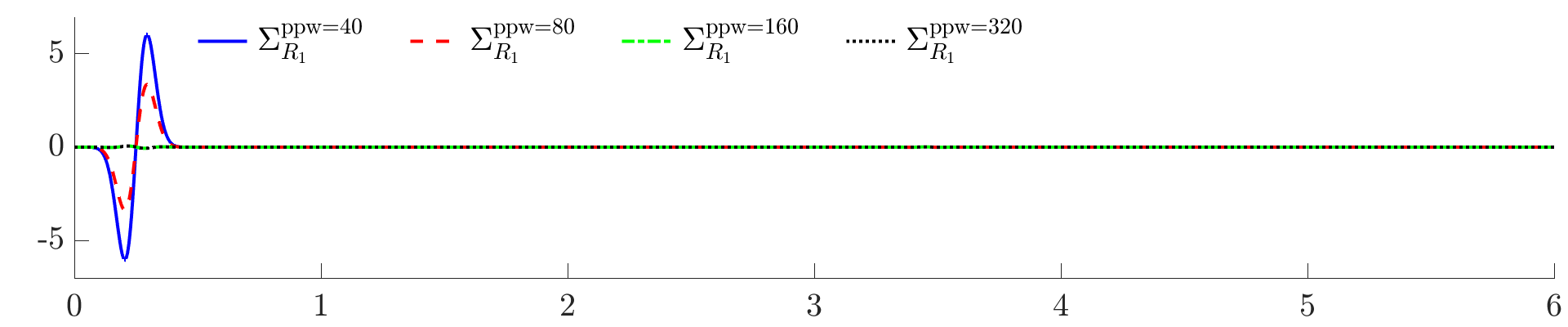}
\caption{Time histories of $\Sigma$ at $R_1$.}
\label{supp_P_1_6_s_extrapolating_40_to_320}
\end{subfigure}\hfill
\\[2ex]
\begin{subfigure}[b]{1\textwidth}
\captionsetup{width=1\textwidth, font=small,labelfont=small}
\centering\includegraphics[scale=0.15]{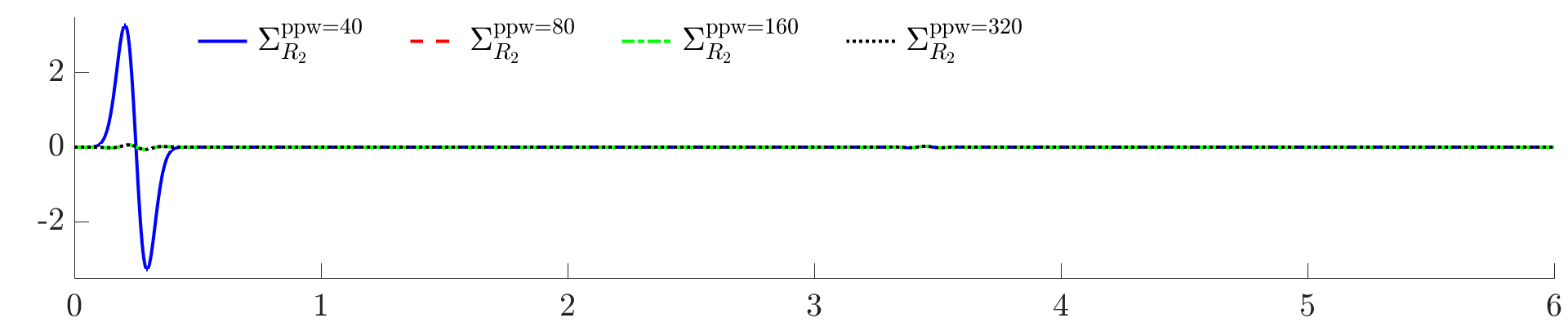}
\caption{Time histories of $\Sigma$ at $R_2$.}
\label{supp_P_2_6_s_extrapolating_40_to_320}
\end{subfigure}\hfill
\caption{Time histories of the solution variable $\Sigma$ at the three receiver locations depicted in Figure \ref{Figure_grid_layout_source_receiver} with the four grid lines now representing the cases of 40, 80, 160, 320 ppw, respectively. 
}
\label{supp_P_extrapolating_weak_40_to_320}
\end{figure}

An additional experiment is presented here to supplement the discussion in section \ref{subsection_1D_numerical_experiments}, Figure \ref{P_extrapolating_weak} in particular.
The experimental setup is identical to that has been illustrated in Figure \ref{Figure_grid_layout_source_receiver}, except that the ppw associated with the four grid lines, from top to bottom, are 40, 80, 160, and 320, respectively.
In other words, the grid spacing is a quarter of that used in the experiment from section \ref{subsection_1D_numerical_experiments}.
The time step length is specified as 5e-5 s for all four simulations, and the simulation time is specified as 6 s, which amounts to 120000 time steps.
Time histories of the solution variable $\Sigma$ at the three receiver locations are displayed in Figure \ref{supp_P_extrapolating_weak_40_to_320}, from which we observe that the violation of the weakly imposed free surface boundary condition persists when the source is placed too close to the free surface, despite the refined grid resolution.

\newpage
\section{Additional 1D experiment with the source placed far away from the free surface}\label{supp_section_interior_source}

\begin{figure}[H]
\captionsetup{width=1\textwidth, font=small,labelfont=small}
\centering
\begin{subfigure}[b]{1\textwidth}
\captionsetup{width=1\textwidth, font=small,labelfont=small}
\centering\includegraphics[scale=0.15]{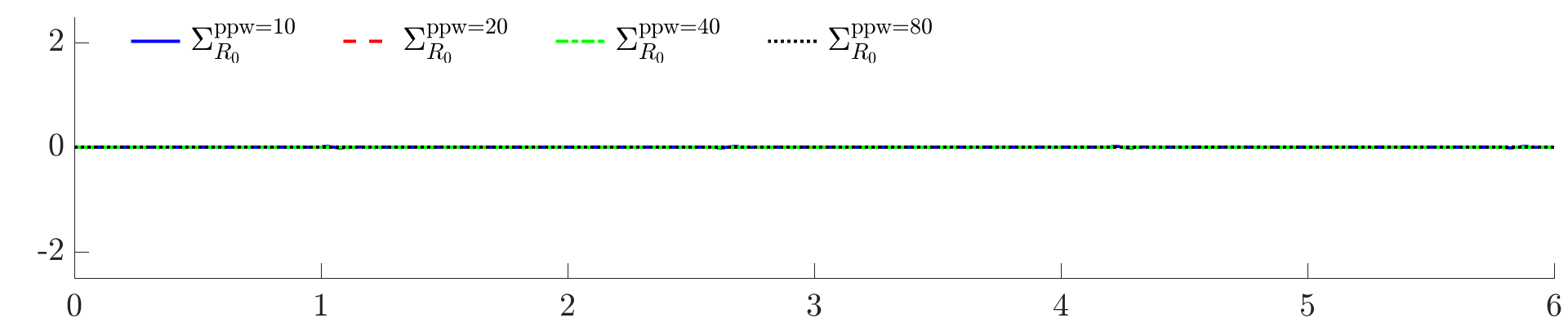}
\caption{Time histories of $\Sigma$ at $R_0$.} 
\label{supp_P_0_6_s_extrapolating_10_to_80_interior}
\end{subfigure}\hfill
\\[2ex]
\begin{subfigure}[b]{1\textwidth}
\captionsetup{width=1\textwidth, font=small,labelfont=small}
\centering\includegraphics[scale=0.15]{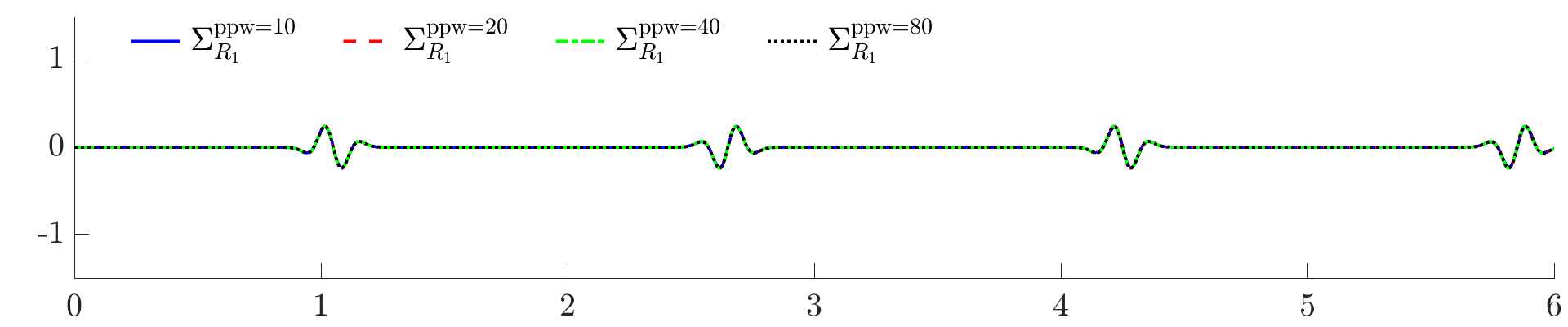}
\caption{Time histories of $\Sigma$ at $R_1$.}
\label{supp_P_1_6_s_extrapolating_10_to_80_interior}
\end{subfigure}\hfill
\\[2ex]
\begin{subfigure}[b]{1\textwidth}
\captionsetup{width=1\textwidth, font=small,labelfont=small}
\centering\includegraphics[scale=0.15]{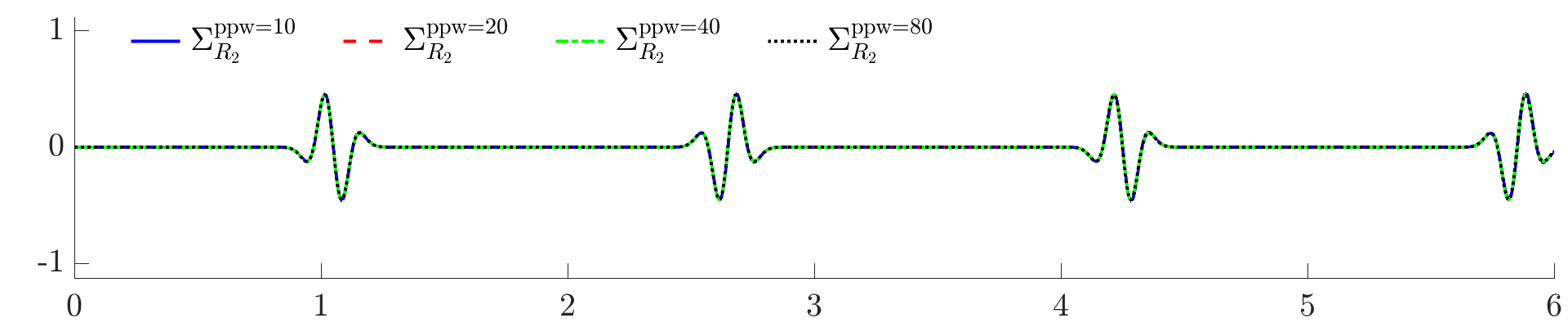}
\caption{Time histories of $\Sigma$ at $R_2$.}
\label{supp_P_2_6_s_extrapolating_10_to_80_interior}
\end{subfigure}\hfill
\caption{Time histories of the solution variable $\Sigma$ at the three receiver locations depicted in Figure \ref{Figure_grid_layout_source_receiver} with the source now placed far ($10$ times ppw grid points) away from the free surface.
}
\label{supp_P_extrapolating_weak_10_to_80_interior}
\end{figure}

An additional experiment is presented here to supplement the discussion in section \ref{subsection_1D_numerical_experiments}, Figure \ref{P_extrapolating_weak} in particular.
The experimental setup is identical to that has been illustrated in Figure \ref{Figure_grid_layout_source_receiver}, except that the source is placed at the middle of the simulation interval (10 minimal wavelength away from both boundaries, which correspond to 100, 200, 400, and 800 grid points for the four simulations, respectively). 
Time histories of the solution variable $\Sigma$ at the three receiver locations are displayed in Figure \ref{supp_P_extrapolating_weak_10_to_80_interior}, from which we observe that there is no severe violation of the free surface boundary condition 
%
and that the four simulation results agree much better (cf. Figure \ref{P_extrapolating_weak}).

\newpage
\section{Cross comparison}
\label{supp_cross_comparison}

Simulation results from three different discretizations are compared against each other to validate the implementation, which includes the finite difference discretization described in section \ref{subsection_1D_background} with weakly imposed free surface boundary condition, the finite difference discretization described in section \ref{subsection_1D_strong_imposition} with strongly imposed free surface boundary condition, and a finite element discretization.

The finite element discretization is based on a different formulation of the wave system from \eqref{1D_wave_equation}. 
Essentially, we reduce \eqref{1D_wave_equation} to a single equation by substituting \eqref{1D_wave_equation_a} into the time derivative of \eqref{1D_wave_equation_b}, leading to an equation that involves only second order temporal and spatial derivatives of $\sigma$. 
Starting with a different formulation gives better confidence in the validation. 
In space, the finite element simulation uses second-order Lagrange polynomials; in time, it uses second-order central difference.

\begin{figure}[H]
\captionsetup{width=1\textwidth, font=small,labelfont=small}
\centering
\begin{subfigure}[b]{1\textwidth}
\captionsetup{width=1\textwidth, font=small,labelfont=small}
\centering\includegraphics[scale=0.15]{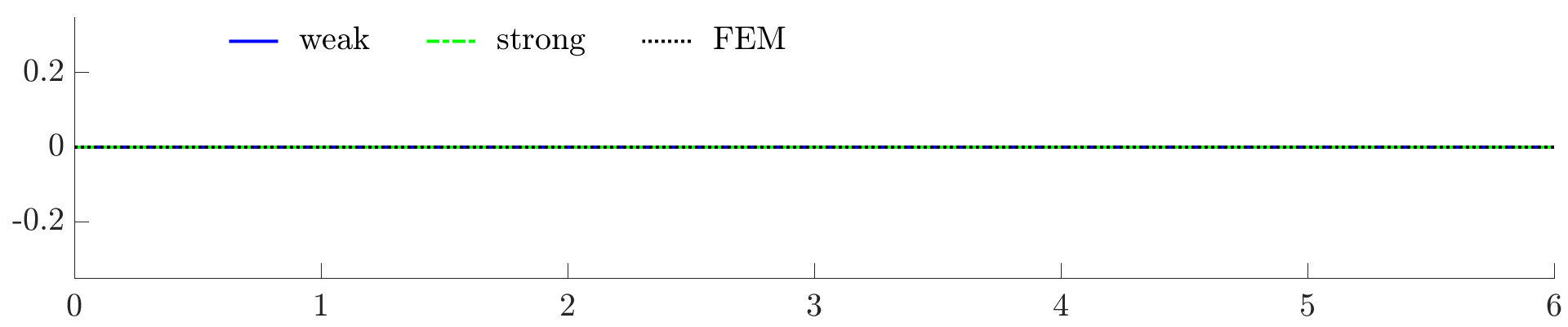}
\caption{Time histories of $\Sigma$ at $R_0$.} 
\label{supp_P_0_6_s_compare}
\end{subfigure}\hfill
\\[2ex]
\begin{subfigure}[b]{1\textwidth}
\captionsetup{width=1\textwidth, font=small,labelfont=small}
\centering\includegraphics[scale=0.15]{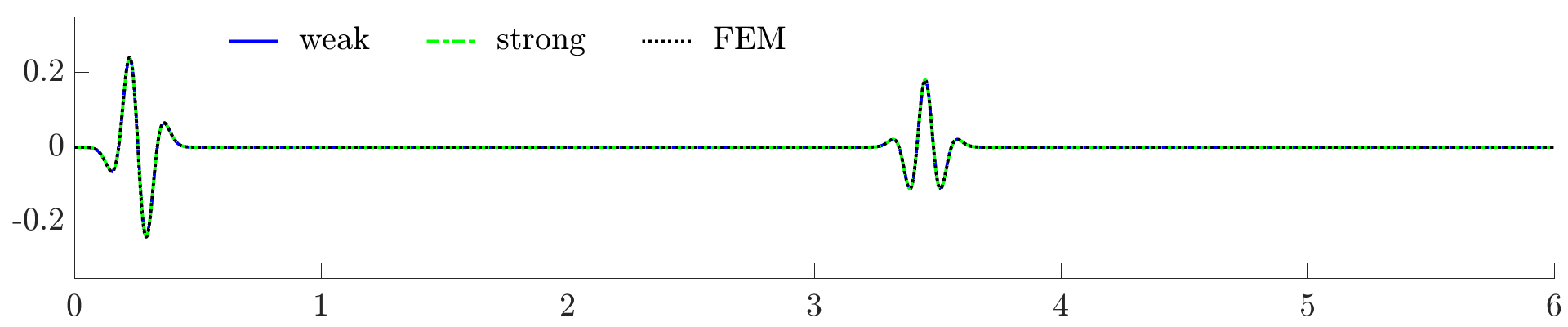}
\caption{Time histories of $\Sigma$ at $R_1$.}
\label{supp_P_1_6_s_compare}
\end{subfigure}\hfill
\\[2ex]
\begin{subfigure}[b]{1\textwidth}
\captionsetup{width=1\textwidth, font=small,labelfont=small}
\centering\includegraphics[scale=0.15]{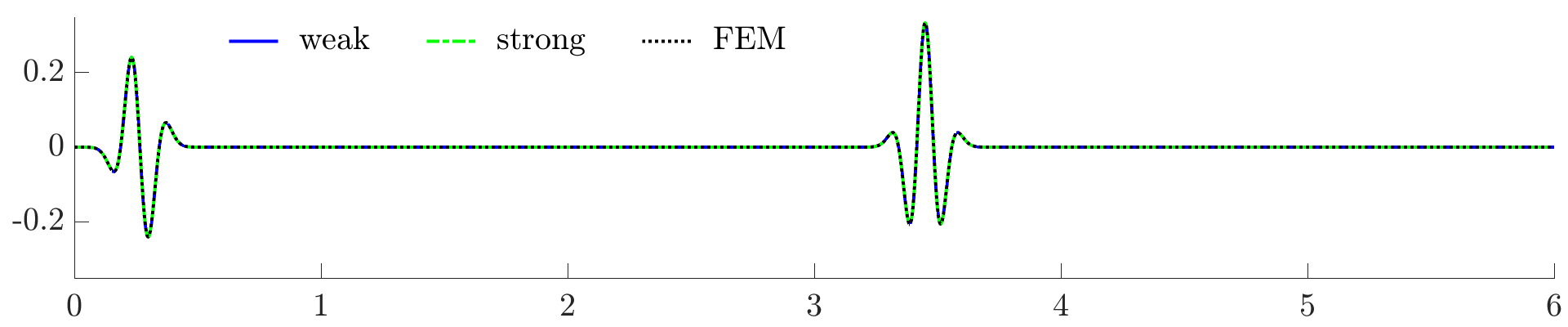}
\caption{Time histories of $\Sigma$ at $R_2$.}
\label{supp_P_2_6_s_compare}
\end{subfigure}\hfill
\caption{Time histories of the solution variable $\Sigma$ at the three receiver locations depicted in Figure \ref{Figure_grid_layout_source_receiver} from three different discretizations, including the finite difference discretization described in section \ref{subsection_1D_background} (denoted as weak), the finite difference discretization described in section \ref{subsection_1D_strong_imposition} (denoted as strong), and a finite element discretization (denoted as FEM). All three simulations use $\text{ppw}=80$.
}
\label{supp_P_compare}
\end{figure}

For clarity, we only compare the simulation results with $\text{ppw}=80$ for all three discretizations. The grid spacing and time step length are the same for all three simulations.
Time histories of the solution variable $\Sigma$ at the three receiver locations depicted in Figure \ref{Figure_grid_layout_source_receiver} are displayed in Figure \ref{supp_P_compare}, from where we observe that results from all three simulations agree well with each other, hence providing confidence in our implementations used for the experiments.

\newpage
\section{Spectral radii}
\label{supp_spectral_radii}
In the following, we consider the 1D wave system \eqref{1D_wave_equation} defined on interval $[0,1]$ with unit density and wave-speed; spectral radii of the semi-discretized systems described in sections \ref{subsection_1D_background} and \ref{subsection_1D_strong_imposition} are evaluated numerically to confirm the claims on CFL restrictions made therein. 

Both semi-discretized systems are represented using the following generic notations:
\begin{equation}
\label{generalized_eigenvalue_system_large}
\frac{d }{d t} 
\left[ 
\begin{array}{c c}
\mathcal A^M & \boldsymbol 0 \\
\boldsymbol 0 & \mathcal A^N 
\end{array}
\right]
\left[ 
\begin{array}{c}
V \\
\Sigma
\end{array}
\right]
\enskip = \enskip
\left[ 
\begin{array}{c c}
\boldsymbol 0 & \mathcal A^M \widetilde{\mathcal D}^N \\
\mathcal A^N \widetilde{\mathcal D}^M & \boldsymbol 0
\end{array}
\right]
\left[ 
\begin{array}{c}
V \\
\Sigma
\end{array}
\right],
\end{equation}
where $\widetilde{\mathcal D}^M$ and $\widetilde{\mathcal D}^N$ are the modified finite difference operators, either by appending SATs to \eqref{1D_semi_discretized_extrapolating_wo_SATs_b} as described in section \ref{subsection_1D_background} (the weak case), or by setting the specific rows and columns to zero as described in section \ref{subsection_1D_strong_imposition} (the strong case). 
In both cases, $\widetilde{\mathcal D}^M$ and $\widetilde{\mathcal D}^N$ satisfy the following relation:
\begin{equation}
\label{relation_modified_operators}
\mathcal A^N \widetilde{\mathcal D}^M \ + \ \left( \mathcal A^M \widetilde{\mathcal D}^N \right)^T \ = \ \boldsymbol 0
\end{equation}
by construction.
The generalized eigenvalues associated with the above system, together with the chosen time integration scheme, determine the stability requirement on time step length.

We note here that since the two matrices from \eqref{generalized_eigenvalue_system_large} are diagonal and skew-symmetric, all generalized eigenvalues are on the imaginary axis.
Moreover, because of the skew-symmetric nature, the generalized eigenvalues calculated based on \eqref{generalized_eigenvalue_system_large} tend to be more reliable than those from the system where the norm matrices $\mathcal A^M$ and $\mathcal A^N$ are omitted on both sides.

A related note is that system \eqref{generalized_eigenvalue_system_large} can be further reduced to 
\begin{equation}
\label{generalized_eigenvalue_system_small}
\mathcal A^M \frac{d^2 V}{d t^2} 
\enskip = \enskip 
- \left( \widetilde{\mathcal D}^M \right)^T \!\! \mathcal A^N \widetilde{\mathcal D}^M V \, ,
\end{equation}
by differentiating with respect to $t$, eliminating $\Sigma$ by substitution, and invoking \eqref{relation_modified_operators}.
It can be shown that the generalized eigenvalues from \eqref{generalized_eigenvalue_system_large} are the square roots of the generalized eigenvalues from \eqref{generalized_eigenvalue_system_small}.
Since \eqref{generalized_eigenvalue_system_small} is smaller in size and involves only symmetric matrices, the associated generalized eigenvalue problem is easier to solve numerically with more reliable outcome.
%

The spectral radii (i.e., maximal magnitude of the generalized eigenvalues) are shown in Table \ref{Table_spectral_radii} for $\Delta x = \nicefrac{1}{40}$,  $\Delta x = \nicefrac{1}{160}$, and $\Delta x = \nicefrac{1}{640}$. 
%
%
All calculated generalized eigenvalues associated with \eqref{generalized_eigenvalue_system_large} are indeed on the imaginary axis. 
For staggered leapfrog scheme, its associated stability region is $[-2\imath, 2\imath]$; see \cite{ghrist2000staggered}.
In this case, the spectral radii dictates the restriction on time step length. 
The spectral radii scaled by $C_\text{max} \cdot \Delta x$ (i.e., the maximally allowed $\Delta t$ according to $C_\text{max}$) are presented in Table \ref{Table_spectral_radii_scaled}, which need to be less than or equal to 2 for stable simulations with the staggered leapfrog scheme. 
The values in Table \ref{Table_spectral_radii_scaled} are right beneath 2, which confirms the claims made in sections \ref{subsection_1D_background} and \ref{subsection_1D_strong_imposition} on time step restrictions.


\begin{table}[H]
\small
\captionsetup{width=1\textwidth, font=small,labelfont=small}
\setlength{\tabcolsep}{3.5mm}
\centering
\begin{tabular}{ | l | c c c | }
\hline
   \multicolumn{1}{|c|}{spectral radii}      				& $\Delta x = \nicefrac{1}{40}$ & $\Delta x = \nicefrac{1}{160}$ & $\Delta x = \nicefrac{1}{640}$ \\[0.75mm] \hline
Extrapolating (weak)   & $125.871385897805$ & $503.485543591221$ & $2013.942174364883$ \\[0mm]
Extrapolating (strong) & $93.256016056512$  & $373.311280516844$ & $1493.327619757113$ \\[0.75mm]
\hline
Periodic 			   & $93.333333333333$  & $373.333333333333$ & $1493.333333333333$ \\[0.75mm]
\hline
\end{tabular}
\caption{Spectral radii corresponding to the semi-discretized systems described in sections \ref{subsection_1D_background} (weak) and \ref{subsection_1D_strong_imposition} (strong).
The spectral radii associated with periodic boundary condition are also presented as references, which amounts to applying the interior stencil on unbounded domain. 
} 
\label{Table_spectral_radii}
\end{table}

\begin{table}[H]
\small
\captionsetup{width=1\textwidth, font=small,labelfont=small}
\setlength{\tabcolsep}{3.5mm}
\centering
\begin{tabular}{ | l | c c c | }
\hline
\multicolumn{1}{|c|}{spectral radii (scaled)}           				& $\Delta x = \nicefrac{1}{40}$ & $\Delta x = \nicefrac{1}{160}$ & $\Delta x = \nicefrac{1}{640}$ \\[0.75mm] \hline
Extrapolating (weak)   & $1.999781643451$ & $1.999781643451$ & $1.999781643451$ \\[0mm]
Extrapolating (strong) & $1.998343201211$ & $1.999881859912$ & $1.999992347889$ \\[0.75mm]
\hline
Periodic 			   & $2.000000000000$ & $2.000000000000$ & $2.000000000000$ \\[0.75mm]
\hline
\end{tabular}
\caption{Spectral radii presented in Table \ref{Table_spectral_radii} scaled by 
$C_\text{max} \cdot \Delta x$, where $C_\text{max}$ is the numerically measured maximal Courant number. Such $C_\text{max}$ are $0.6355$, $\nicefrac{6}{7}$, and $\nicefrac{6}{7}$, as mentioned in \eqref{Interior_CFL_limit}, \eqref{1D_extrapolating_CFL_limit_weak}, and \eqref{1D_extrapolating_CFL_limit_strong} of section \ref{section_1D_case}.
} 
\label{Table_spectral_radii_scaled}
\end{table}

\newpage
\section{Acoustic wave equation and the semi-discretized system}
\label{supp_acoustic_background}

Below, we present the acoustic wave equation and its semi-discretization as the background information omitted in the main text, section \ref{section_2D_acoustic_case}.

The acoustic wave equation can be expressed as 
\begin{subequations}
\label{supp_2D_acoustic_wave_equation}
\begin{empheq}[left=\empheqlbrace]{alignat = 2}
\displaystyle \enskip \rho \frac{\partial v_x}{\partial t} \enskip &= \enskip \displaystyle \! \frac{\partial \sigma}{\partial x}; \label{supp_2D_acoustic_wave_equation_vx} \\
\displaystyle \enskip \rho \frac{\partial v_y}{\partial t} \enskip &= \enskip \displaystyle \! \frac{\partial \sigma}{\partial y}; \label{supp_2D_acoustic_wave_equation_vy} \\
\displaystyle \enskip \beta \frac{\partial \sigma}{\partial t} \enskip &= \enskip \displaystyle \! \frac{\partial v_x}{\partial x} + \frac{\partial v_y}{\partial y}, \label{supp_2D_acoustic_wave_equation_sigma}
\end{empheq}
\end{subequations}
where $\rho$ and $\beta$ are the physical parameters, representing density and compressibility, respectively, $v_x$ and $v_y$ are the particle velocities, $\sigma$ is the negative of pressure. 
We note here that for consistency within this study, we use $\sigma$, i.e., the negative of pressure, instead of pressure itself, so that all wave equations under discussion will have the same sign in front of the spatial derivatives. 
Moreover, $\beta = \nicefrac{1}{\rho c^2}$, where $c$ is wave-speed.

The semi-discretization of \eqref{supp_2D_acoustic_wave_equation} without concerning the boundary conditions can be expressed as 
\begin{subequations}
\label{Semi_discretized_acoustic_wave_equation_2D}
\begin{empheq}[left=\empheqlbrace]{alignat = 2}
\displaystyle \mathcal A^{V_x} \boldsymbol{\rho}^{V_x} \frac{d V_x}{d t} &\enskip = \enskip \displaystyle \mathcal A^{V_x} \mathcal D^{\Sigma}_x \, \Sigma \, ; 
\label{Semi_discretized_acoustic_wave_equation_2D_Vx} \\
\displaystyle \mathcal A^{V_y} \boldsymbol{\rho}^{V_y} \frac{d V_y}{d t} &\enskip = \enskip \displaystyle \mathcal A^{V_y} \mathcal D^{\Sigma}_y \, \Sigma \, ;
\label{Semi_discretized_acoustic_wave_equation_2D_Vy} \\
\displaystyle \mathcal A^{\Sigma} \boldsymbol{\beta}^{\Sigma} \frac{d \Sigma}{d t} &\enskip = \enskip \displaystyle \mathcal A^{\Sigma} \mathcal D^{V_x}_x V_x + \mathcal A^{\Sigma} \mathcal D^{V_y}_y V_y \, . 
\label{Semi_discretized_acoustic_wave_equation_2D_Sigma} 
\end{empheq}
\end{subequations}
Following the standard procedure, the 2D SBP operators appearing in \eqref{Semi_discretized_acoustic_wave_equation_2D} are constructed from their 1D counterparts via tensor product as

\begin{subequations}
\label{supp_2D_acoustic_SBP_operators}
\begin{alignat}{1}
\omit\hfill
$
\mathcal A^{P}   \! = \mathcal A^N_x \otimes \mathcal A^N_y, \enskip
\mathcal A^{V_x} \! = \mathcal A^M_x \otimes \mathcal A^N_y, \enskip 
\mathcal A^{V_y} \! = \mathcal A^N_x \otimes \mathcal A^M_y; 
$
\hfill
\\
\omit\hfill
$
\mathcal D_{x}^{P}   \! = \mathcal D^N_x \otimes \mathcal I^N_y, \enskip 
\mathcal D_{x}^{V_x} \! = \mathcal D^M_x \otimes \mathcal I^N_y, \enskip 
\mathcal D_{y}^{P}   \! = \mathcal I^N_x \otimes \mathcal D^N_y, \enskip
\mathcal D_{y}^{V_y} \! = \mathcal I^N_x \otimes \mathcal D^M_y.
$
\hfill
\end{alignat}
\end{subequations}

For the numerical experiment considered in section \ref{section_2D_acoustic_case} of the main text, where the free surface boundary condition is associated with all boundaries, the above difference operators need to be modified to account for these boundary conditions.

For the strong case, these modifications involve setting the corresponding rows and columns in the 1D operators to zero as follows:
\begin{itemize}
\setlength\itemsep{0.em}
\item[] setting the first and last rows of $D_x^M$ to zero;
\item[] setting the first and last rows of $D_y^M$ to zero;
\item[] setting the first and last columns of $D_x^N$ to zero;
\item[] setting the first and last columns of $D_y^N$ to zero;
\item[] setting the first and last entries of $I_x^N$ to zero;
\item[] setting the first and last entries of $I_y^N$ to zero.
\end{itemize}

For the weak case, these modifications involve adjusting for the penalty terms as follows:
\begin{subequations}
\label{acoustic_weak_modifications}
\begin{alignat}{3}
\mathcal D_{x}^{P} \enspace \longleftarrow \enspace \mathcal D_{x}^{P} 
& + 
\left[ \left( \left( \mathcal A^M_x \right)^{-1} \mathcal E^L_x \right) \otimes \mathcal I^N_y \right] 
\left[ \left( \mathcal P^L_x \right)^{T} \otimes \mathcal I^N_y - \boldsymbol 0 \right] 
\label{acoustic_weak_modifications_a}
\\
& -
\left[ \left( \left( \mathcal A^M_x \right)^{-1} \mathcal E^R_x \right) \otimes \mathcal I^N_y \right] 
\left[ \left( \mathcal P^R_x \right)^{T} \otimes \mathcal I^N_y  - \boldsymbol 0 \right];
\label{acoustic_weak_modifications_b}
\\
\mathcal D_{y}^{P} \enspace \longleftarrow \enspace \mathcal D_{y}^{P} 
& + 
\left[ \mathcal I^N_x \otimes \left( \left( \mathcal A^M_y \right)^{-1} \mathcal E^L_y \right) \right] 
\left[ \mathcal I^N_x \otimes \left( \mathcal P^L_y \right)^{T}  - \boldsymbol 0 \right] 
\label{acoustic_weak_modifications_c}
\\
& -
\left[ \mathcal I^N_x \otimes \left( \left( \mathcal A^M_y \right)^{-1} \mathcal E^R_y \right) \right] 
\left[ \mathcal I^N_x \otimes \left( \mathcal P^R_y \right)^{T}  - \boldsymbol 0 \right],
\label{acoustic_weak_modifications_d}
\end{alignat}
\end{subequations}
where \eqref{acoustic_weak_modifications_a} - \eqref{acoustic_weak_modifications_d} accounts for the left, right, bottom, and top boundaries, respectively.

\newpage
\section{Additional figures for the 2D acoustic experiments}
\label{supp_additional_figures_2D_acoustic}

Below, we present the time histories of the solution variables $V_x$ and $V_y$ for the 2D acoustic experiment from section \ref{section_2D_acoustic_case} of the main text.
For $V_x$, the strong and weak cases are presented in Figures \ref{Acoustic_full_FSBC_strong_Vx} and \ref{Acoustic_full_FSBC_weak_Vx}, respectively.
For $V_y$, the strong and weak cases are presented in Figures \ref{Acoustic_full_FSBC_strong_Vy} and \ref{Acoustic_full_FSBC_weak_Vy}, respectively.

Because of grid staggering, the receiver locations for $V_x$ is shifted rightwards for half grid spacing from those for $\Sigma$ in the main text; the receiver locations for $V_y$ is shifted downwards for half grid spacing from those for $\Sigma$ in the main text. 
We note here that the first receiver of $V_x$ is on the surface, which is supposed to remain zero.

\input{input_supp_acoustic}

\newpage
\section{1D Convergence tests}\label{supp_convergence_tests}
In the following, 1D convergence test results using a smooth manufactured solution are shown for the experiment settings presented in section \ref{subsection_1D_numerical_experiments} (weak), section \ref{subsection_1D_strong_imposition} (strong) of the main text, and for the operators presented in \eqref{SBP_matrices_1D_intertwined} of the main text applied on an $\mathbb M$ grid column of Figure \ref{Elastic_2D_grid_layout_full_FSBC_shaded} of the main text.
These operators applied on an $\mathbb N$ grid column delivers the same results as the strong case corresponding to section \ref{subsection_1D_strong_imposition} since after resetting, the operators become the same, and are omitted below.

Specifically, the manufactured solution reads:
\begin{equation}
\label{smooth_manufactured_solution}
\begin{array}{rcl}
p & \!\!\! = \!\!\! & \displaystyle \phantom{-}  \sin(8 \pi x) \sin(8 \pi t); \\
v & \!\!\! = \!\!\! & \displaystyle - \cos(8 \pi x) \cos(8 \pi t).
\end{array}
\end{equation}
where $x \in [0,1]$. 
%
%
For all simulations involved in the following tests, the time step length is fixed at 1e-6~s, while the total number of time steps is chosen as 666667, leading to a final time at approximately \nicefrac{2}{3}~s.
The error between the simulated solution and the manufactured solution at the final time is measured using the energy norm (i.e., the weighted $L_2$ norm with the norm matrices serving as the weights) and displayed below.

\begin{table}[H]
\small
\captionsetup{width=1\textwidth, font=small,labelfont=small}
\setlength{\tabcolsep}{3.5mm}
\centering
\begin{tabular}{ | c | c  c  c  c c | }
\hline
		& 10 ppw     & 20 ppw     & 40 ppw     & 80 ppw     & 160 ppw     \\ \hline
error  	& 2.5174e-02 & 2.1005e-03 & 1.6514e-04 & 1.4079e-05 & 1.2281e-06  \\
rate    & ---        & 3.5832     & 3.6689     & 3.5521     & 3.5191 	  \\
\hline
\end{tabular}
\caption{Error and rate for the weak case, corresponding to section \ref{subsection_1D_numerical_experiments} of the main text.} 
\label{Convergence_test_1D_weak}
\end{table}
\begin{table}[H]
\small
\captionsetup{width=1\textwidth, font=small,labelfont=small}
\setlength{\tabcolsep}{3.5mm}
\centering
\begin{tabular}{ | c | c  c  c  c c | }
\hline
   		& 10 ppw     & 20 ppw     & 40 ppw     & 80 ppw     & 160 ppw     \\ \hline
error   & 2.4014e-02 & 1.2913e-03 & 6.2780e-05 & 3.4777e-06 & 2.2958e-07  \\
rate    & ---        & 4.2169     & 4.3624     & 4.1741     & 3.9210 	  \\
\hline
\end{tabular}
\caption{Error and rate for the strong case, corresponding to section \ref{subsection_1D_strong_imposition} of the main text.} 
\label{Convergence_test_1D_strong}
\end{table}
\begin{table}[H]
\small
\captionsetup{width=1\textwidth, font=small,labelfont=small}
\setlength{\tabcolsep}{3.5mm}
\centering
\begin{tabular}{ | c | c  c  c  c c | }
\hline
   		& 10 ppw     & 20 ppw     & 40 ppw     & 80 ppw     & 160 ppw     \\ \hline
error   & 1.0505e-02 & 5.5674e-04 & 3.2667e-05 & 2.5932e-06 & 2.1248e-07  \\
rate    & ---        & 4.2379     & 4.0911     & 3.6550     & 3.6093	  \\
\hline
\end{tabular}
\caption{Error and rate for the strong case, corresponding to applying the operators presented in \eqref{SBP_matrices_1D_intertwined} of the main text on an $\mathbb M$ grid column of Figure \ref{Elastic_2D_grid_layout_full_FSBC_shaded} of the main text.}
\label{Convergence_test_1D_strong}
\end{table}

\newpage
\section{Additional figures for the 2D elastic experiments}
\label{supp_additional_figures_2D_elastic}

Below, we present the time histories of the solution variables $\Sigma_{xx}$, $\Sigma_{yy}$, $V_x$, and $V_y$ for the 2D elastic experiment from section \ref{section_numerical_example_elastic} of the main text, corresponding to Figures \ref{Elastic_full_FSBC_Sxy_strong} and \ref{Elastic_full_FSBC_Sxy_weak} therein, where a compressional source (i.e., applied on $\Sigma_{xx}$ and $\Sigma_{yy}$) is considered.

For $S_{xx}$, the strong and weak cases are presented in Figures \ref{Elastic_full_FSBC_Sxx_strong} and \ref{Elastic_full_FSBC_Sxx_weak}, respectively.
For $S_{yy}$, the strong and weak cases are presented in Figures \ref{Elastic_full_FSBC_Syy_strong} and \ref{Elastic_full_FSBC_Syy_weak}, respectively.
For $V_x$, the strong and weak cases are presented in Figures \ref{Elastic_full_FSBC_Vx_strong} and \ref{Elastic_full_FSBC_Vx_weak}, respectively.
For $V_y$, the strong and weak cases are presented in Figures \ref{Elastic_full_FSBC_Vy_strong} and \ref{Elastic_full_FSBC_Vy_weak}, respectively.

\input{input_source_on_Sxx_Syy}

\newpage
\section{Additional figures for the 2D elastic ``stress test'' with directional source}
\label{supp_additional_figures_2D_elastic_stress_test_Vy}

Below, we present the time histories of the solution variables $\Sigma_{xx}$, $\Sigma_{yy}$, $V_x$, and $V_y$ for the 2D elastic experiment from section \ref{section_numerical_example_elastic} of the main text, corresponding to Figures \ref{stress_test_Vy_elastic_full_FSBC_Sxy_strong} and \ref{stress_test_Vy_elastic_full_FSBC_Sxy_weak} therein, where a directional source (i.e., applied on $V_y$) is considered and placed at the corner.

For $S_{xx}$, the strong and weak cases are presented in Figures \ref{stress_test_Vy_elastic_full_FSBC_Sxx_strong} and \ref{stress_test_Vy_elastic_full_FSBC_Sxx_weak}, respectively.
For $S_{yy}$, the strong and weak cases are presented in Figures \ref{stress_test_Vy_elastic_full_FSBC_Syy_strong} and \ref{stress_test_Vy_elastic_full_FSBC_Syy_weak}, respectively.
For $V_x$, the strong and weak cases are presented in Figures \ref{stress_test_Vy_elastic_full_FSBC_Vx_strong} and \ref{stress_test_Vy_elastic_full_FSBC_Vx_weak}, respectively.
For $V_y$, the strong and weak cases are presented in Figures \ref{stress_test_Vy_elastic_full_FSBC_Vy_strong} and \ref{stress_test_Vy_elastic_full_FSBC_Vy_weak}, respectively.

\input{input_source_on_Vy_corner}

\newpage
\section{Another 2D elastic ``stress test'' with source applied on $\Sigma_{xy}$}
\label{supp_additional_figures_2D_elastic_stress_test_Sxy}

Below, we present the time histories of another ``stress test'' for the strong and weak approaches. 
Specifically, the source is applied on $\Sigma_{xy}$ (although not common, this can still happen in practice, e.g., when complex source pattern is being inverted) and placed at 1 grid spacing away from the top and left surfaces, i.e., the first interior $\Sigma_{xy}$ point at the top left corner.
The receiver is placed at 5 grid spacing away from the left surface and at 0, 1, and 2 grid spacing at depth for $\Sigma_{xy}$ and $V_y$, with an additional half grid downward shift for the other components due to grid staggering.

The recorded time histories are presented in Figures \ref{stress_test_Sxy_Elastic_full_FSBC_Sxx_strong} - \ref{stress_test_Sxy_Elastic_full_FSBC_Vy_weak} , in the order of $\Sigma_{xx}$, $\Sigma_{yy}$, $\Sigma_{xy}$, $V_x$, and $V_y$, from which we observe that the strong approach performs markedly better than the weak case for all five solution components.

\input{input_source_on_Sxy_corner}
\newpage
\section{2D convergence tests}\label{supp_2D_convergence_tests}
In the following, 
2D convergence test results using a smooth manufactured solution for the elastic wave system (see \eqref{2D_elastic_wave_system} of the main text) are shown for the strong (Table \ref{Convergence_test_2D_strong}) and weak (Table \ref{Convergence_test_2D_weak}) approaches to accompany the discussion and comparison made in section \ref{section_2D_elastic_case} of the main text.

Specifically, the physical domain is a square with unit length on both $x$- and $y$-directions.
The manufactured solution reads:
\begin{linenomath}
\begin{subequations}
\label{2D_elastic_manufactured_solution}
\begin{empheq}[left=\empheqlbrace]{alignat = 4}
\sigma_{xx} = \ &   & & 2 k \sin(w t) \sin(k x) \sin(k y);  \\
\sigma_{xy} = \ &   & & 0; 									\\
\sigma_{yy} = \ & - & & 2 k \sin(w t) \sin(k y) \sin(k x);  \\
v_x 		= \ & - & & w \cos(w t) \sin(k y) \cos(k x); 	\\
v_y 		= \ &   & & w \cos(w t) \cos(k y) \sin(k x),
\end{empheq}
\end{subequations}
\end{linenomath}
where $k$ is set to $2\pi$ and $w$ is set to $\sqrt{2}k$.

For all simulations involved in the following tests, the time step length is fixed at 1e-6~s, while the total number of time steps is chosen as 666667, leading to a final time at approximately \nicefrac{2}{3}~s.
The error between the simulated solution and the manufactured solution at the final time is measured using the energy norm, which is derived from \eqref{Discrete_energy_2D_elastic_wave_system} of the main text.

\begin{table}[H]
\small
\captionsetup{width=1\textwidth, font=small,labelfont=small}
\setlength{\tabcolsep}{3.5mm}
\centering
\begin{tabular}{ | c | c  c  c  c c | }
\hline
		& 10 ppw     & 20 ppw     & 40 ppw     & 80 ppw     & 160 ppw     \\ \hline
error  	& 1.7389e-01 & 1.1894e-02 & 7.9771e-04 & 6.3477e-05 & 5.4416e-06  \\
rate    & ---        & 3.8699     & 3.8982     & 3.6515     & 3.5441 	  \\
\hline
\end{tabular}
\caption{Error and rate for the strong approach} 
\label{Convergence_test_2D_strong}
\end{table}

\begin{table}[H]
\small
\captionsetup{width=1\textwidth, font=small,labelfont=small}
\setlength{\tabcolsep}{3.5mm}
\centering
\begin{tabular}{ | c | c  c  c  c c | }
\hline
		& 10 ppw     & 20 ppw     & 40 ppw     & 80 ppw     & 160 ppw     \\ \hline
error  	& 2.8617e-01 & 3.3790e-02 & 5.1094e-03 & 8.6026e-04 & 1.5148e-04  \\
rate    & ---        & 3.0822     & 2.7253     & 2.5703     & 2.5057 	  \\
\hline
\end{tabular}
\caption{Error and rate for the weak approach.} 
\label{Convergence_test_2D_weak}
\end{table}



%% file: input_supp_acoustic.tex
\hspace{-0.0625\textwidth}
\begin{minipage}[t]{.5\textwidth}
\begin{figure}[H]
\captionsetup{width=1\textwidth, font=footnotesize,labelfont=footnotesize}
\centering
\begin{subfigure}[b]{1\textwidth}
\captionsetup{width=1\textwidth, font=footnotesize,labelfont=footnotesize}
\centering\includegraphics[scale=0.1]{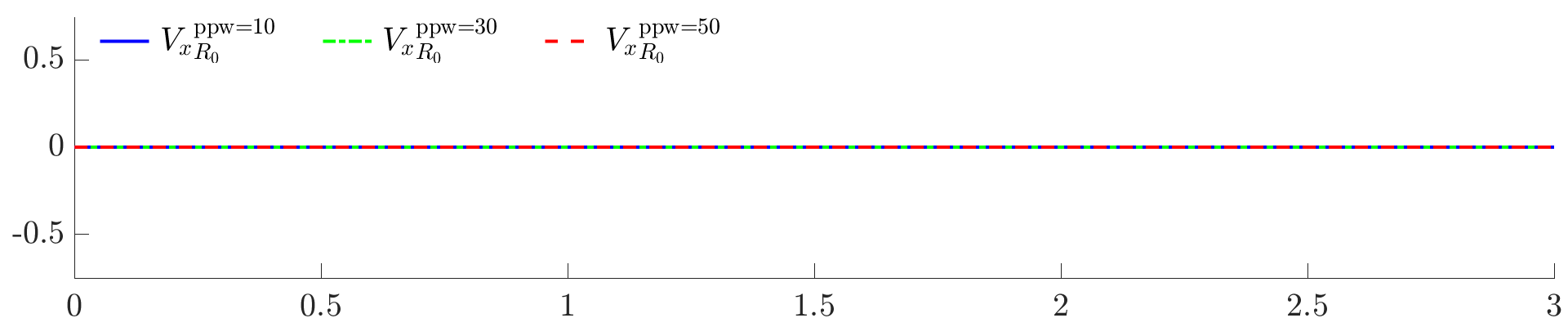}
\caption{Time history of $V_x$ at the surface.}
\label{Acoustic_full_FSBC_Vx_0_strong}
\end{subfigure}\hfill
\\[2ex]
\begin{subfigure}[b]{1\textwidth}
\captionsetup{width=1\textwidth, font=footnotesize,labelfont=footnotesize}
\centering\includegraphics[scale=0.1]{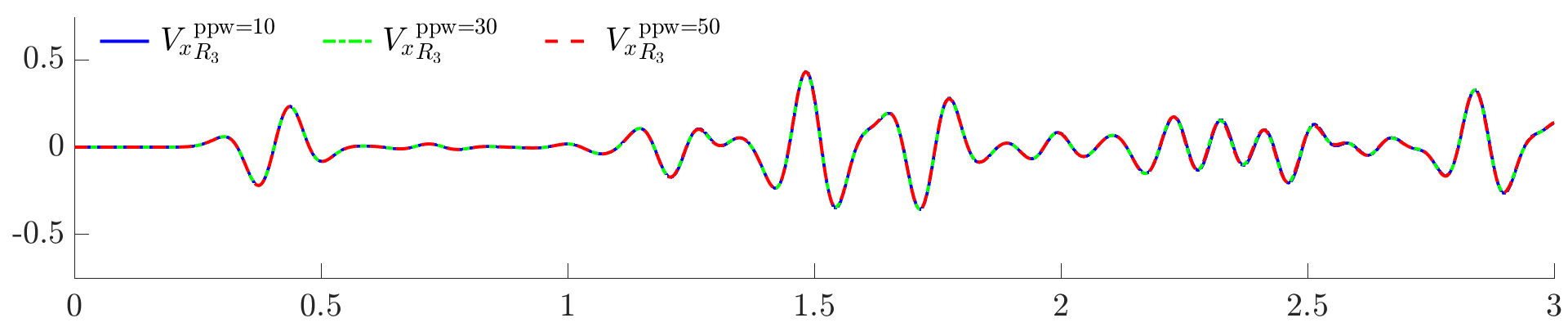}
\caption{Time history of $V_x$ at 3 grid points below the surface for the case ppw = 10.}
\label{Acoustic_full_FSBC_Vx_1_strong}
\end{subfigure}\hfill
\\[2ex]
\begin{subfigure}[b]{1\textwidth}
\captionsetup{width=1\textwidth, font=footnotesize,labelfont=footnotesize}
\centering\includegraphics[scale=0.1]{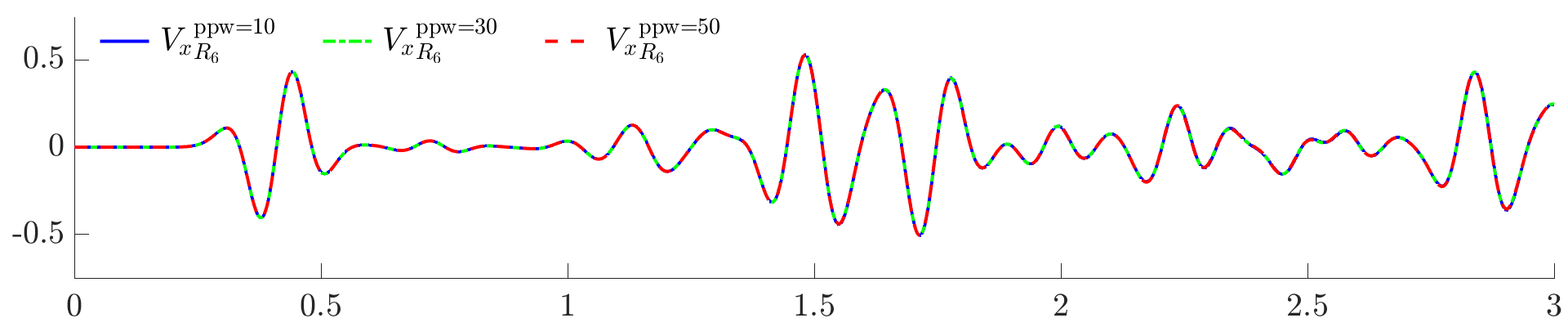}
\caption{Time history of $V_x$ at 6 grid points below the surface for the case ppw = 10.}
\label{Acoustic_full_FSBC_Vx_2_strong}
\end{subfigure}\hfill
\caption{Time histories of the 2D acoustic experiments ($V_x$ component). The free surface boundary condition is imposed strongly.
}
\label{Acoustic_full_FSBC_strong_Vx}
\end{figure}
\end{minipage}
\hspace{0.0075\textwidth}
\begin{minipage}[t]{.5\textwidth}
\begin{figure}[H]
\captionsetup{width=1\textwidth, font=footnotesize,labelfont=footnotesize}
\centering
\begin{subfigure}[b]{1\textwidth}
\captionsetup{width=1\textwidth, font=footnotesize,labelfont=footnotesize}
\centering\includegraphics[scale=0.1]{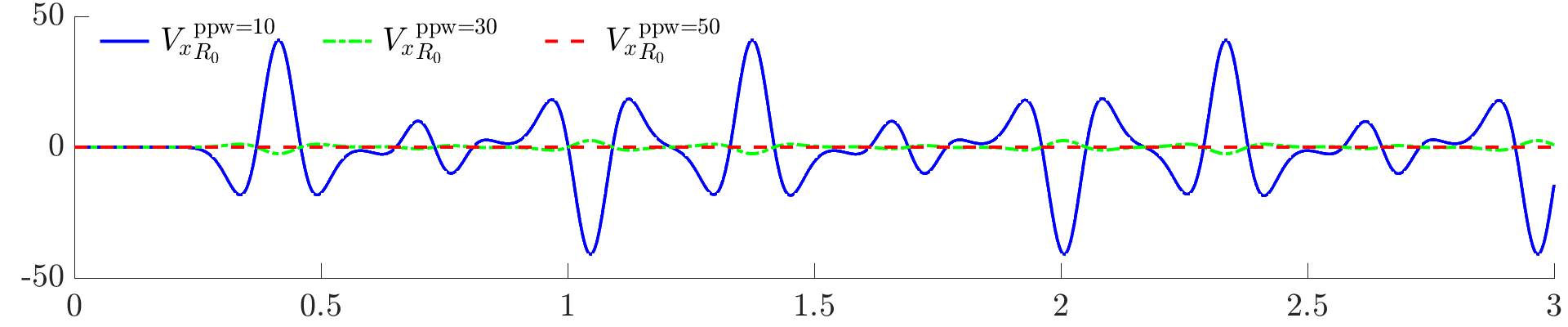}
\caption{Time history of $Vx$ at the surface.}
\label{Acoustic_full_FSBC_Vx_0_weak}
\end{subfigure}\hfill
\\[2ex]
\begin{subfigure}[b]{1\textwidth}
\captionsetup{width=1\textwidth, font=footnotesize,labelfont=footnotesize}
\centering\includegraphics[scale=0.1]{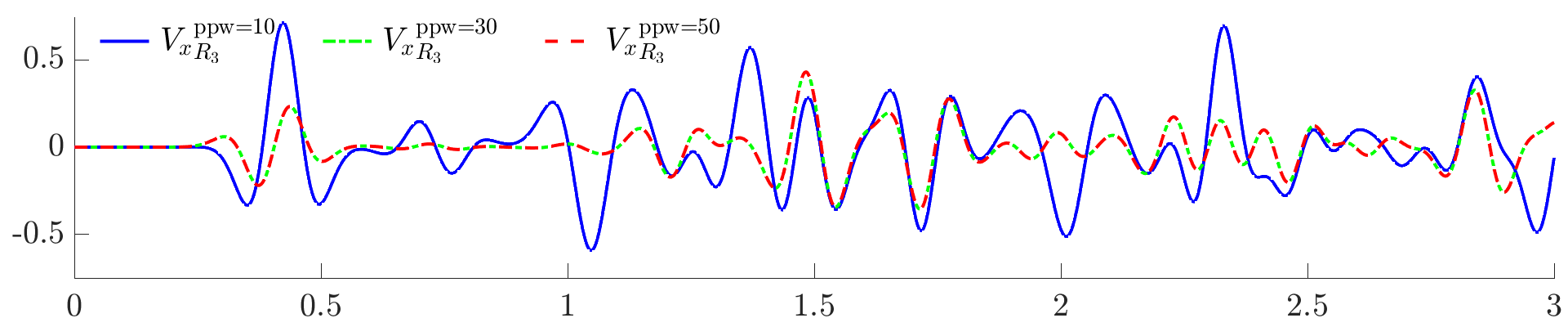}
\caption{Time history of $Vx$ at 3 grid points below the surface for the case ppw = 10.}
\label{Acoustic_full_FSBC_Vx_1_weak}
\end{subfigure}\hfill
\\[2ex]
\begin{subfigure}[b]{1\textwidth}
\captionsetup{width=1\textwidth, font=footnotesize,labelfont=footnotesize}
\centering\includegraphics[scale=0.1]{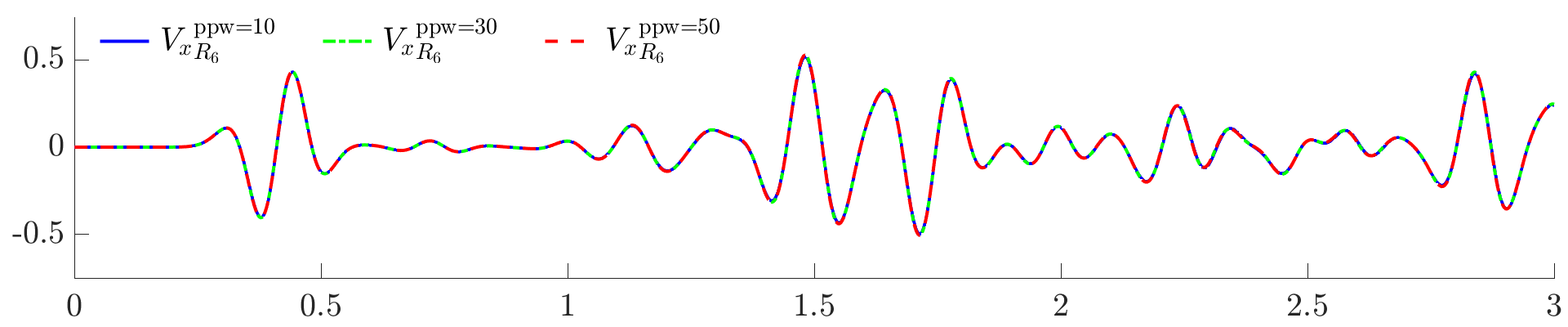}
\caption{Time history of $Vx$ at 6 grid points below the surface for the case ppw = 10.}
\label{Acoustic_full_FSBC_Vx_2_weak}
\end{subfigure}\hfill
\caption{Time histories of the 2D acoustic experiments ($V_x$ component). The free surface boundary condition is imposed weakly.
}
\label{Acoustic_full_FSBC_weak_Vx}
\end{figure}
\end{minipage}
\ \\
\ \\

\hspace{-0.0625\textwidth}
\begin{minipage}[t]{.5\textwidth}
\begin{figure}[H]
\captionsetup{width=1\textwidth, font=footnotesize,labelfont=footnotesize}
\centering
\begin{subfigure}[b]{1\textwidth}
\captionsetup{width=1\textwidth, font=footnotesize,labelfont=footnotesize}
\centering\includegraphics[scale=0.1]{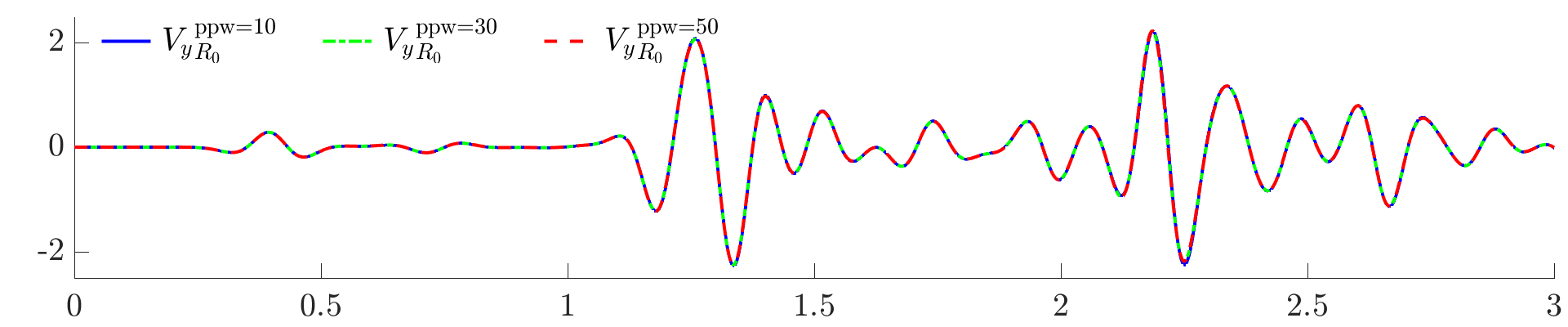}
\caption{Time history of $V_y$ at the surface.}
\label{Acoustic_full_FSBC_Vy_0_strong}
\end{subfigure}\hfill
\\[2ex]
\begin{subfigure}[b]{1\textwidth}
\captionsetup{width=1\textwidth, font=footnotesize,labelfont=footnotesize}
\centering\includegraphics[scale=0.1]{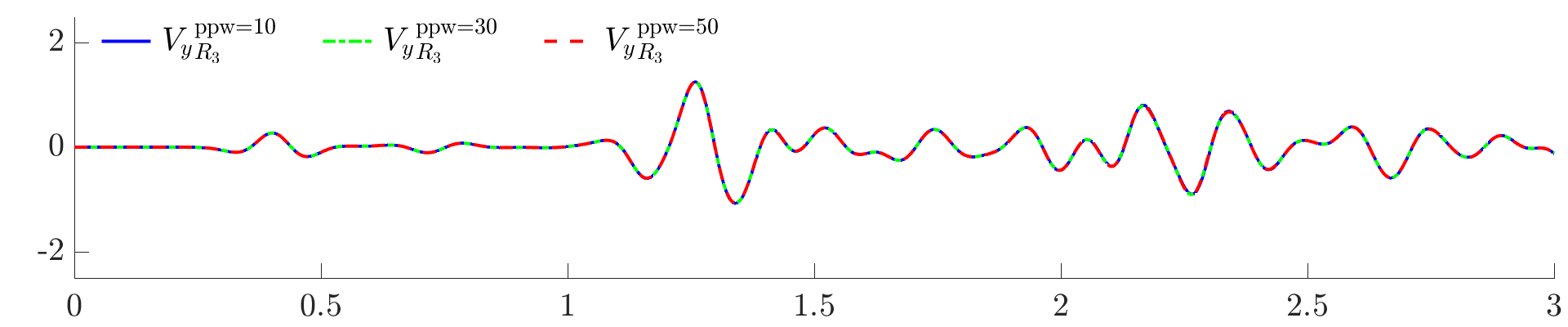}
\caption{Time history of $V_y$ at 3 grid points below the surface for the case ppw = 10.}
\label{Acoustic_full_FSBC_Vy_1_strong}
\end{subfigure}\hfill
\\[2ex]
\begin{subfigure}[b]{1\textwidth}
\captionsetup{width=1\textwidth, font=footnotesize,labelfont=footnotesize}
\centering\includegraphics[scale=0.1]{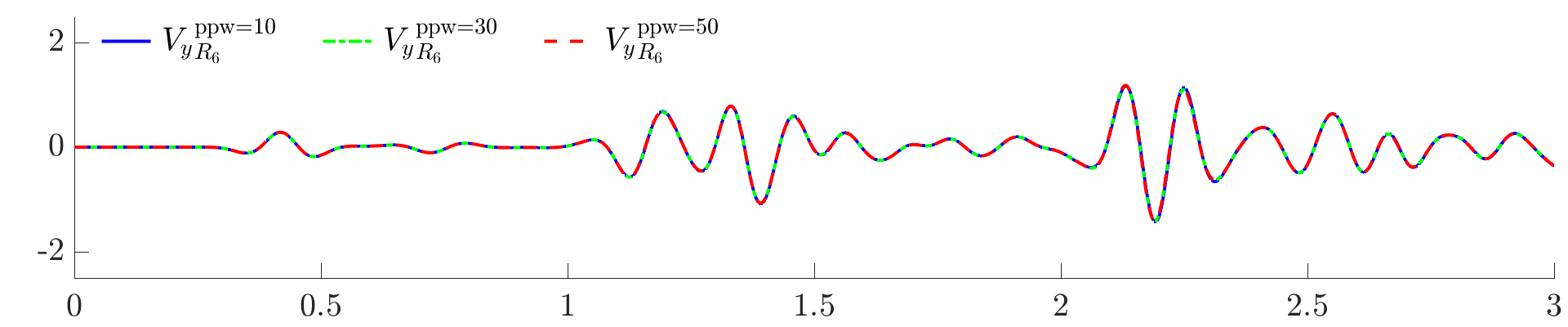}
\caption{Time history of $V_y$ at 6 grid points below the surface for the case ppw = 10.}
\label{Acoustic_full_FSBC_Vy_2_strong}
\end{subfigure}\hfill
\caption{Time histories of the 2D acoustic experiments ($V_x$ component). The free surface boundary condition is imposed strongly.
}
\label{Acoustic_full_FSBC_strong_Vy}
\end{figure}
\end{minipage}
\hspace{0.0075\textwidth}
\begin{minipage}[t]{.5\textwidth}
\begin{figure}[H]
\captionsetup{width=1\textwidth, font=footnotesize,labelfont=footnotesize}
\centering
\begin{subfigure}[b]{1\textwidth}
\captionsetup{width=1\textwidth, font=footnotesize,labelfont=footnotesize}
\centering\includegraphics[scale=0.1]{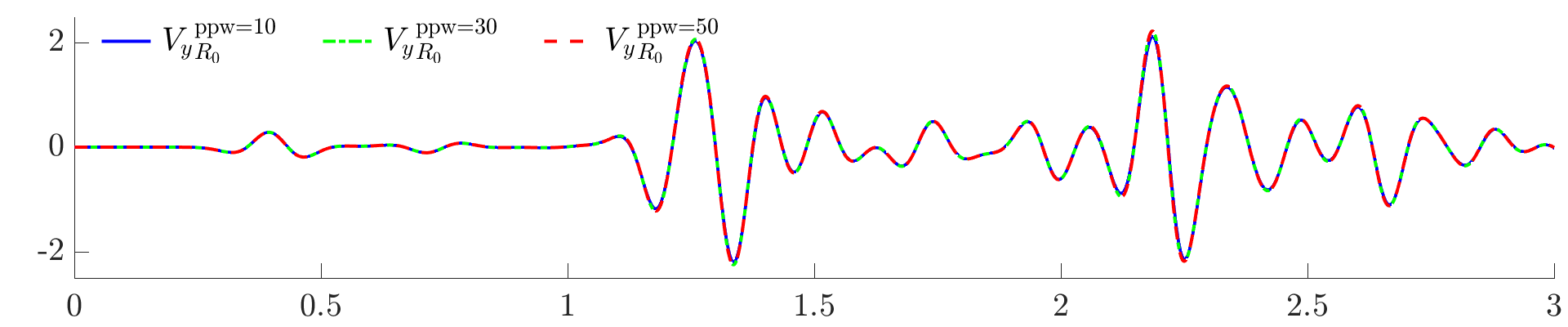}
\caption{Time history of $Vy$ at the surface.}
\label{Acoustic_full_FSBC_Vy_0_weak}
\end{subfigure}\hfill
\\[2ex]
\begin{subfigure}[b]{1\textwidth}
\captionsetup{width=1\textwidth, font=footnotesize,labelfont=footnotesize}
\centering\includegraphics[scale=0.1]{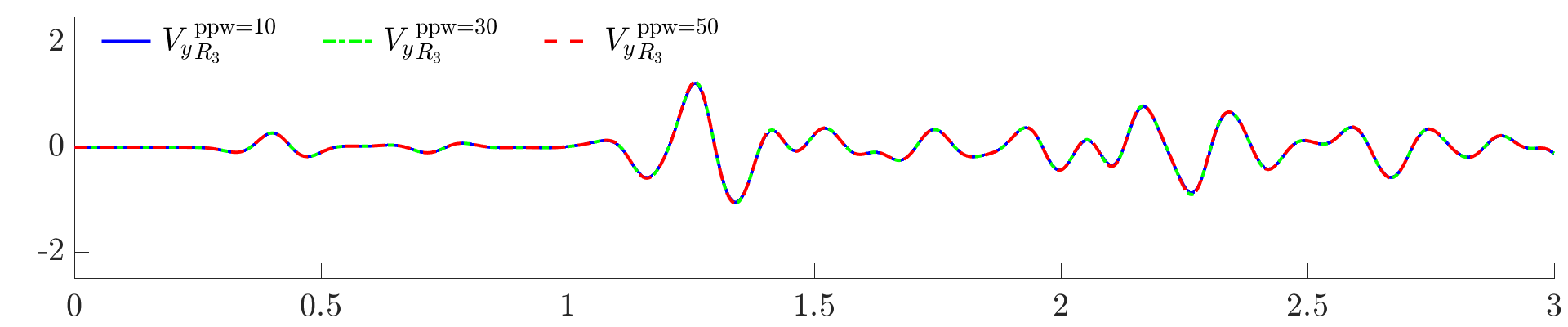}
\caption{Time history of $Vy$ at 3 grid points below the surface for the case ppw = 10.}
\label{Acoustic_full_FSBC_Vy_1_weak}
\end{subfigure}\hfill
\\[2ex]
\begin{subfigure}[b]{1\textwidth}
\captionsetup{width=1\textwidth, font=footnotesize,labelfont=footnotesize}
\centering\includegraphics[scale=0.1]{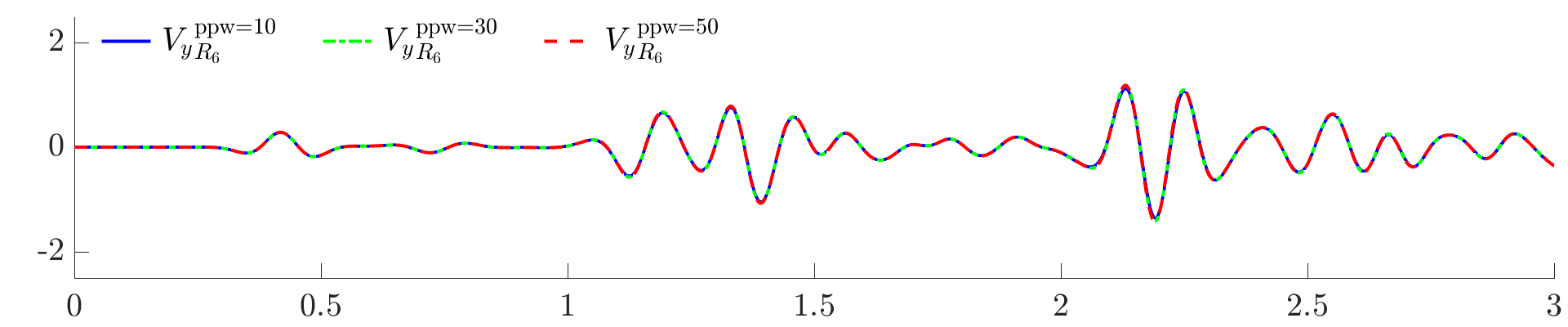}
\caption{Time history of $Vx$ at 6 grid points below the surface for the case ppw = 10.}
\label{Acoustic_full_FSBC_Vy_2_weak}
\end{subfigure}\hfill
\caption{Time histories of the 2D acoustic experiments ($V_y$ component). The free surface boundary condition is imposed weakly.
}
\label{Acoustic_full_FSBC_weak_Vy}
\end{figure}
\end{minipage}

%% file: input_source_on_Sxx_Syy.tex
\hspace{-0.0625\textwidth}
\begin{minipage}[t]{.5\textwidth}
\begin{figure}[H]
\captionsetup{width=1\textwidth, font=small,labelfont=small}
\centering
\begin{subfigure}[b]{1\textwidth}
\captionsetup{width=1\textwidth, font=small,labelfont=small}
\centering\includegraphics[scale=0.2]{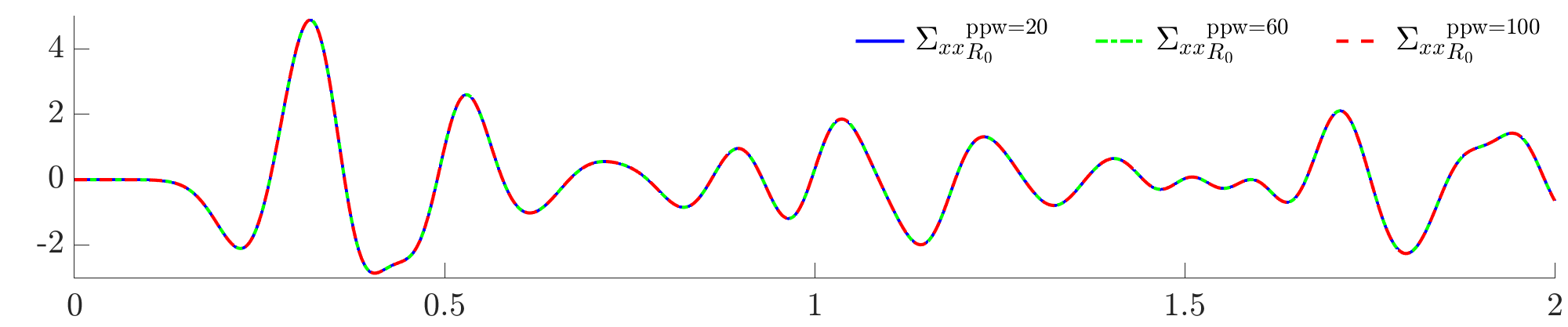}
\caption{Time history of $\Sigma_{xx}$ at the surface.}
\label{Elastic_full_FSBC_Sxx_0_strong}
\end{subfigure}\hfill
\\[2ex]
\begin{subfigure}[b]{1\textwidth}
\captionsetup{width=1\textwidth, font=small,labelfont=small}
\centering\includegraphics[scale=0.2]{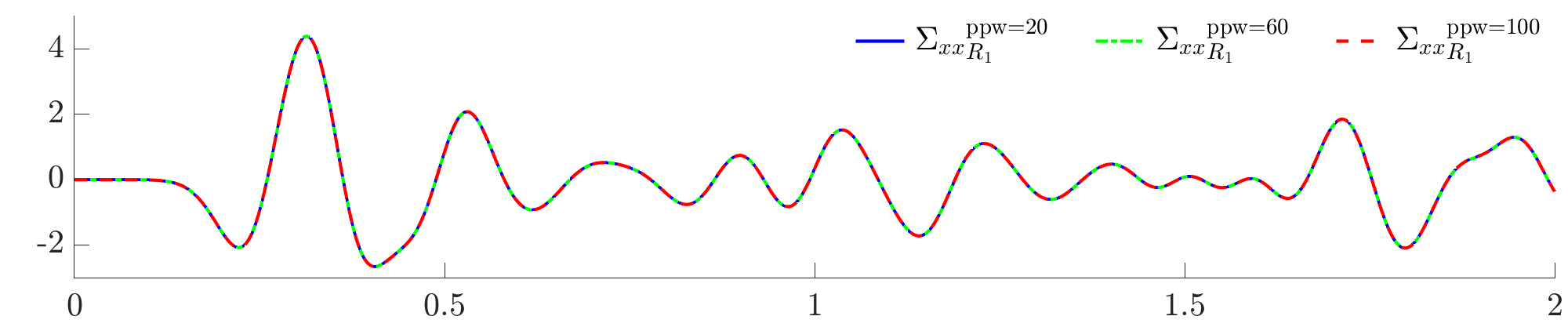}
\caption{Time history of $\Sigma_{xx}$ at 1 grid points below the surface for the case ppw = 20.}
\label{Elastic_full_FSBC_Sxx_1_strong}
\end{subfigure}\hfill
\\[2ex]
\begin{subfigure}[b]{1\textwidth}
\captionsetup{width=1\textwidth, font=small,labelfont=small}
\centering\includegraphics[scale=0.2]{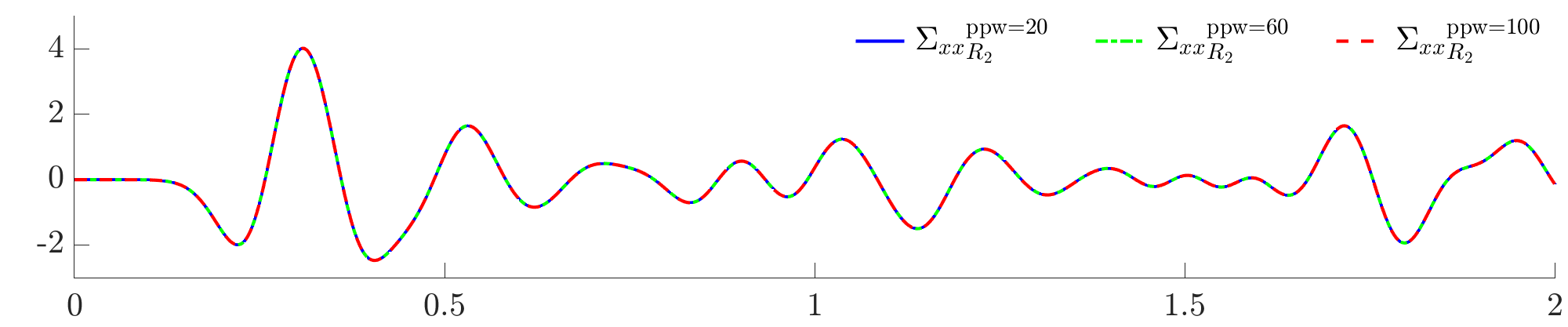}
\caption{Time history of $\Sigma_{xx}$ at 2 grid points below the surface for the case ppw = 20.}
\label{Elastic_full_FSBC_Sxx_2_strong}
\end{subfigure}\hfill
\caption{Time histories of the 2D elastic experiments. The free surface boundary condition is imposed strongly.
}
\label{Elastic_full_FSBC_Sxx_strong}
\end{figure}
\end{minipage}
\hspace{0.0075\textwidth}
\begin{minipage}[t]{.5\textwidth}
\begin{figure}[H]
\captionsetup{width=1\textwidth, font=small,labelfont=small}
\centering
\begin{subfigure}[b]{1\textwidth}
\captionsetup{width=1\textwidth, font=small,labelfont=small}
\centering\includegraphics[scale=0.2]{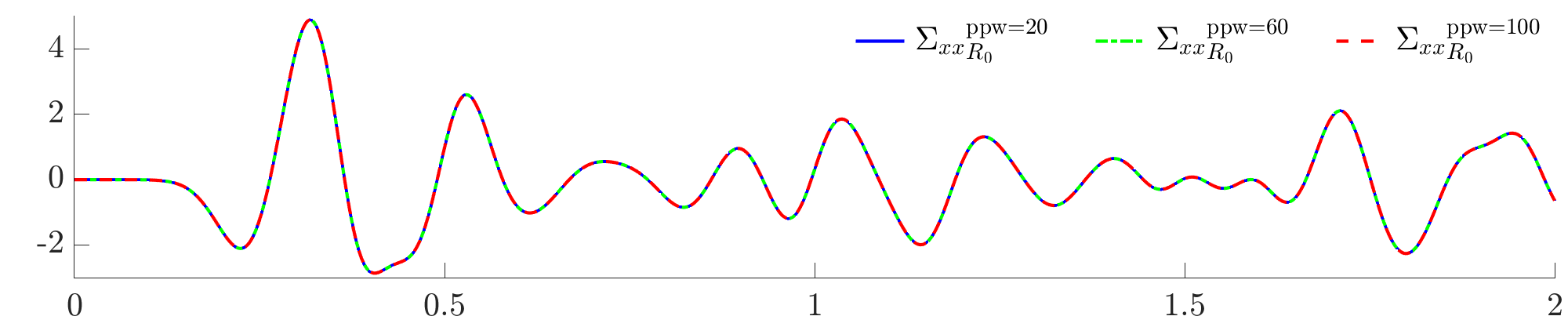}
\caption{Time history of $\Sigma_{xx}$ at the surface.}
\label{Elastic_full_FSBC_Sxx_0_weak}
\end{subfigure}\hfill
\\[2ex]
\begin{subfigure}[b]{1\textwidth}
\captionsetup{width=1\textwidth, font=small,labelfont=small}
\centering\includegraphics[scale=0.2]{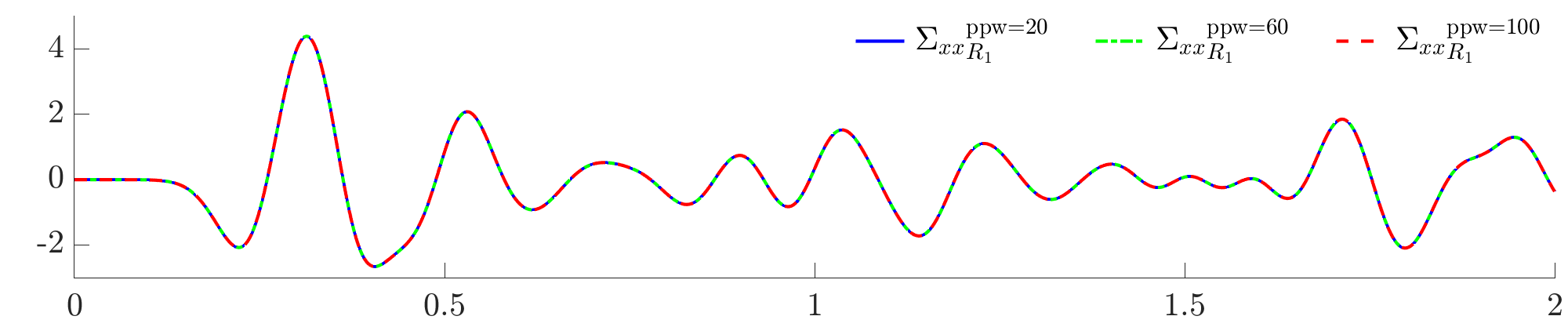}
\caption{Time history of $\Sigma_{xx}$ at 1 grid points below the surface for the case ppw = 20.}
\label{Elastic_full_FSBC_Sxx_1_weak}
\end{subfigure}\hfill
\\[2ex]
\begin{subfigure}[b]{1\textwidth}
\captionsetup{width=1\textwidth, font=small,labelfont=small}
\centering\includegraphics[scale=0.2]{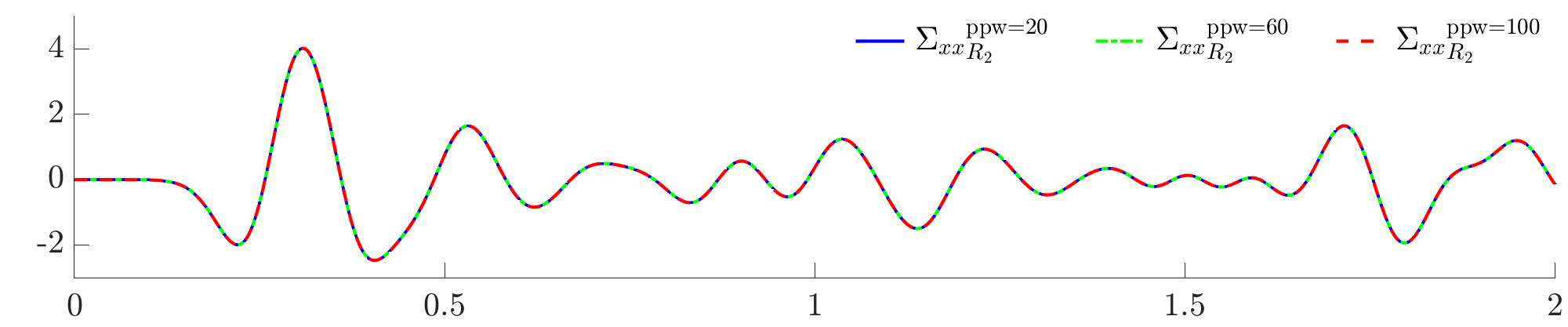}
\caption{Time history of $\Sigma_{xx}$ at 2 grid points below the surface for the case ppw = 20.}
\label{Elastic_full_FSBC_Sxx_2_weak}
\end{subfigure}\hfill
\caption{Time histories of the 2D elastic experiments. The free surface boundary condition is imposed weakly.
}
\label{Elastic_full_FSBC_Sxx_weak}
\end{figure}
\end{minipage}
\ \\
\ \\

\hspace{-0.0625\textwidth}
\begin{minipage}[t]{.5\textwidth}
\begin{figure}[H]
\captionsetup{width=1\textwidth, font=small,labelfont=small}
\centering
\begin{subfigure}[b]{1\textwidth}
\captionsetup{width=1\textwidth, font=small,labelfont=small}
\centering\includegraphics[scale=0.2]{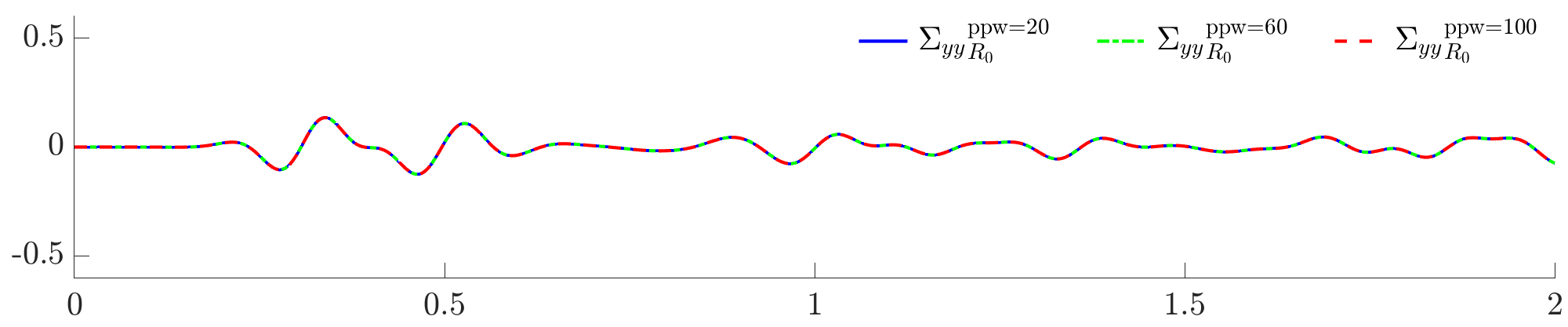}
\caption{Time history of $\Sigma_{yy}$ at the surface.}
\label{Elastic_full_FSBC_Syy_0_strong}
\end{subfigure}\hfill
\\[2ex]
\begin{subfigure}[b]{1\textwidth}
\captionsetup{width=1\textwidth, font=small,labelfont=small}
\centering\includegraphics[scale=0.2]{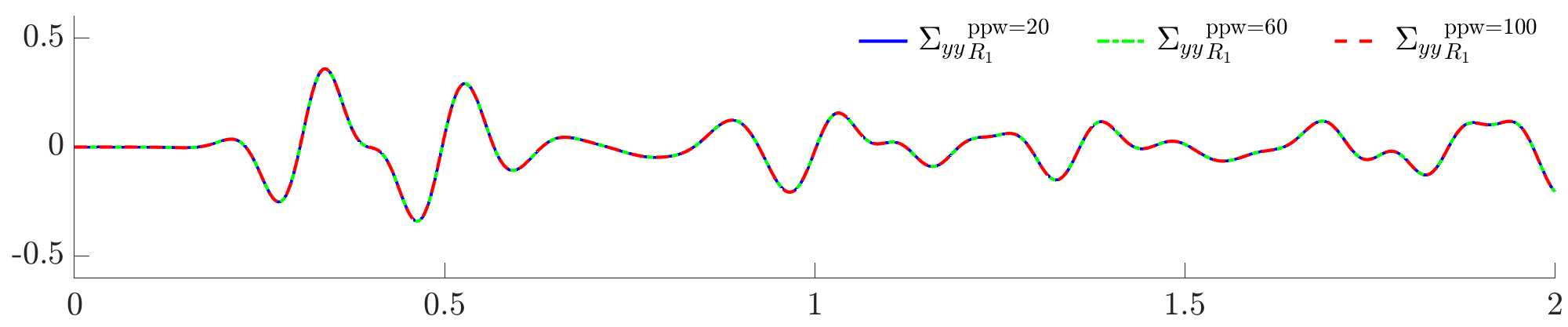}
\caption{Time history of $\Sigma_{yy}$ at 1 grid points below the surface for the case ppw = 20.}
\label{Elastic_full_FSBC_Syy_1_strong}
\end{subfigure}\hfill
\\[2ex]
\begin{subfigure}[b]{1\textwidth}
\captionsetup{width=1\textwidth, font=small,labelfont=small}
\centering\includegraphics[scale=0.2]{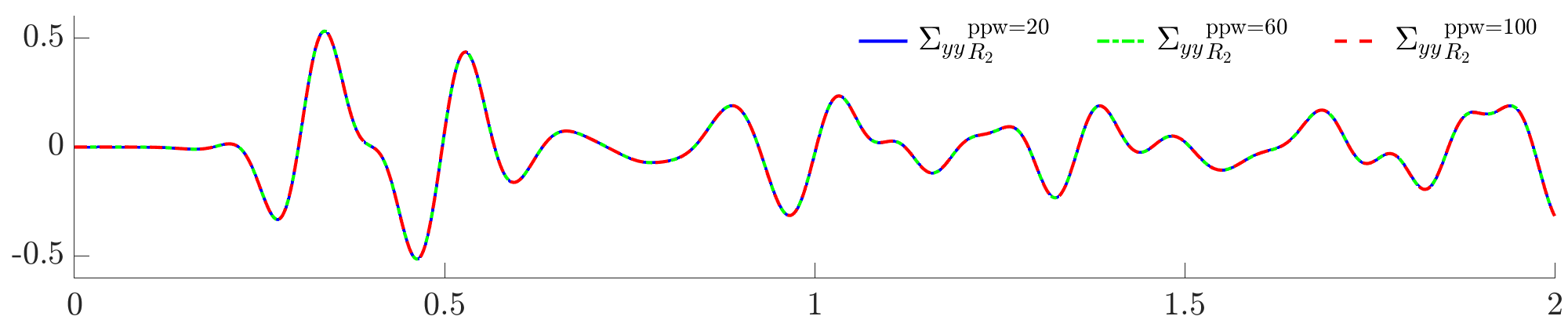}
\caption{Time history of $\Sigma_{yy}$ at 2 grid points below the surface for the case ppw = 20.}
\label{Elastic_full_FSBC_Syy_2_strong}
\end{subfigure}\hfill
\caption{Time histories of the 2D elastic experiments. The free surface boundary condition is imposed strongly.
}
\label{Elastic_full_FSBC_Syy_strong}
\end{figure}
\end{minipage}
\hspace{0.0075\textwidth}
\begin{minipage}[t]{.5\textwidth}
\begin{figure}[H]
\captionsetup{width=1\textwidth, font=small,labelfont=small}
\centering
\begin{subfigure}[b]{1\textwidth}
\captionsetup{width=1\textwidth, font=small,labelfont=small}
\centering\includegraphics[scale=0.2]{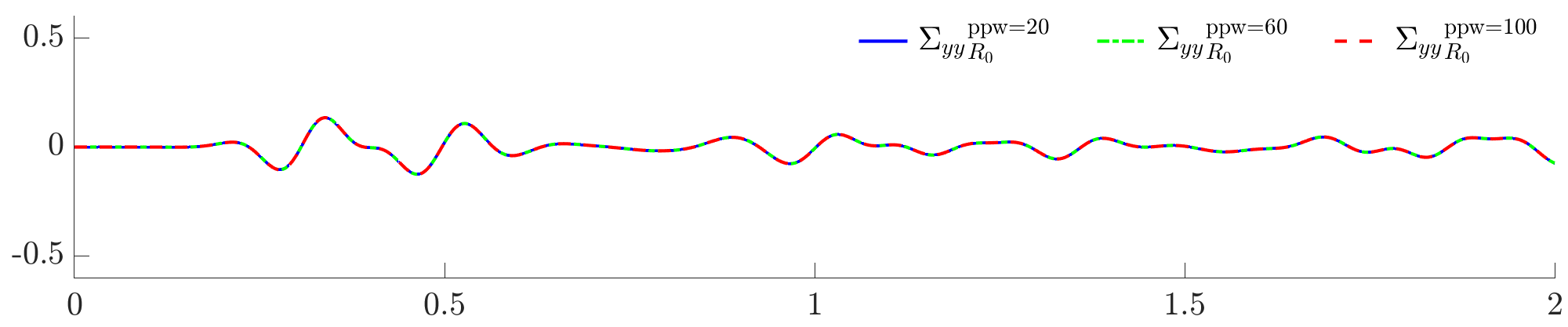}
\caption{Time history of $\Sigma_{yy}$ at the surface.}
\label{Elastic_full_FSBC_Syy_0_weak}
\end{subfigure}\hfill
\\[2ex]
\begin{subfigure}[b]{1\textwidth}
\captionsetup{width=1\textwidth, font=small,labelfont=small}
\centering\includegraphics[scale=0.2]{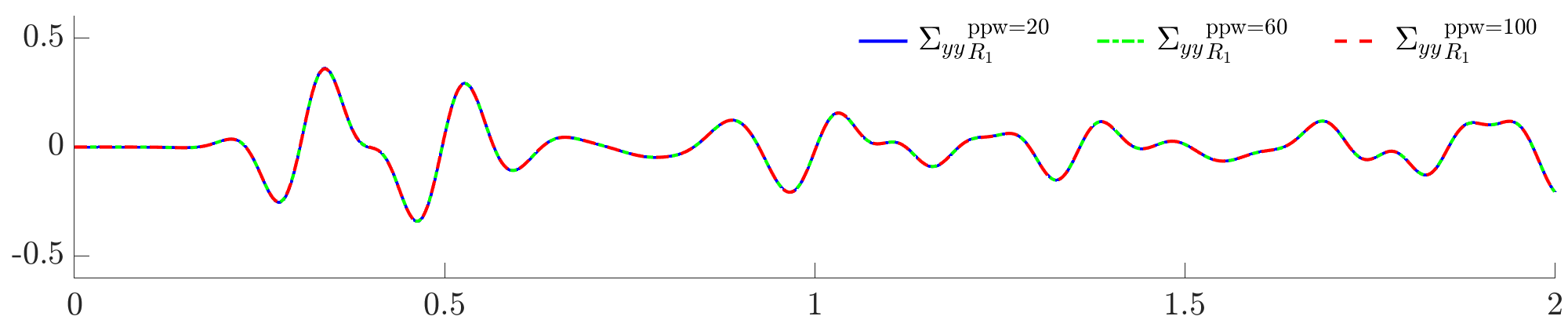}
\caption{Time history of $\Sigma_{yy}$ at 1 grid points below the surface for the case ppw = 20.}
\label{Elastic_full_FSBC_Syy_1_weak}
\end{subfigure}\hfill
\\[2ex]
\begin{subfigure}[b]{1\textwidth}
\captionsetup{width=1\textwidth, font=small,labelfont=small}
\centering\includegraphics[scale=0.2]{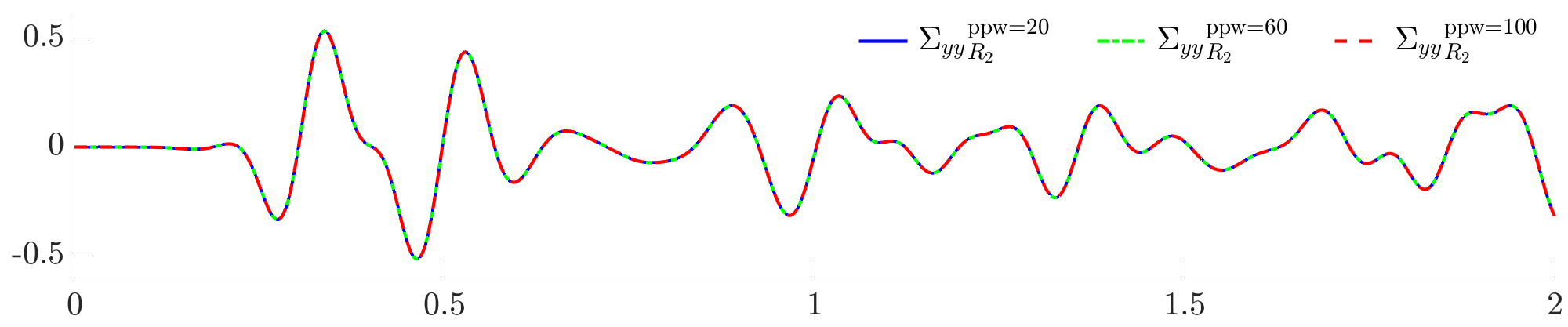}
\caption{Time history of $\Sigma_{yy}$ at 2 grid points below the surface for the case ppw = 20.}
\label{Elastic_full_FSBC_Syy_2_weak}
\end{subfigure}\hfill
\caption{Time histories of the 2D elastic experiments. The free surface boundary condition is imposed weakly.
}
\label{Elastic_full_FSBC_Syy_weak}
\end{figure}
\end{minipage}
\ \\
\ \\

\hspace{-0.0625\textwidth}
\begin{minipage}[t]{.5\textwidth}
\begin{figure}[H]
\captionsetup{width=1\textwidth, font=small,labelfont=small}
\centering
\begin{subfigure}[b]{1\textwidth}
\captionsetup{width=1\textwidth, font=small,labelfont=small}
\centering\includegraphics[scale=0.2]{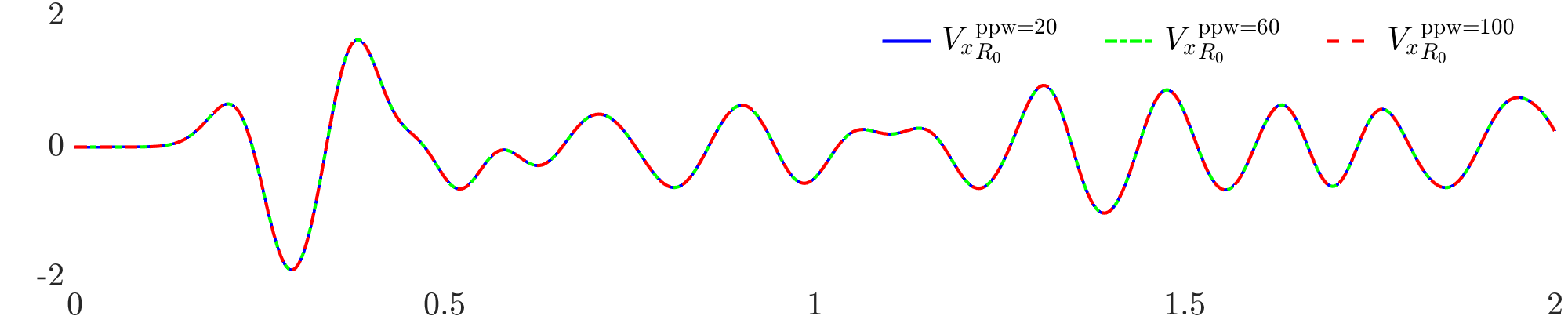}
\caption{Time history of $V_x$ at the surface.}
\label{Elastic_full_FSBC_Vx_0_strong}
\end{subfigure}\hfill
\\[2ex]
\begin{subfigure}[b]{1\textwidth}
\captionsetup{width=1\textwidth, font=small,labelfont=small}
\centering\includegraphics[scale=0.2]{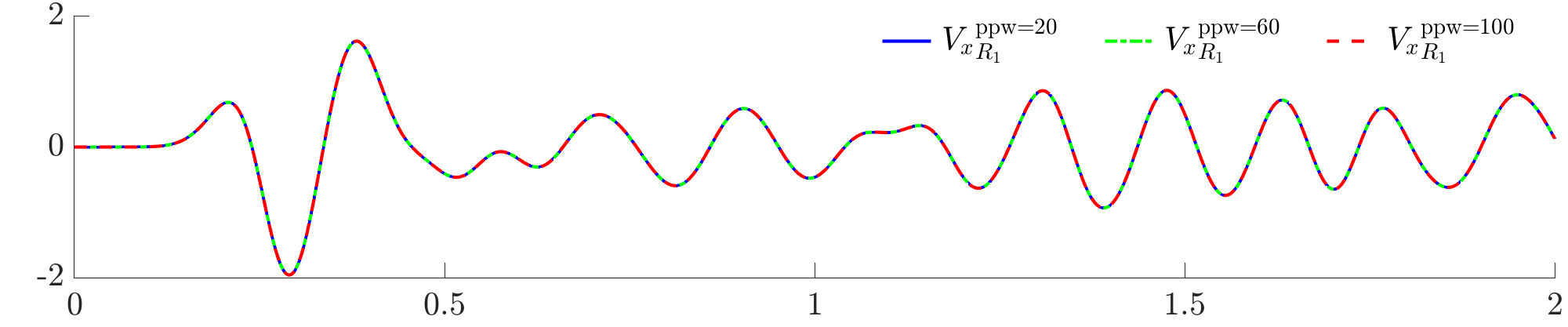}
\caption{Time history of $V_x$ at 1 grid points below the surface for the case ppw = 20.}
\label{Elastic_full_FSBC_Vx_1_strong}
\end{subfigure}\hfill
\\[2ex]
\begin{subfigure}[b]{1\textwidth}
\captionsetup{width=1\textwidth, font=small,labelfont=small}
\centering\includegraphics[scale=0.2]{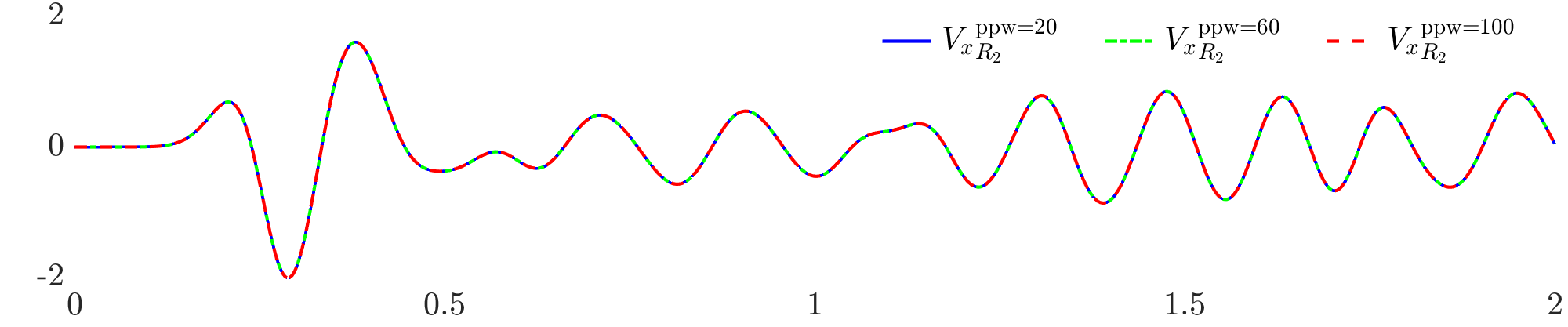}
\caption{Time history of $V_x$ at 2 grid points below the surface for the case ppw = 20.}
\label{Elastic_full_FSBC_Vx_2_strong}
\end{subfigure}\hfill
\caption{Time histories of the 2D elastic experiments. The free surface boundary condition is imposed strongly.
}
\label{Elastic_full_FSBC_Vx_strong}
\end{figure}
\end{minipage}
\hspace{0.0075\textwidth}
\begin{minipage}[t]{.5\textwidth}
\begin{figure}[H]
\captionsetup{width=1\textwidth, font=small,labelfont=small}
\centering
\begin{subfigure}[b]{1\textwidth}
\captionsetup{width=1\textwidth, font=small,labelfont=small}
\centering\includegraphics[scale=0.2]{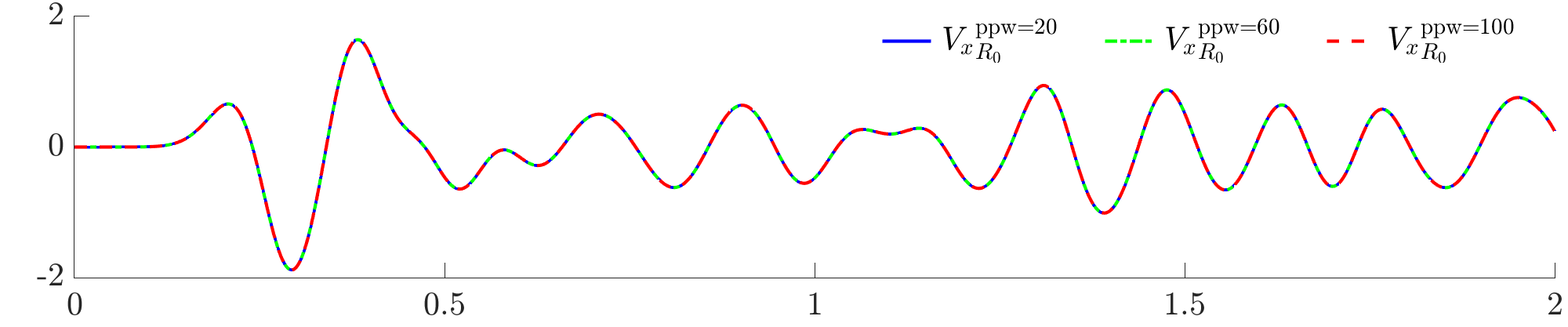}
\caption{Time history of $V_x$ at the surface.}
\label{Elastic_full_FSBC_Vx_0_weak}
\end{subfigure}\hfill
\\[2ex]
\begin{subfigure}[b]{1\textwidth}
\captionsetup{width=1\textwidth, font=small,labelfont=small}
\centering\includegraphics[scale=0.2]{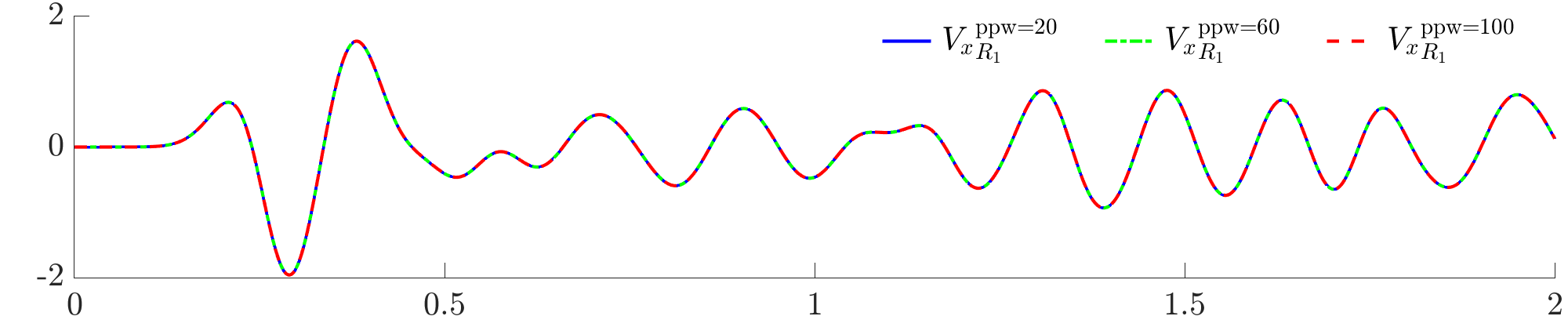}
\caption{Time history of $V_x$ at 1 grid points below the surface for the case ppw = 20.}
\label{Elastic_full_FSBC_Vx_1_weak}
\end{subfigure}\hfill
\\[2ex]
\begin{subfigure}[b]{1\textwidth}
\captionsetup{width=1\textwidth, font=small,labelfont=small}
\centering\includegraphics[scale=0.2]{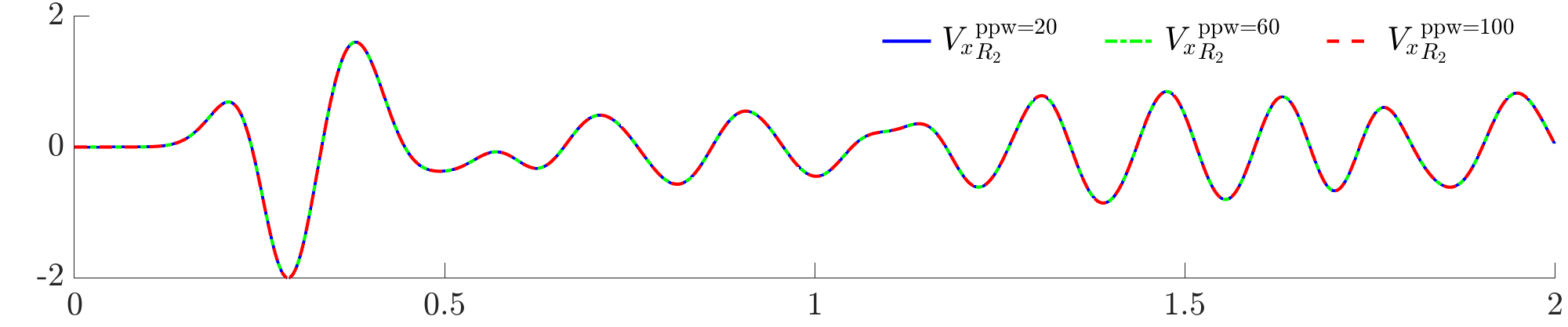}
\caption{Time history of $V_x$ at 2 grid points below the surface for the case ppw = 20.}
\label{Elastic_full_FSBC_Vx_2_weak}
\end{subfigure}\hfill
\caption{Time histories of the 2D elastic experiments. The free surface boundary condition is imposed weakly.
}
\label{Elastic_full_FSBC_Vx_weak}
\end{figure}
\end{minipage}
\ \\
\ \\

\hspace{-0.0625\textwidth}
\begin{minipage}[t]{.5\textwidth}
\begin{figure}[H]
\captionsetup{width=1\textwidth, font=small,labelfont=small}
\centering
\begin{subfigure}[b]{1\textwidth}
\captionsetup{width=1\textwidth, font=small,labelfont=small}
\centering\includegraphics[scale=0.2]{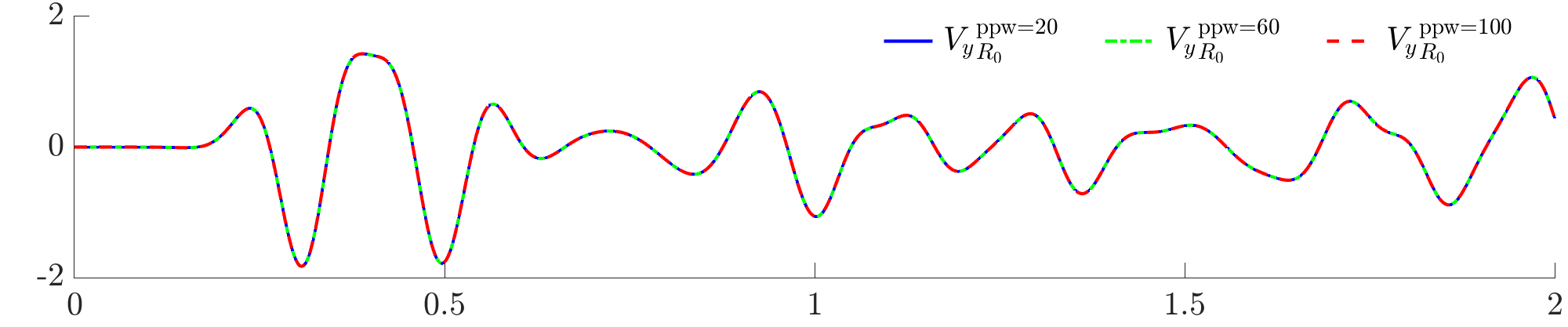}
\caption{Time history of $V_y$ at the surface.}
\label{Elastic_full_FSBC_Vy_0_strong}
\end{subfigure}\hfill
\\[2ex]
\begin{subfigure}[b]{1\textwidth}
\captionsetup{width=1\textwidth, font=small,labelfont=small}
\centering\includegraphics[scale=0.2]{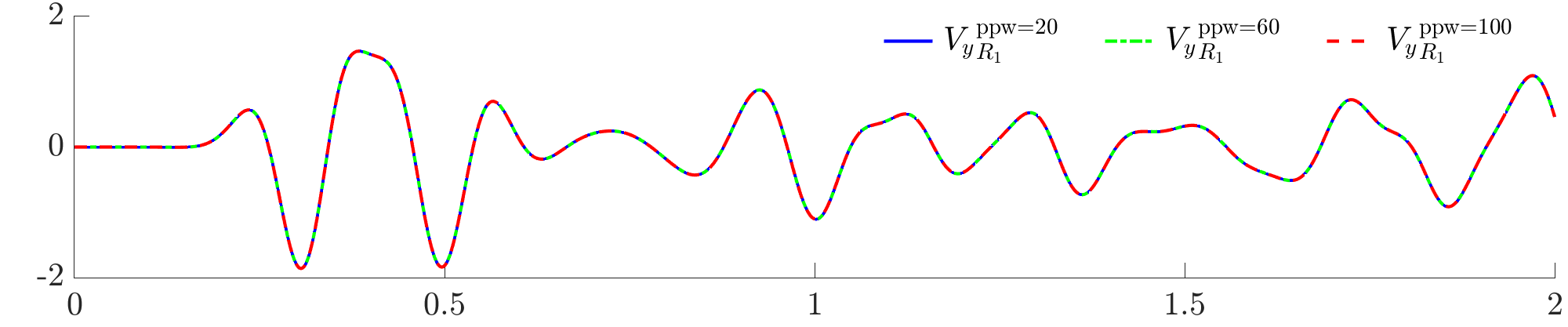}
\caption{Time history of $V_y$ at 1 grid points below the surface for the case ppw = 20.}
\label{Elastic_full_FSBC_Vy_1_strong}
\end{subfigure}\hfill
\\[2ex]
\begin{subfigure}[b]{1\textwidth}
\captionsetup{width=1\textwidth, font=small,labelfont=small}
\centering\includegraphics[scale=0.2]{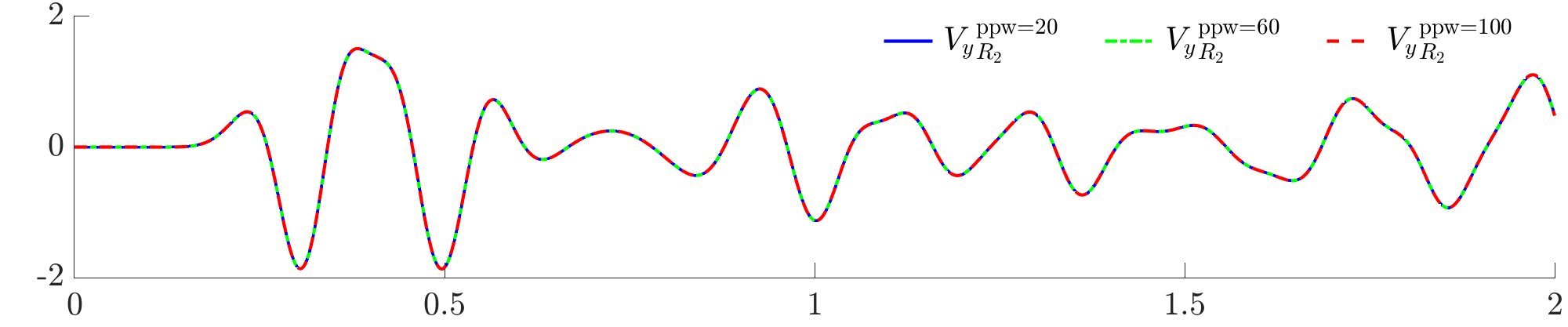}
\caption{Time history of $V_y$ at 2 grid points below the surface for the case ppw = 20.}
\label{Elastic_full_FSBC_Vy_2_strong}
\end{subfigure}\hfill
\caption{Time histories of the 2D elastic experiments. The free surface boundary condition is imposed strongly.
}
\label{Elastic_full_FSBC_Vy_strong}
\end{figure}
\end{minipage}
\hspace{0.0075\textwidth}
\begin{minipage}[t]{.5\textwidth}
\begin{figure}[H]
\captionsetup{width=1\textwidth, font=small,labelfont=small}
\centering
\begin{subfigure}[b]{1\textwidth}
\captionsetup{width=1\textwidth, font=small,labelfont=small}
\centering\includegraphics[scale=0.2]{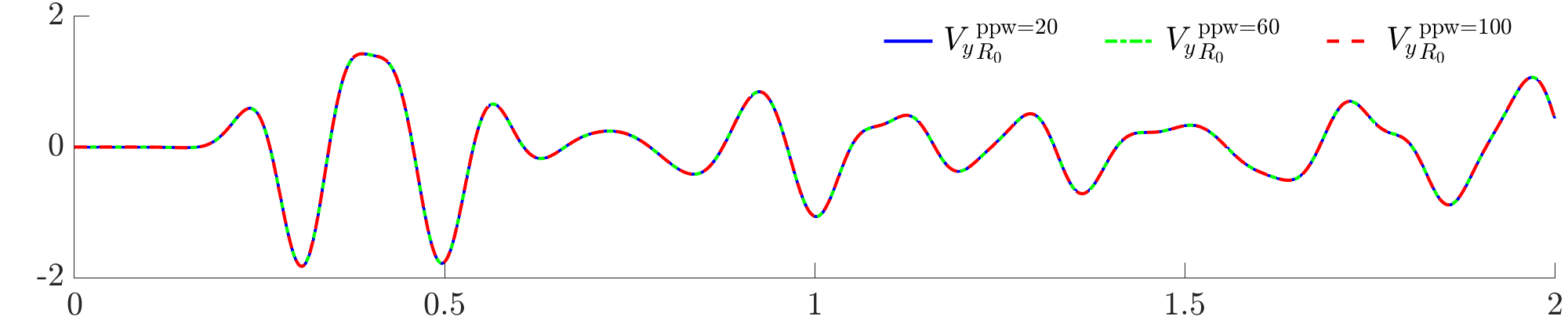}
\caption{Time history of $V_y$ at the surface.}
\label{Elastic_full_FSBC_Vy_0_weak}
\end{subfigure}\hfill
\\[2ex]
\begin{subfigure}[b]{1\textwidth}
\captionsetup{width=1\textwidth, font=small,labelfont=small}
\centering\includegraphics[scale=0.2]{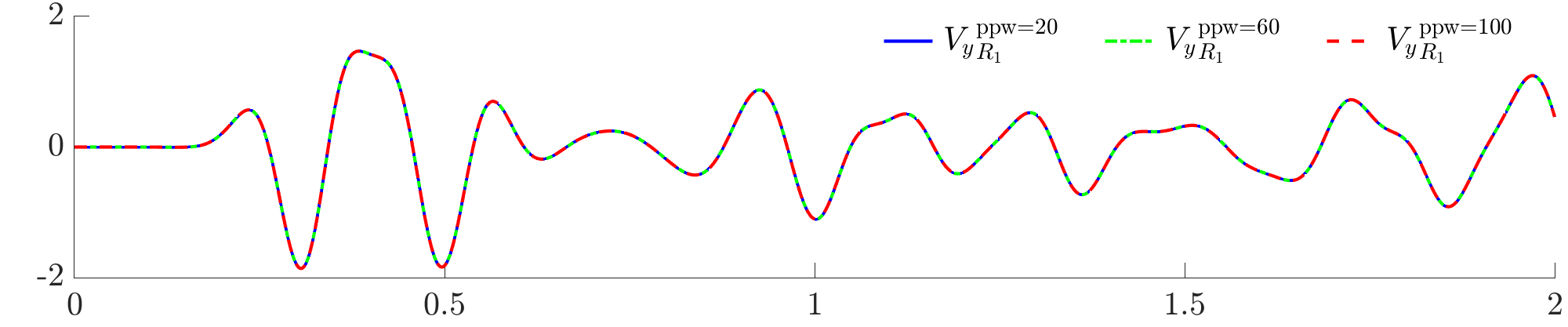}
\caption{Time history of $V_y$ at 1 grid points below the surface for the case ppw = 20.}
\label{Elastic_full_FSBC_Vy_1_weak}
\end{subfigure}\hfill
\\[2ex]
\begin{subfigure}[b]{1\textwidth}
\captionsetup{width=1\textwidth, font=small,labelfont=small}
\centering\includegraphics[scale=0.2]{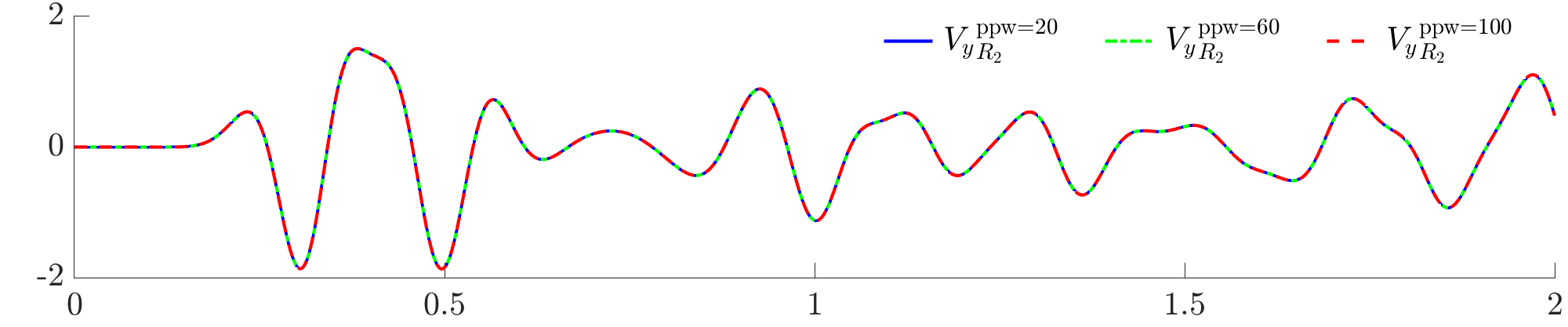}
\caption{Time history of $V_y$ at 2 grid points below the surface for the case ppw = 20.}
\label{Elastic_full_FSBC_Vy_2_weak}
\end{subfigure}\hfill
\caption{Time histories of the 2D elastic experiments. The free surface boundary condition is imposed weakly.
}
\label{Elastic_full_FSBC_Vy_weak}
\end{figure}
\end{minipage}

%% file: input_source_on_Vy_corner.tex
\hspace{-0.0625\textwidth}
\begin{minipage}[t]{.5\textwidth}
\begin{figure}[H]
\captionsetup{width=1\textwidth, font=footnotesize,labelfont=footnotesize}
\centering
\begin{subfigure}[b]{1\textwidth}
\captionsetup{width=1\textwidth, font=footnotesize,labelfont=footnotesize}
\centering\includegraphics[scale=0.2]{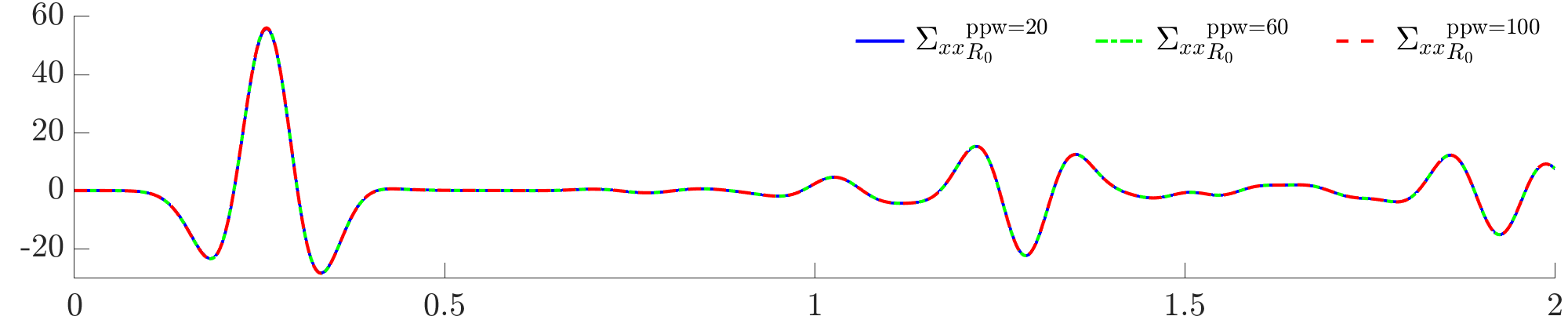}
\caption{Time history of $\Sigma_{xx}$ at the surface.}
\label{stress_test_Vy_elastic_full_FSBC_Sxx_0_strong}
\end{subfigure}\hfill
\\[2ex]
\begin{subfigure}[b]{1\textwidth}
\captionsetup{width=1\textwidth, font=footnotesize,labelfont=footnotesize}
\centering\includegraphics[scale=0.2]{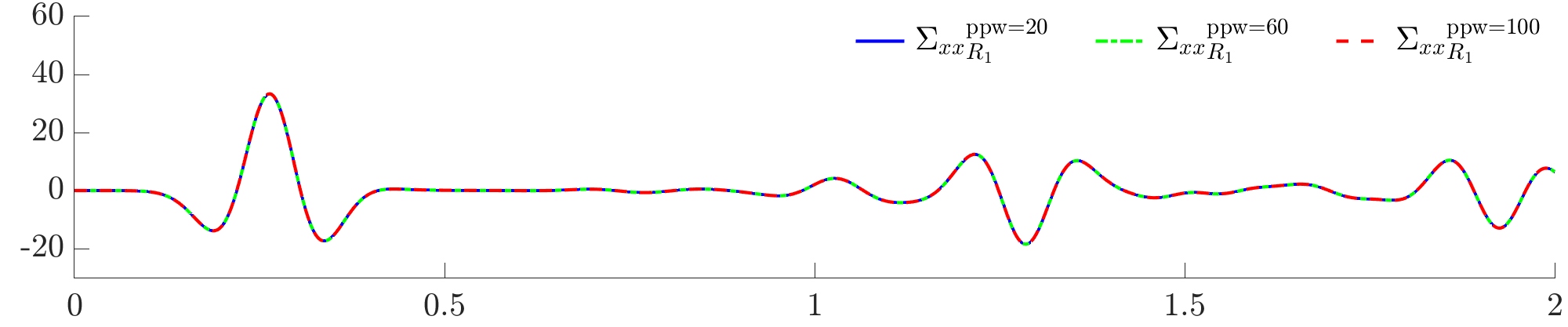}
\caption{Time history of $\Sigma_{xx}$ at 1 grid points below the surface for the case ppw = 20.}
\label{stress_test_Vy_elastic_full_FSBC_Sxx_1_strong}
\end{subfigure}\hfill
\\[2ex]
\begin{subfigure}[b]{1\textwidth}
\captionsetup{width=1\textwidth, font=footnotesize,labelfont=footnotesize}
\centering\includegraphics[scale=0.2]{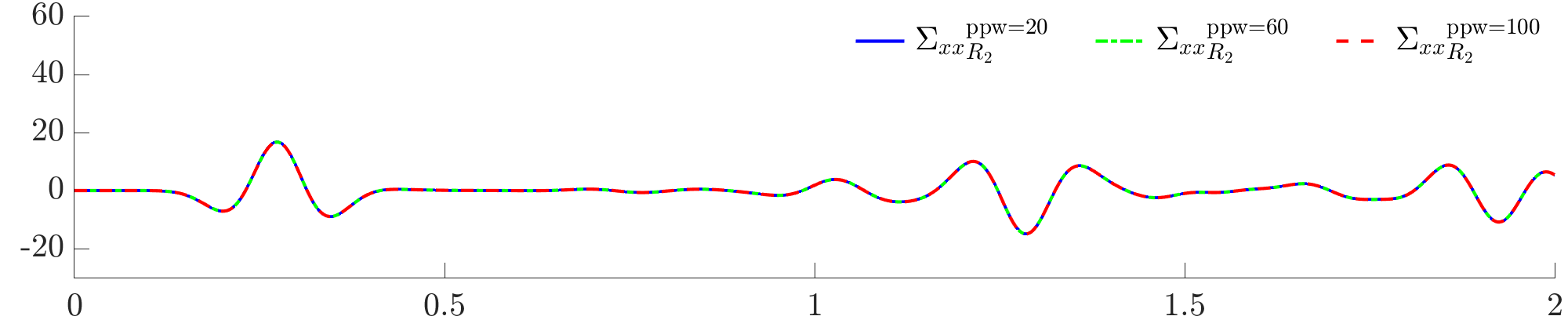}
\caption{Time history of $\Sigma_{xx}$ at 2 grid points below the surface for the case ppw = 20.}
\label{stress_test_Vy_elastic_full_FSBC_Sxx_2_strong}
\end{subfigure}\hfill
\caption{Time histories of the 2D elastic experiments. The free surface boundary condition is imposed strongly.
}
\label{stress_test_Vy_elastic_full_FSBC_Sxx_strong}
\end{figure}
\end{minipage}
\hspace{0.0075\textwidth}
\begin{minipage}[t]{.5\textwidth}
\begin{figure}[H]
\captionsetup{width=1\textwidth, font=footnotesize,labelfont=footnotesize}
\centering
\begin{subfigure}[b]{1\textwidth}
\captionsetup{width=1\textwidth, font=footnotesize,labelfont=footnotesize}
\centering\includegraphics[scale=0.2]{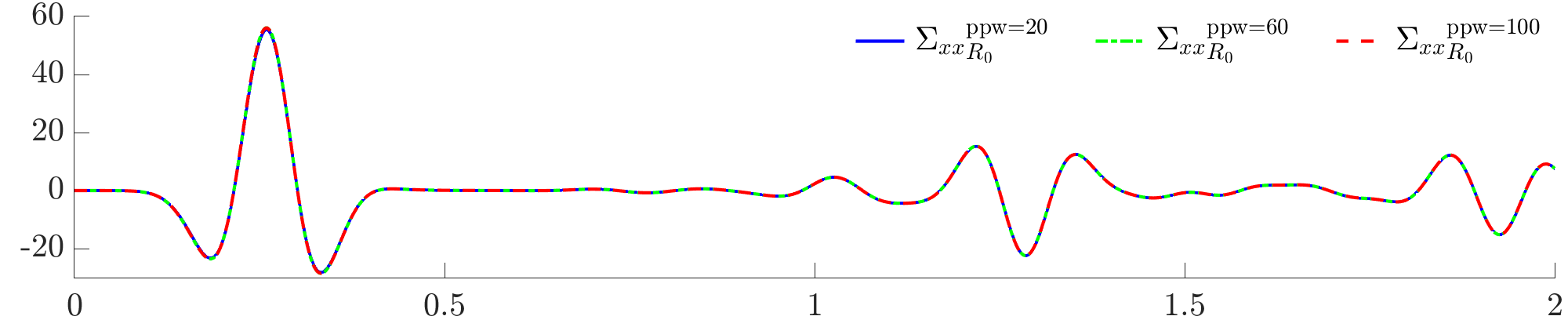}
\caption{Time history of $\Sigma_{xx}$ at the surface.}
\label{stress_test_Vy_elastic_full_FSBC_Sxx_0_weak}
\end{subfigure}\hfill
\\[2ex]
\begin{subfigure}[b]{1\textwidth}
\captionsetup{width=1\textwidth, font=footnotesize,labelfont=footnotesize}
\centering\includegraphics[scale=0.2]{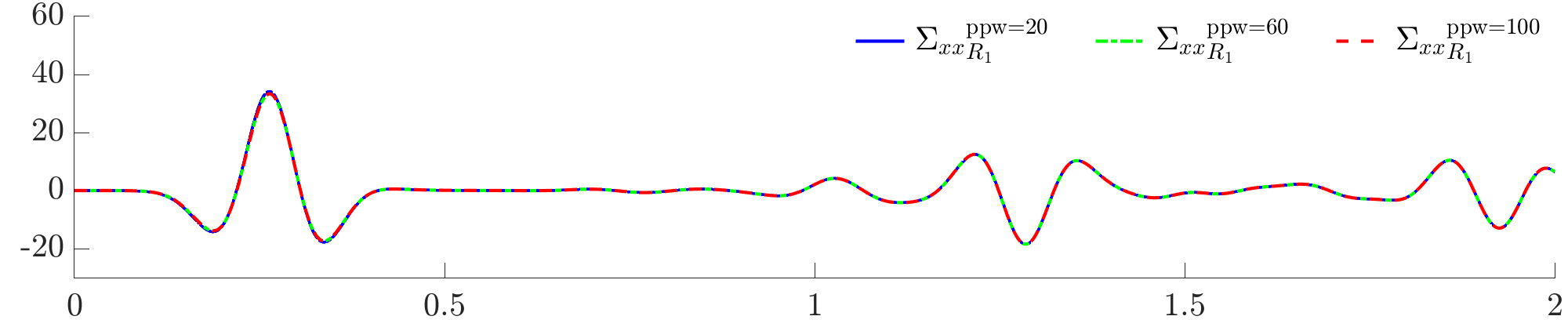}
\caption{Time history of $\Sigma_{xx}$ at 1 grid points below the surface for the case ppw = 20.}
\label{stress_test_Vy_elastic_full_FSBC_Sxx_1_weak}
\end{subfigure}\hfill
\\[2ex]
\begin{subfigure}[b]{1\textwidth}
\captionsetup{width=1\textwidth, font=footnotesize,labelfont=footnotesize}
\centering\includegraphics[scale=0.2]{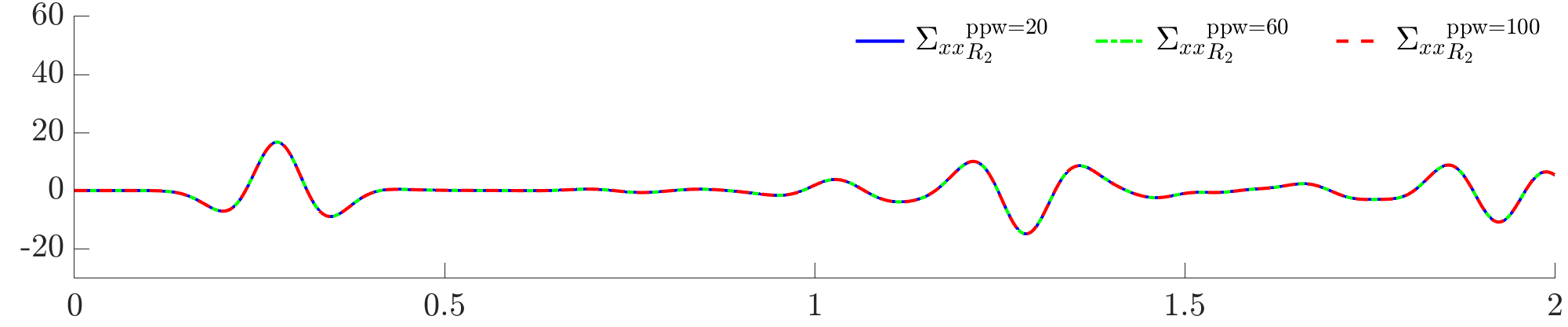}
\caption{Time history of $\Sigma_{xx}$ at 2 grid points below the surface for the case ppw = 20.}
\label{stress_test_Vy_elastic_full_FSBC_Sxx_2_weak}
\end{subfigure}\hfill
\caption{Time histories of the 2D elastic experiments. The free surface boundary condition is imposed weakly.
}
\label{stress_test_Vy_elastic_full_FSBC_Sxx_weak}
\end{figure}
\end{minipage}
\ \\
\ \\

\hspace{-0.0625\textwidth}
\begin{minipage}[t]{.5\textwidth}
\begin{figure}[H]
\captionsetup{width=1\textwidth, font=footnotesize,labelfont=footnotesize}
\centering
\begin{subfigure}[b]{1\textwidth}
\captionsetup{width=1\textwidth, font=footnotesize,labelfont=footnotesize}
\centering\includegraphics[scale=0.2]{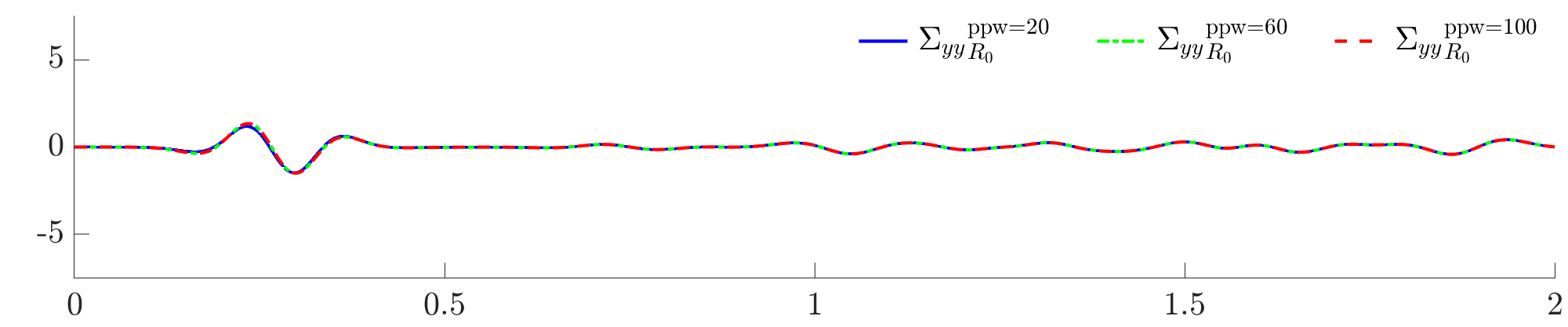}
\caption{Time history of $\Sigma_{yy}$ at the surface.}
\label{stress_test_Vy_elastic_full_FSBC_Syy_0_strong}
\end{subfigure}\hfill
\\[2ex]
\begin{subfigure}[b]{1\textwidth}
\captionsetup{width=1\textwidth, font=footnotesize,labelfont=footnotesize}
\centering\includegraphics[scale=0.2]{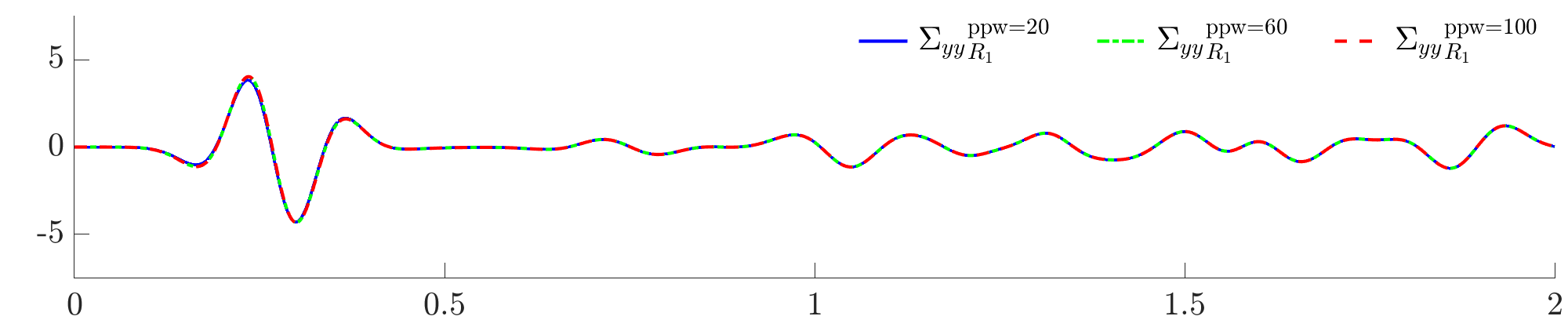}
\caption{Time history of $\Sigma_{yy}$ at 1 grid points below the surface for the case ppw = 20.}
\label{stress_test_Vy_elastic_full_FSBC_Syy_1_strong}
\end{subfigure}\hfill
\\[2ex]
\begin{subfigure}[b]{1\textwidth}
\captionsetup{width=1\textwidth, font=footnotesize,labelfont=footnotesize}
\centering\includegraphics[scale=0.2]{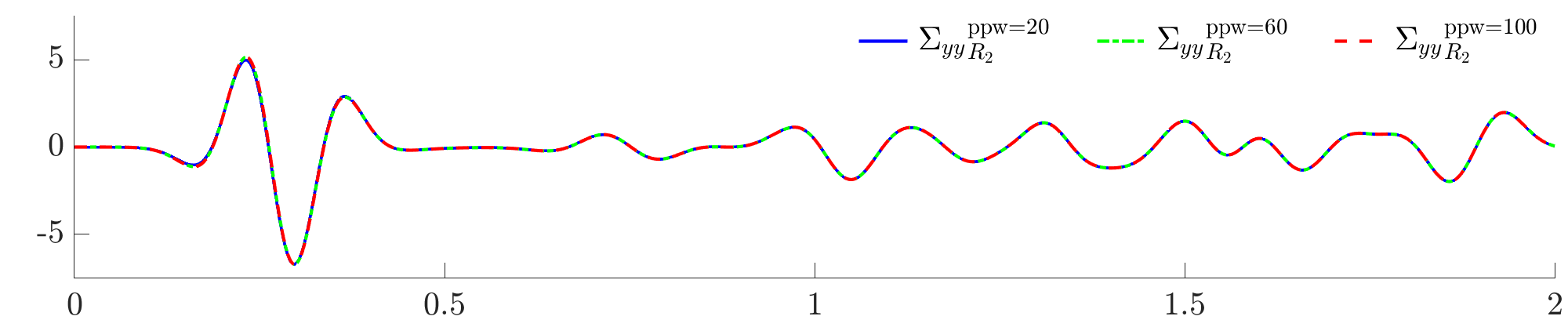}
\caption{Time history of $\Sigma_{yy}$ at 2 grid points below the surface for the case ppw = 20.}
\label{stress_test_Vy_elastic_full_FSBC_Syy_2_strong}
\end{subfigure}\hfill
\caption{Time histories of the 2D elastic experiments. The free surface boundary condition is imposed strongly.
}
\label{stress_test_Vy_elastic_full_FSBC_Syy_strong}
\end{figure}
\end{minipage}
\hspace{0.0075\textwidth}
\begin{minipage}[t]{.5\textwidth}
\begin{figure}[H]
\captionsetup{width=1\textwidth, font=footnotesize,labelfont=footnotesize}
\centering
\begin{subfigure}[b]{1\textwidth}
\captionsetup{width=1\textwidth, font=footnotesize,labelfont=footnotesize}
\centering\includegraphics[scale=0.2]{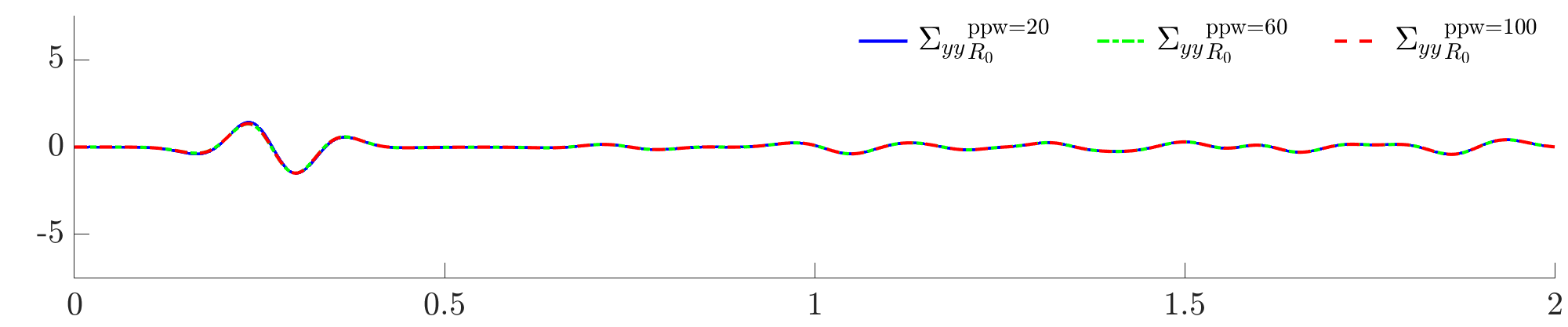}
\caption{Time history of $\Sigma_{yy}$ at the surface.}
\label{stress_test_Vy_elastic_full_FSBC_Syy_0_weak}
\end{subfigure}\hfill
\\[2ex]
\begin{subfigure}[b]{1\textwidth}
\captionsetup{width=1\textwidth, font=footnotesize,labelfont=footnotesize}
\centering\includegraphics[scale=0.2]{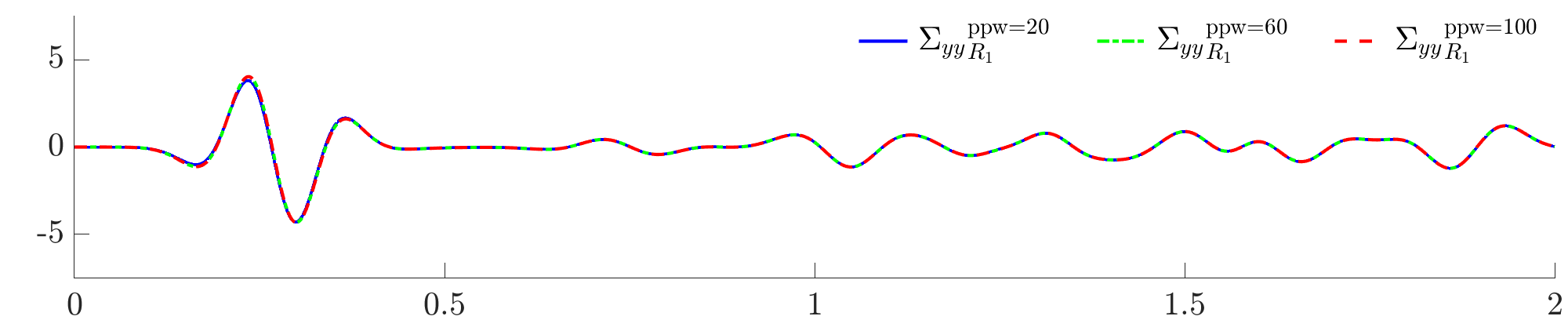}
\caption{Time history of $\Sigma_{yy}$ at 1 grid points below the surface for the case ppw = 20.}
\label{stress_test_Vy_elastic_full_FSBC_Syy_1_weak}
\end{subfigure}\hfill
\\[2ex]
\begin{subfigure}[b]{1\textwidth}
\captionsetup{width=1\textwidth, font=footnotesize,labelfont=footnotesize}
\centering\includegraphics[scale=0.2]{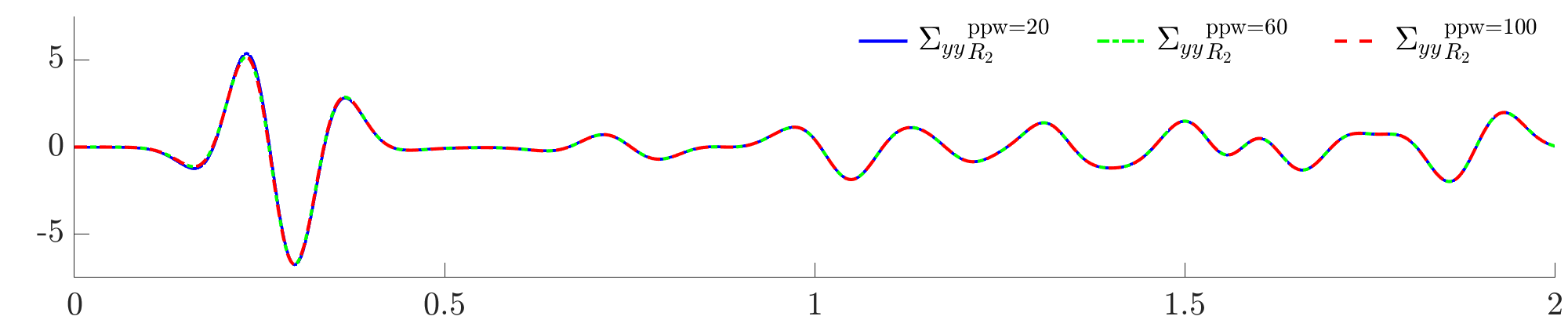}
\caption{Time history of $\Sigma_{yy}$ at 2 grid points below the surface for the case ppw = 20.}
\label{stress_test_Vy_elastic_full_FSBC_Syy_2_weak}
\end{subfigure}\hfill
\caption{Time histories of the 2D elastic experiments. The free surface boundary condition is imposed weakly.
}
\label{stress_test_Vy_elastic_full_FSBC_Syy_weak}
\end{figure}
\end{minipage}
\ \\
\ \\

\hspace{-0.0625\textwidth}
\begin{minipage}[t]{.5\textwidth}
\begin{figure}[H]
\captionsetup{width=1\textwidth, font=footnotesize,labelfont=footnotesize}
\centering
\begin{subfigure}[b]{1\textwidth}
\captionsetup{width=1\textwidth, font=footnotesize,labelfont=footnotesize}
\centering\includegraphics[scale=0.2]{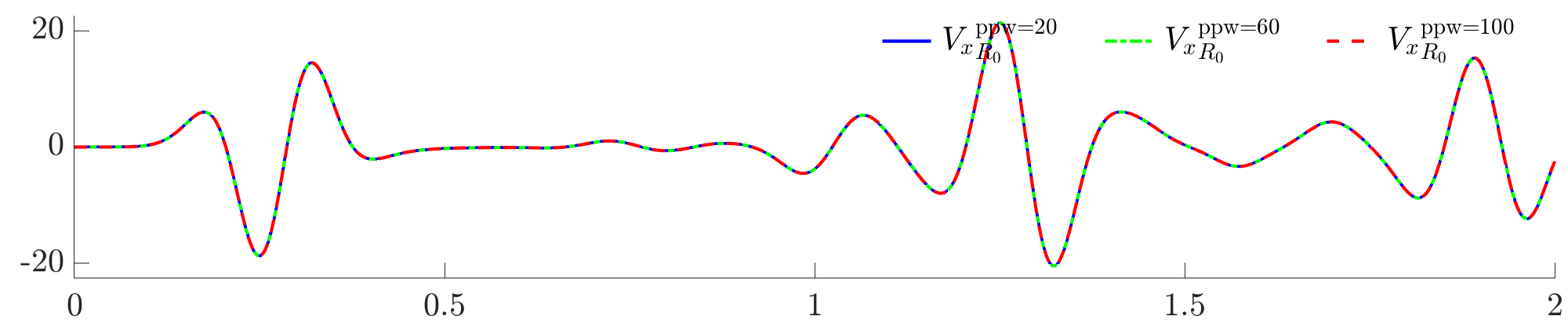}
\caption{Time history of $V_x$ at the surface.}
\label{stress_test_Vy_elastic_full_FSBC_Vx_0_strong}
\end{subfigure}\hfill
\\[2ex]
\begin{subfigure}[b]{1\textwidth}
\captionsetup{width=1\textwidth, font=footnotesize,labelfont=footnotesize}
\centering\includegraphics[scale=0.2]{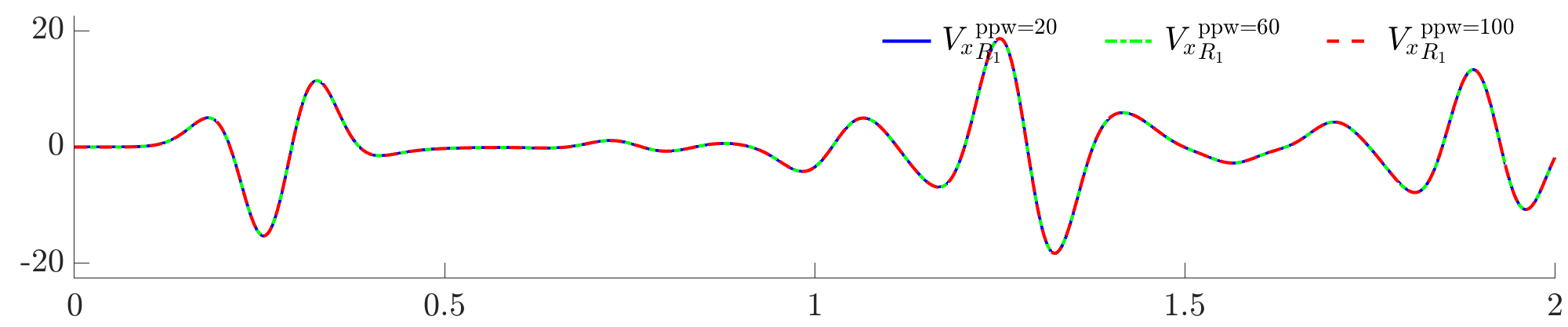}
\caption{Time history of $V_x$ at 1 grid points below the surface for the case ppw = 20.}
\label{stress_test_Vy_elastic_full_FSBC_Vx_1_strong}
\end{subfigure}\hfill
\\[2ex]
\begin{subfigure}[b]{1\textwidth}
\captionsetup{width=1\textwidth, font=footnotesize,labelfont=footnotesize}
\centering\includegraphics[scale=0.2]{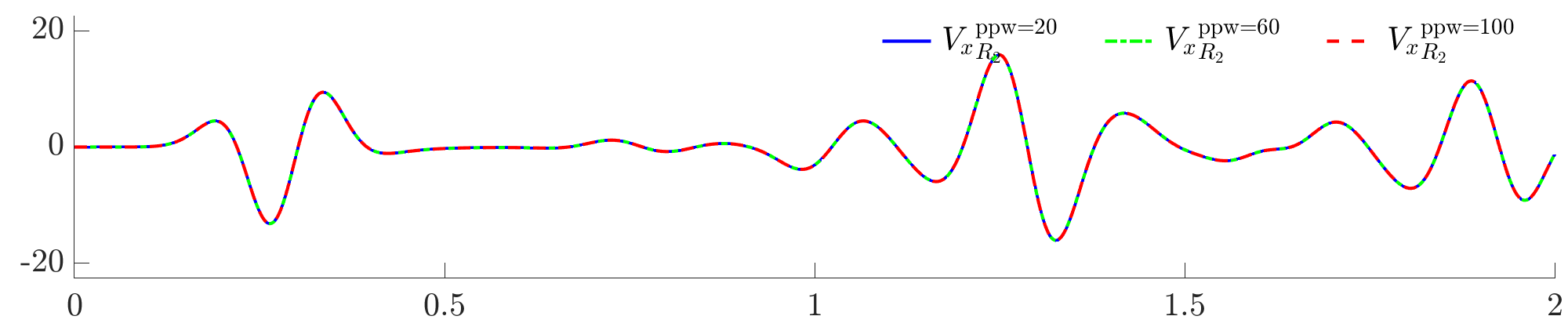}
\caption{Time history of $V_x$ at 2 grid points below the surface for the case ppw = 20.}
\label{stress_test_Vy_elastic_full_FSBC_Vx_2_strong}
\end{subfigure}\hfill
\caption{Time histories of the 2D elastic experiments. The free surface boundary condition is imposed strongly.
}
\label{stress_test_Vy_elastic_full_FSBC_Vx_strong}
\end{figure}
\end{minipage}
\hspace{0.0075\textwidth}
\begin{minipage}[t]{.5\textwidth}
\begin{figure}[H]
\captionsetup{width=1\textwidth, font=footnotesize,labelfont=footnotesize}
\centering
\begin{subfigure}[b]{1\textwidth}
\captionsetup{width=1\textwidth, font=footnotesize,labelfont=footnotesize}
\centering\includegraphics[scale=0.2]{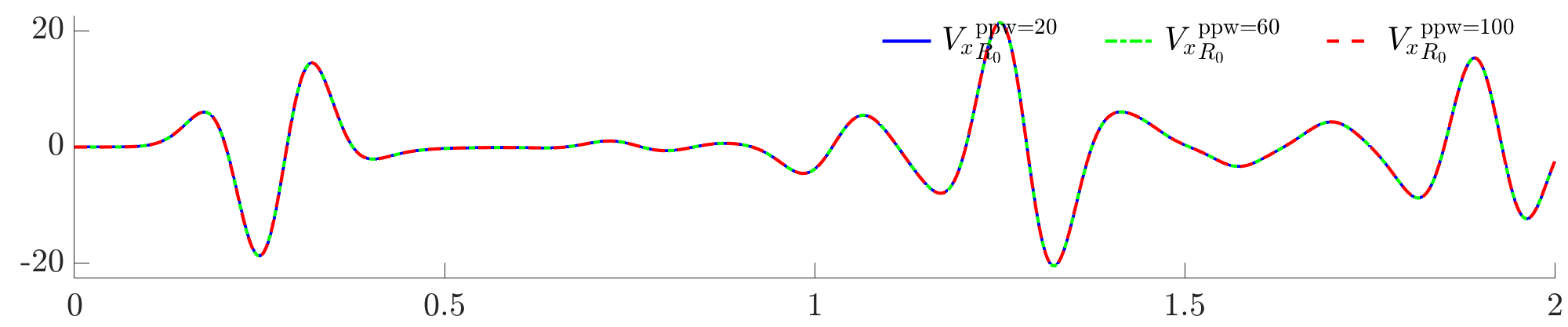}
\caption{Time history of $V_x$ at the surface.}
\label{stress_test_Vy_elastic_full_FSBC_Vx_0_weak}
\end{subfigure}\hfill
\\[2ex]
\begin{subfigure}[b]{1\textwidth}
\captionsetup{width=1\textwidth, font=footnotesize,labelfont=footnotesize}
\centering\includegraphics[scale=0.2]{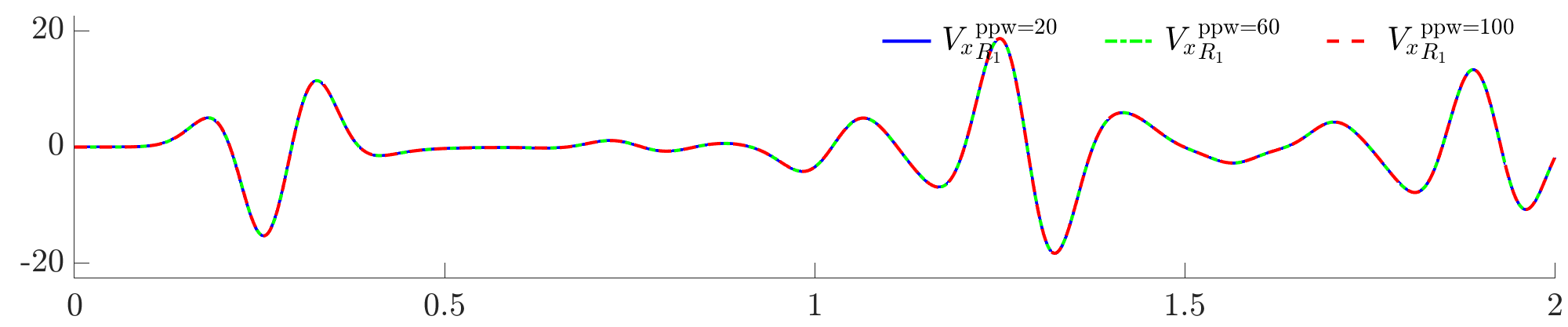}
\caption{Time history of $V_x$ at 1 grid points below the surface for the case ppw = 20.}
\label{stress_test_Vy_elastic_full_FSBC_Vx_1_weak}
\end{subfigure}\hfill
\\[2ex]
\begin{subfigure}[b]{1\textwidth}
\captionsetup{width=1\textwidth, font=footnotesize,labelfont=footnotesize}
\centering\includegraphics[scale=0.2]{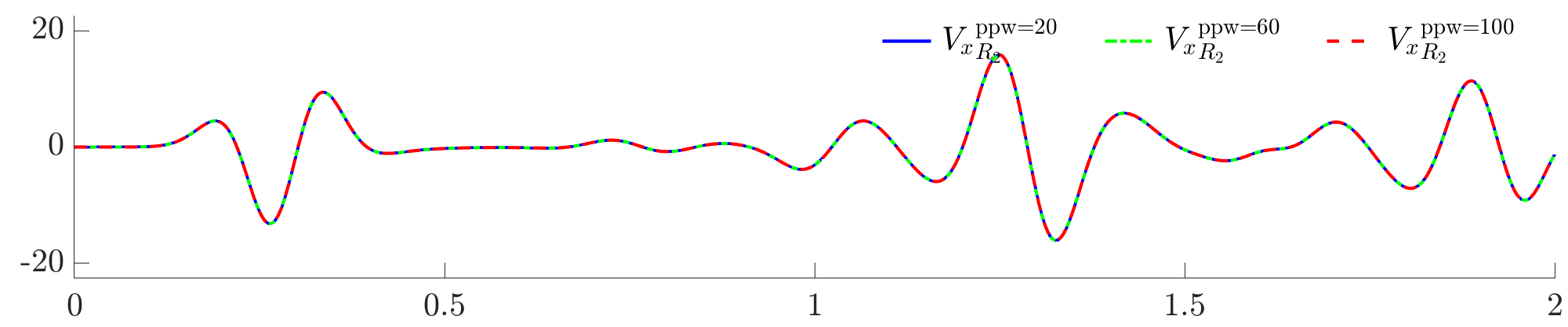}
\caption{Time history of $V_x$ at 2 grid points below the surface for the case ppw = 20.}
\label{stress_test_Vy_elastic_full_FSBC_Vx_2_weak}
\end{subfigure}\hfill
\caption{Time histories of the 2D elastic experiments. The free surface boundary condition is imposed weakly.
}
\label{stress_test_Vy_elastic_full_FSBC_Vx_weak}
\end{figure}
\end{minipage}
\ \\
\ \\

\hspace{-0.0625\textwidth}
\begin{minipage}[t]{.5\textwidth}
\begin{figure}[H]
\captionsetup{width=1\textwidth, font=footnotesize,labelfont=footnotesize}
\centering
\begin{subfigure}[b]{1\textwidth}
\captionsetup{width=1\textwidth, font=footnotesize,labelfont=footnotesize}
\centering\includegraphics[scale=0.2]{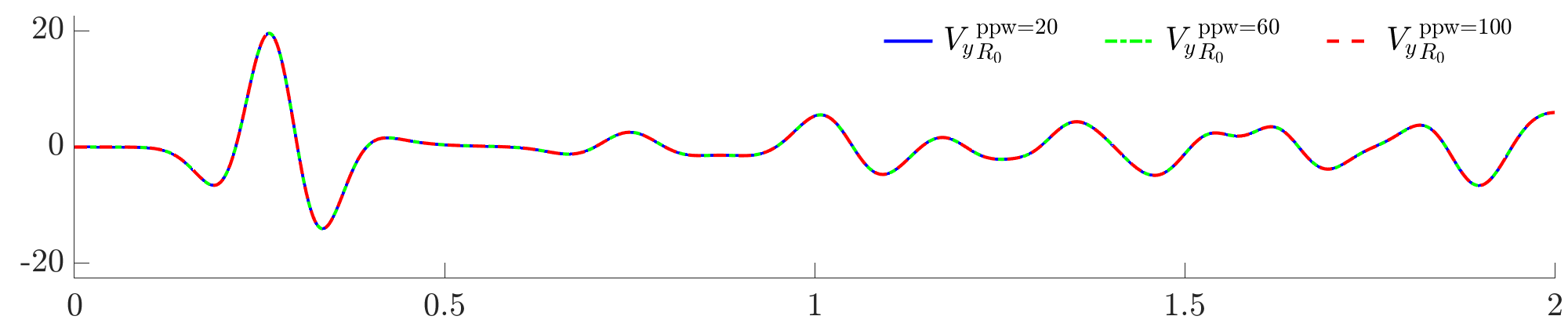}
\caption{Time history of $V_y$ at the surface.}
\label{stress_test_Vy_elastic_full_FSBC_Vy_0_strong}
\end{subfigure}\hfill
\\[2ex]
\begin{subfigure}[b]{1\textwidth}
\captionsetup{width=1\textwidth, font=footnotesize,labelfont=footnotesize}
\centering\includegraphics[scale=0.2]{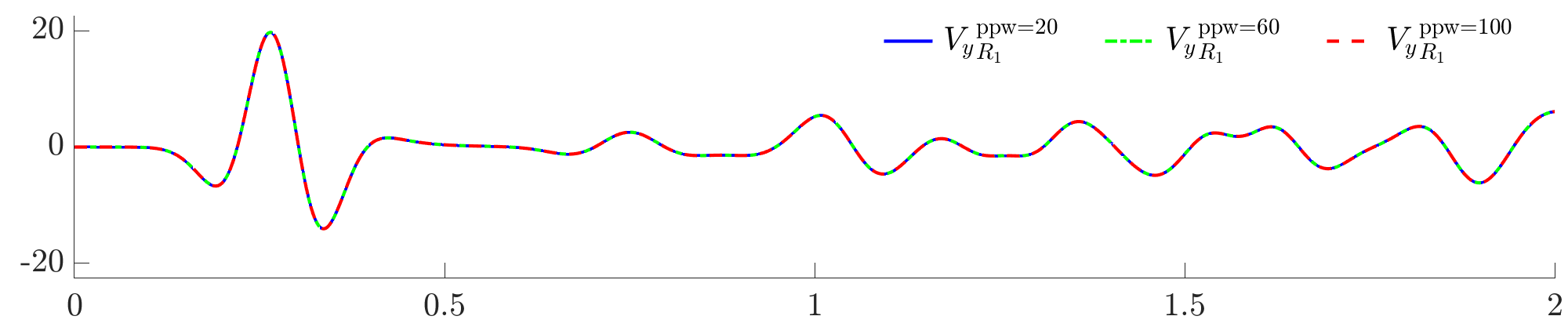}
\caption{Time history of $V_y$ at 1 grid points below the surface for the case ppw = 20.}
\label{stress_test_Vy_elastic_full_FSBC_Vy_1_strong}
\end{subfigure}\hfill
\\[2ex]
\begin{subfigure}[b]{1\textwidth}
\captionsetup{width=1\textwidth, font=footnotesize,labelfont=footnotesize}
\centering\includegraphics[scale=0.2]{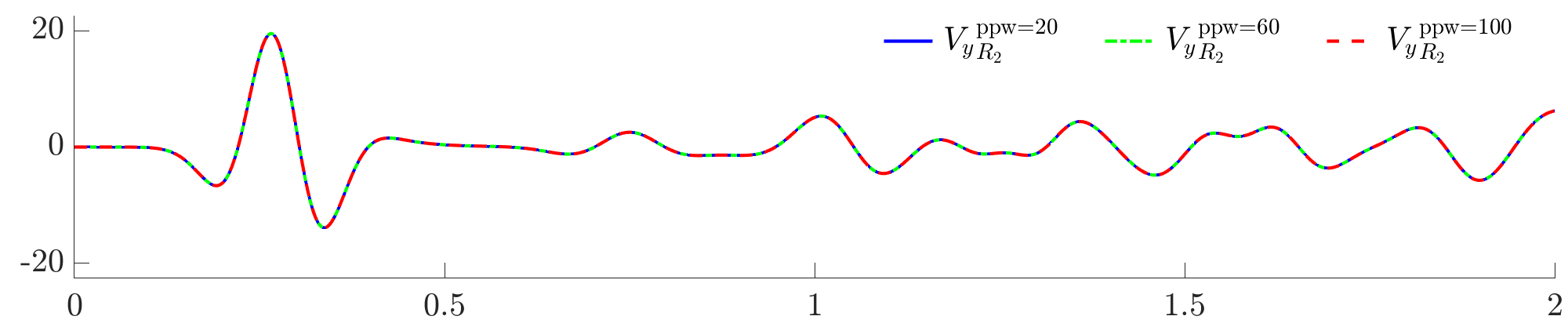}
\caption{Time history of $V_y$ at 2 grid points below the surface for the case ppw = 20.}
\label{stress_test_Vy_elastic_full_FSBC_Vy_2_strong}
\end{subfigure}\hfill
\caption{Time histories of the 2D elastic experiments. The free surface boundary condition is imposed strongly.
}
\label{stress_test_Vy_elastic_full_FSBC_Vy_strong}
\end{figure}
\end{minipage}
\hspace{0.0075\textwidth}
\begin{minipage}[t]{.5\textwidth}
\begin{figure}[H]
\captionsetup{width=1\textwidth, font=footnotesize,labelfont=footnotesize}
\centering
\begin{subfigure}[b]{1\textwidth}
\captionsetup{width=1\textwidth, font=footnotesize,labelfont=footnotesize}
\centering\includegraphics[scale=0.2]{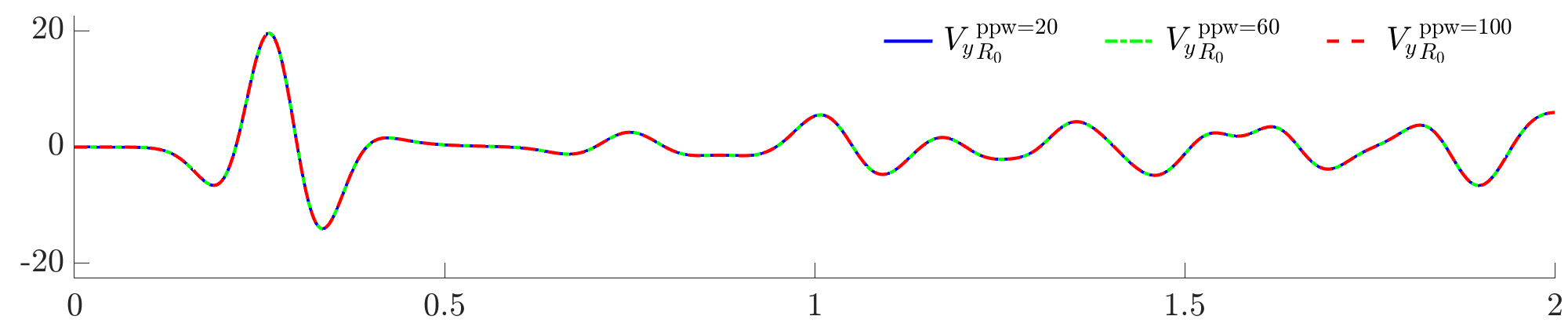}
\caption{Time history of $V_y$ at the surface.}
\label{stress_test_Vy_elastic_full_FSBC_Vy_0_weak}
\end{subfigure}\hfill
\\[2ex]
\begin{subfigure}[b]{1\textwidth}
\captionsetup{width=1\textwidth, font=footnotesize,labelfont=footnotesize}
\centering\includegraphics[scale=0.2]{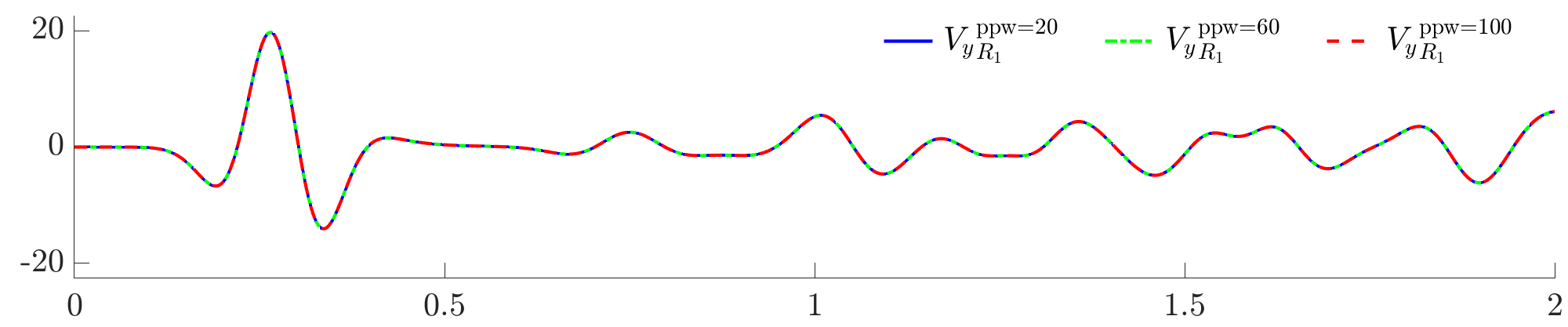}
\caption{Time history of $V_y$ at 1 grid points below the surface for the case ppw = 20.}
\label{stress_test_Vy_elastic_full_FSBC_Vy_1_weak}
\end{subfigure}\hfill
\\[2ex]
\begin{subfigure}[b]{1\textwidth}
\captionsetup{width=1\textwidth, font=footnotesize,labelfont=footnotesize}
\centering\includegraphics[scale=0.2]{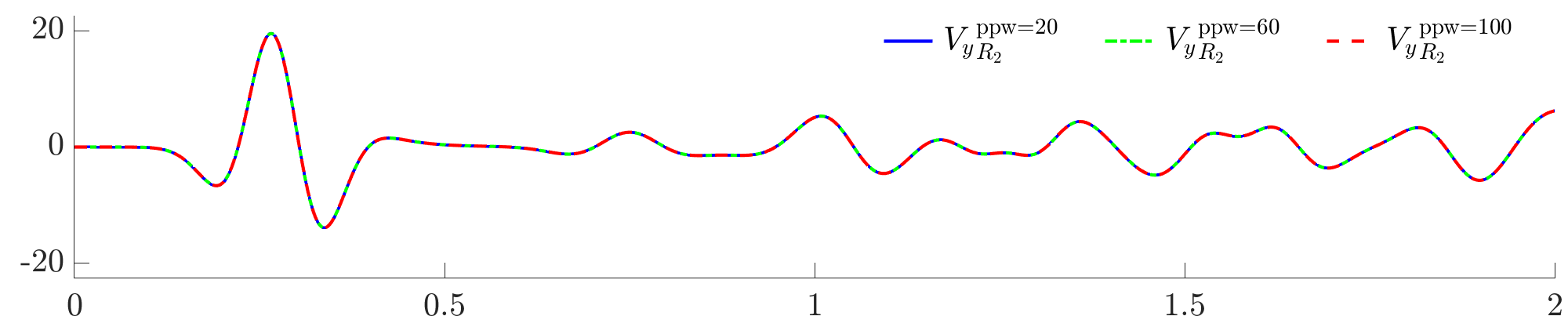}
\caption{Time history of $V_y$ at 2 grid points below the surface for the case ppw = 20.}
\label{stress_test_Vy_elastic_full_FSBC_Vy_2_weak}
\end{subfigure}\hfill
\caption{Time histories of the 2D elastic experiments. The free surface boundary condition is imposed weakly.
}
\label{stress_test_Vy_elastic_full_FSBC_Vy_weak}
\end{figure}
\end{minipage}

%% file: input_source_on_Sxy_corner.tex
\hspace{-0.0625\textwidth}
\begin{minipage}[t]{.5\textwidth}
\begin{figure}[H]
\captionsetup{width=1\textwidth, font=footnotesize,labelfont=footnotesize}
\centering
\begin{subfigure}[b]{1\textwidth}
\captionsetup{width=1\textwidth, font=footnotesize,labelfont=footnotesize}
\centering\includegraphics[scale=0.2]{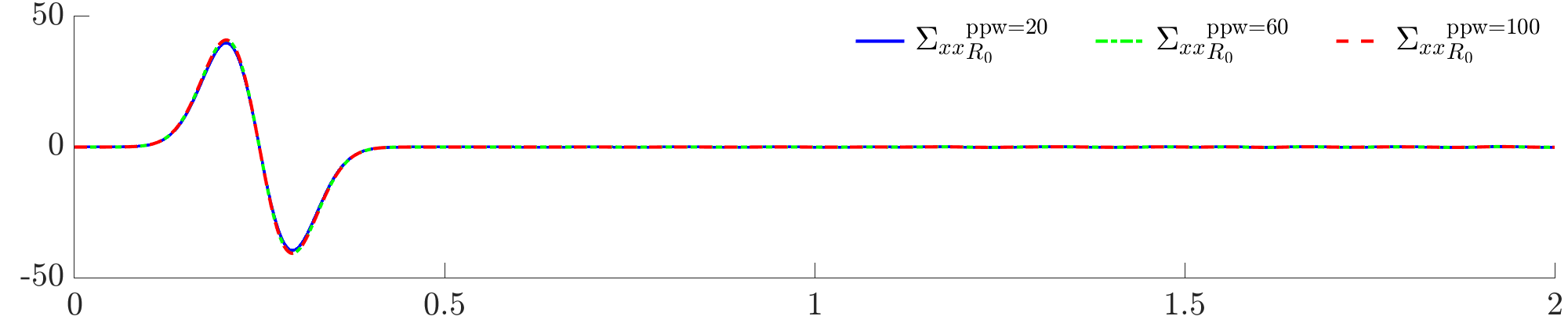}
\caption{Time history of $\Sigma_{xx}$ at the surface.}
\label{stress_test_Sxy_Elastic_full_FSBC_Sxx_0_strong}
\end{subfigure}\hfill
\\[2ex]
\begin{subfigure}[b]{1\textwidth}
\captionsetup{width=1\textwidth, font=footnotesize,labelfont=footnotesize}
\centering\includegraphics[scale=0.2]{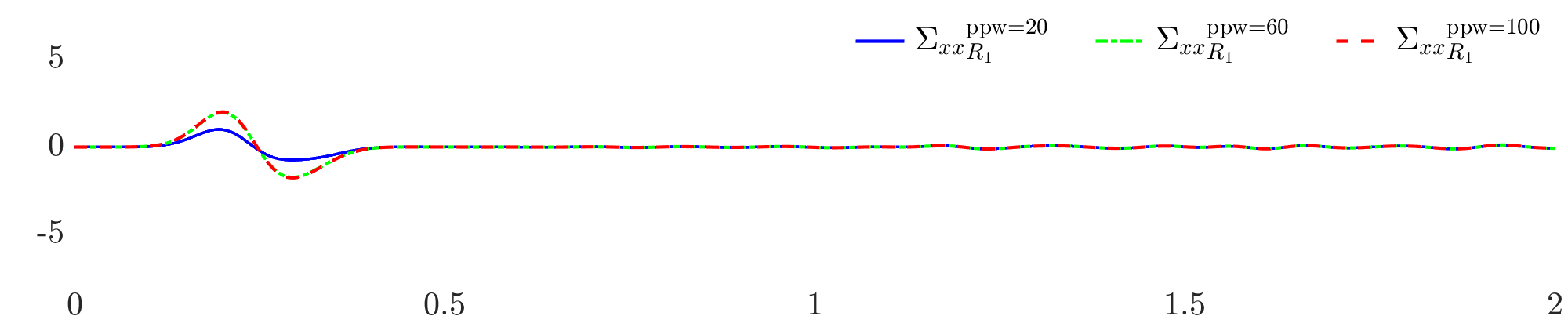}
\caption{Time history of $\Sigma_{xx}$ at 1 grid points below the surface for the case ppw = 20.}
\label{stress_test_Sxy_Elastic_full_FSBC_Sxx_1_strong}
\end{subfigure}\hfill
\\[2ex]
\begin{subfigure}[b]{1\textwidth}
\captionsetup{width=1\textwidth, font=footnotesize,labelfont=footnotesize}
\centering\includegraphics[scale=0.2]{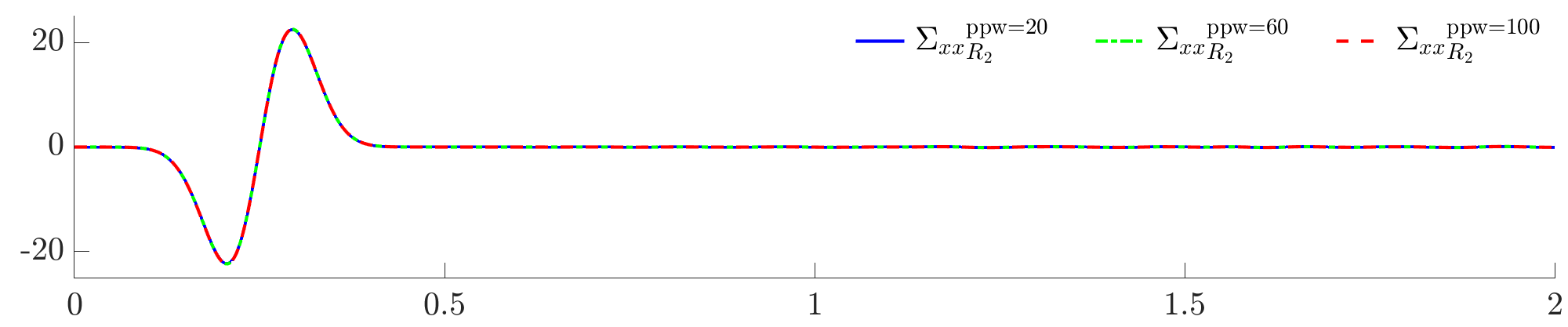}
\caption{Time history of $\Sigma_{xx}$ at 2 grid points below the surface for the case ppw = 20.}
\label{stress_test_Sxy_Elastic_full_FSBC_Sxx_2_strong}
\end{subfigure}\hfill
\caption{Time histories of the 2D elastic experiments. The free surface boundary condition is imposed strongly.
}
\label{stress_test_Sxy_Elastic_full_FSBC_Sxx_strong}
\end{figure}
\end{minipage}
\hspace{0.0075\textwidth}
\begin{minipage}[t]{.5\textwidth}
\begin{figure}[H]
\captionsetup{width=1\textwidth, font=footnotesize,labelfont=footnotesize}
\centering
\begin{subfigure}[b]{1\textwidth}
\captionsetup{width=1\textwidth, font=footnotesize,labelfont=footnotesize}
\centering\includegraphics[scale=0.2]{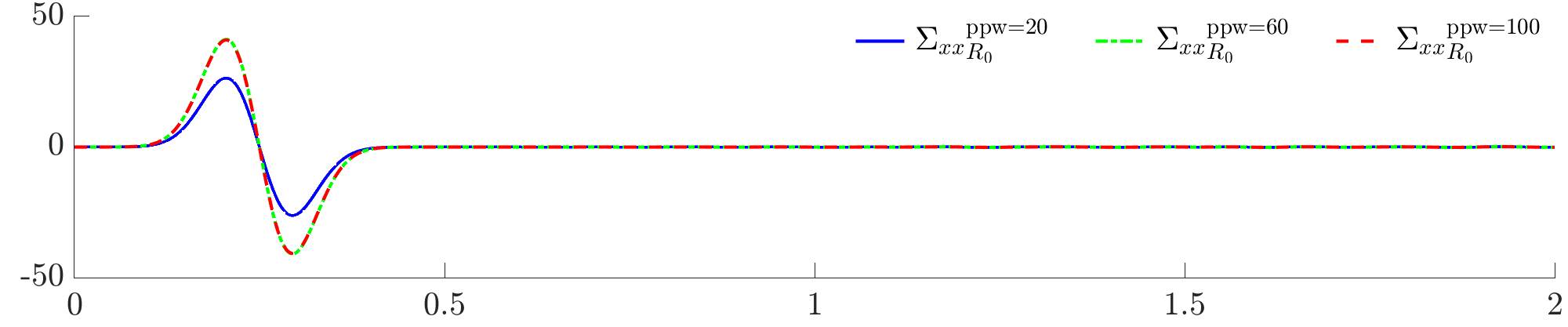}
\caption{Time history of $\Sigma_{xx}$ at the surface.}
\label{stress_test_Sxy_Elastic_full_FSBC_Sxx_0_weak}
\end{subfigure}\hfill
\\[2ex]
\begin{subfigure}[b]{1\textwidth}
\captionsetup{width=1\textwidth, font=footnotesize,labelfont=footnotesize}
\centering\includegraphics[scale=0.2]{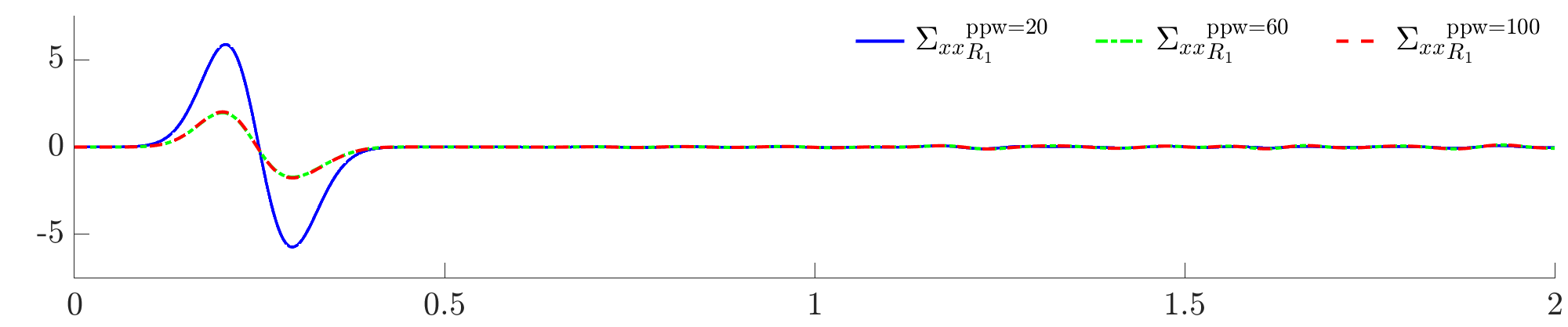}
\caption{Time history of $\Sigma_{xx}$ at 1 grid points below the surface for the case ppw = 20.}
\label{stress_test_Sxy_Elastic_full_FSBC_Sxx_1_weak}
\end{subfigure}\hfill
\\[2ex]
\begin{subfigure}[b]{1\textwidth}
\captionsetup{width=1\textwidth, font=footnotesize,labelfont=footnotesize}
\centering\includegraphics[scale=0.2]{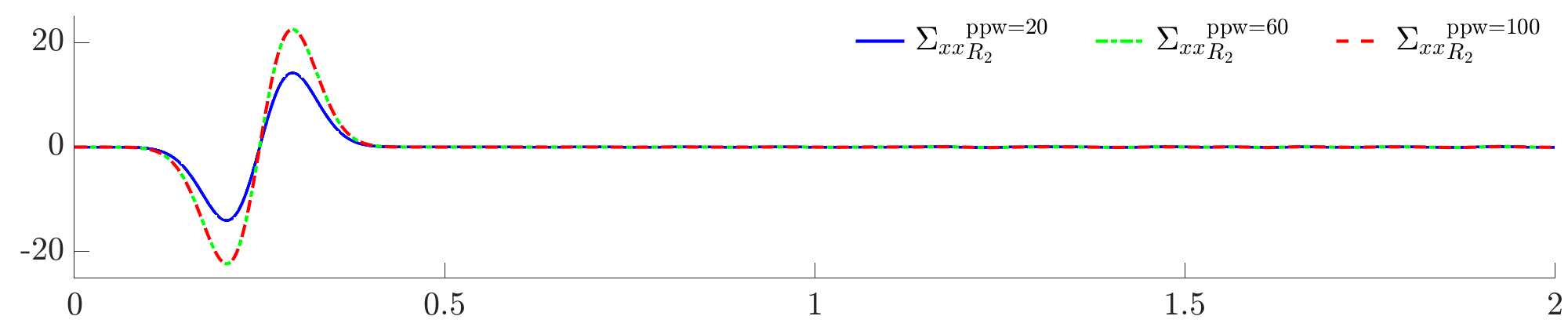}
\caption{Time history of $\Sigma_{xx}$ at 2 grid points below the surface for the case ppw = 20.}
\label{stress_test_Sxy_Elastic_full_FSBC_Sxx_2_weak}
\end{subfigure}\hfill
\caption{Time histories of the 2D elastic experiments. The free surface boundary condition is imposed weakly.
}
\label{stress_test_Sxy_Elastic_full_FSBC_Sxx_weak}
\end{figure}
\end{minipage}
\ \\
\ \\

\hspace{-0.0625\textwidth}
\begin{minipage}[t]{.5\textwidth}
\begin{figure}[H]
\captionsetup{width=1\textwidth, font=footnotesize,labelfont=footnotesize}
\centering
\begin{subfigure}[b]{1\textwidth}
\captionsetup{width=1\textwidth, font=footnotesize,labelfont=footnotesize}
\centering\includegraphics[scale=0.2]{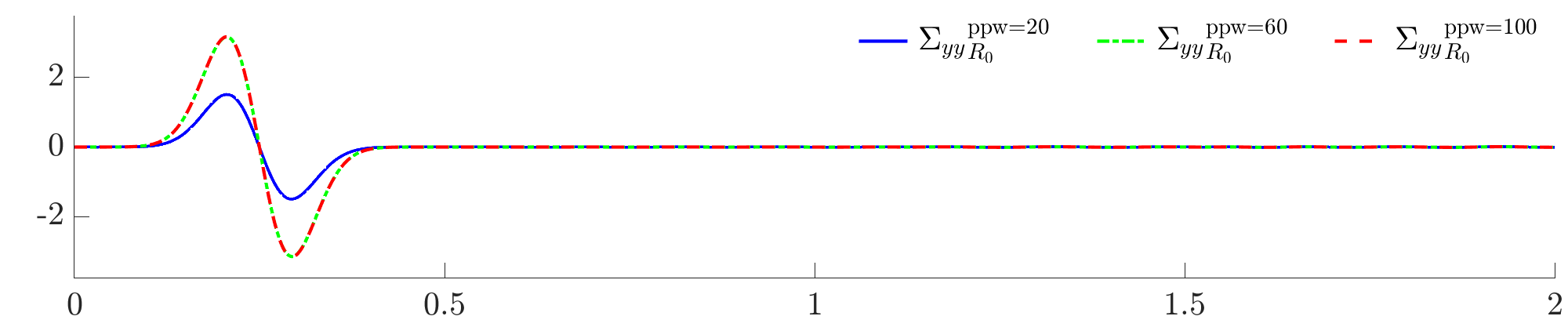}
\caption{Time history of $\Sigma_{yy}$ at the surface.}
\label{stress_test_Sxy_Elastic_full_FSBC_Syy_0_strong}
\end{subfigure}\hfill
\\[2ex]
\begin{subfigure}[b]{1\textwidth}
\captionsetup{width=1\textwidth, font=footnotesize,labelfont=footnotesize}
\centering\includegraphics[scale=0.2]{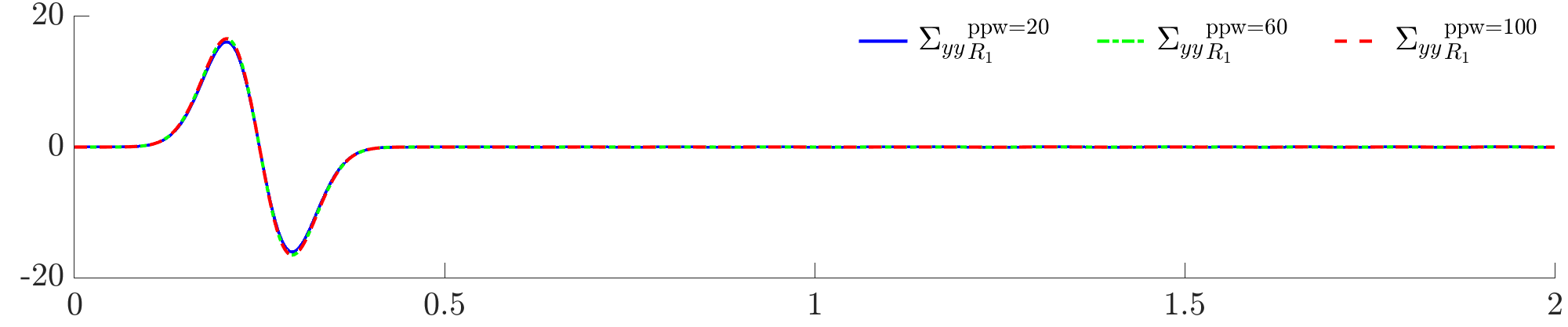}
\caption{Time history of $\Sigma_{yy}$ at 1 grid points below the surface for the case ppw = 20.}
\label{stress_test_Sxy_Elastic_full_FSBC_Syy_1_strong}
\end{subfigure}\hfill
\\[2ex]
\begin{subfigure}[b]{1\textwidth}
\captionsetup{width=1\textwidth, font=footnotesize,labelfont=footnotesize}
\centering\includegraphics[scale=0.2]{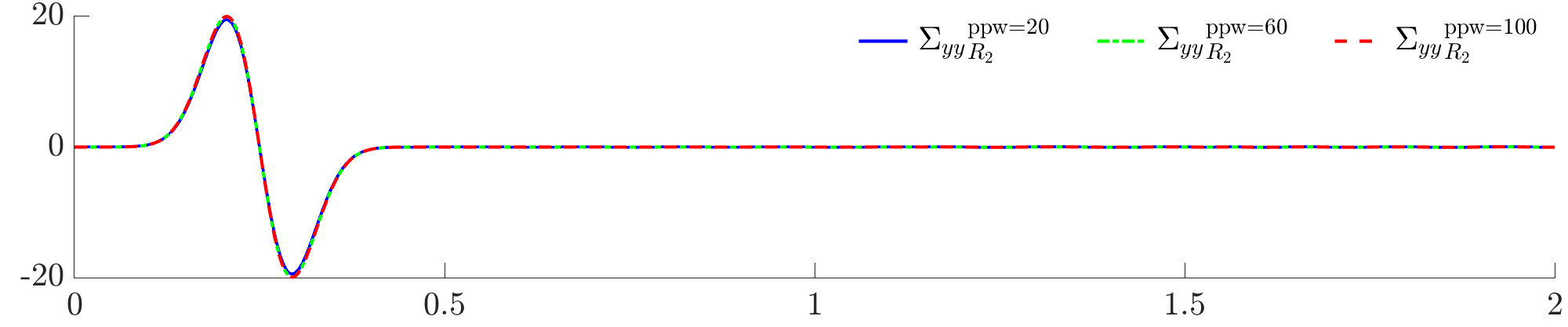}
\caption{Time history of $\Sigma_{yy}$ at 2 grid points below the surface for the case ppw = 20.}
\label{stress_test_Sxy_Elastic_full_FSBC_Syy_2_strong}
\end{subfigure}\hfill
\caption{Time histories of the 2D elastic experiments. The free surface boundary condition is imposed strongly.
}
\label{stress_test_Sxy_Elastic_full_FSBC_Syy_strong}
\end{figure}
\end{minipage}
\hspace{0.0075\textwidth}
\begin{minipage}[t]{.5\textwidth}
\begin{figure}[H]
\captionsetup{width=1\textwidth, font=footnotesize,labelfont=footnotesize}
\centering
\begin{subfigure}[b]{1\textwidth}
\captionsetup{width=1\textwidth, font=footnotesize,labelfont=footnotesize}
\centering\includegraphics[scale=0.2]{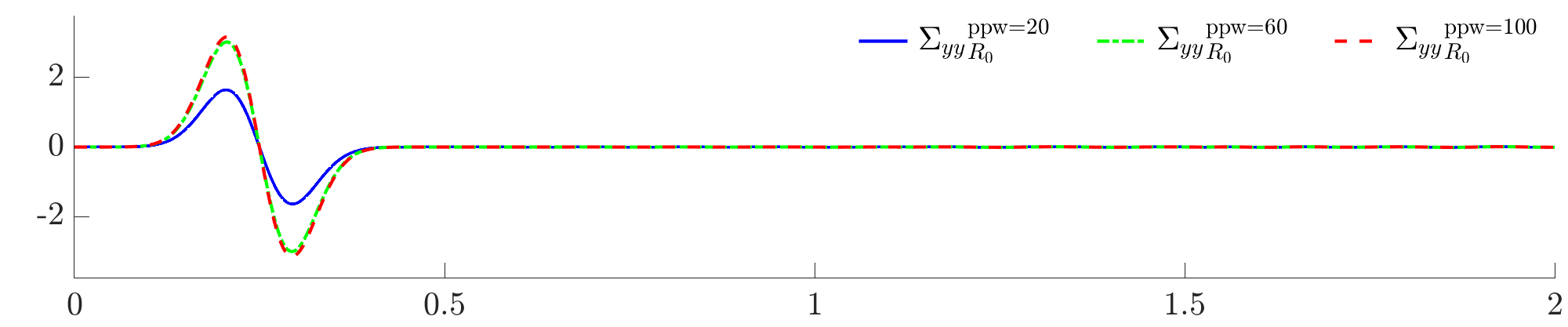}
\caption{Time history of $\Sigma_{yy}$ at the surface.}
\label{stress_test_Sxy_Elastic_full_FSBC_Syy_0_weak}
\end{subfigure}\hfill
\\[2ex]
\begin{subfigure}[b]{1\textwidth}
\captionsetup{width=1\textwidth, font=footnotesize,labelfont=footnotesize}
\centering\includegraphics[scale=0.2]{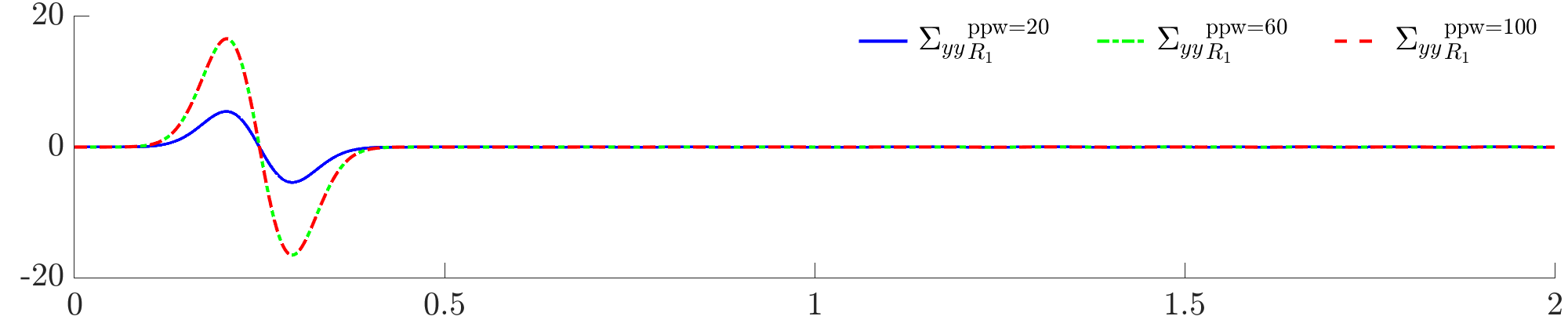}
\caption{Time history of $\Sigma_{yy}$ at 1 grid points below the surface for the case ppw = 20.}
\label{stress_test_Sxy_Elastic_full_FSBC_Syy_1_weak}
\end{subfigure}\hfill
\\[2ex]
\begin{subfigure}[b]{1\textwidth}
\captionsetup{width=1\textwidth, font=footnotesize,labelfont=footnotesize}
\centering\includegraphics[scale=0.2]{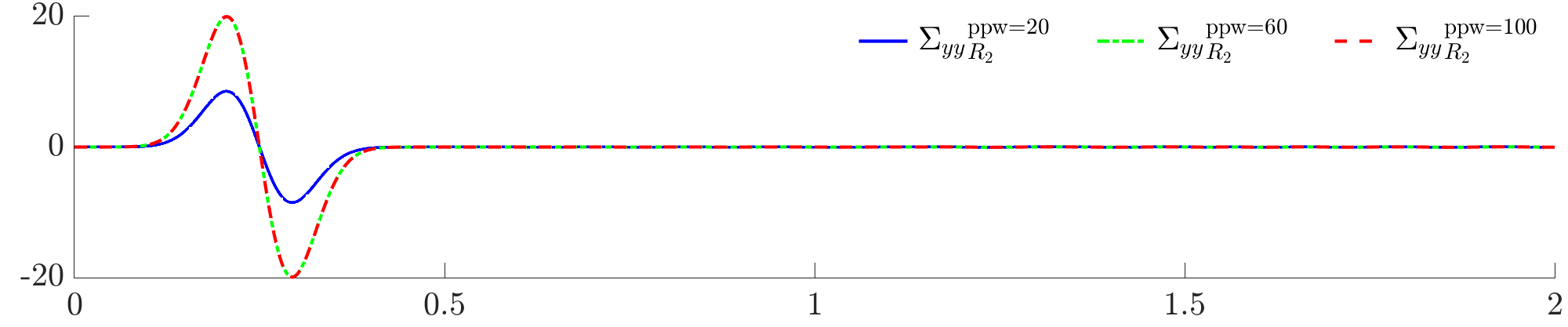}
\caption{Time history of $\Sigma_{yy}$ at 2 grid points below the surface for the case ppw = 20.}
\label{stress_test_Sxy_Elastic_full_FSBC_Syy_2_weak}
\end{subfigure}\hfill
\caption{Time histories of the 2D elastic experiments. The free surface boundary condition is imposed weakly.
}
\label{stress_test_Sxy_Elastic_full_FSBC_Syy_weak}
\end{figure}
\end{minipage}
\ \\
\ \\

\hspace{-0.0625\textwidth}
\begin{minipage}[t]{.5\textwidth}
\begin{figure}[H]
\captionsetup{width=1\textwidth, font=footnotesize,labelfont=footnotesize}
\centering
\begin{subfigure}[b]{1\textwidth}
\captionsetup{width=1\textwidth, font=footnotesize,labelfont=footnotesize}
\centering\includegraphics[scale=0.2]{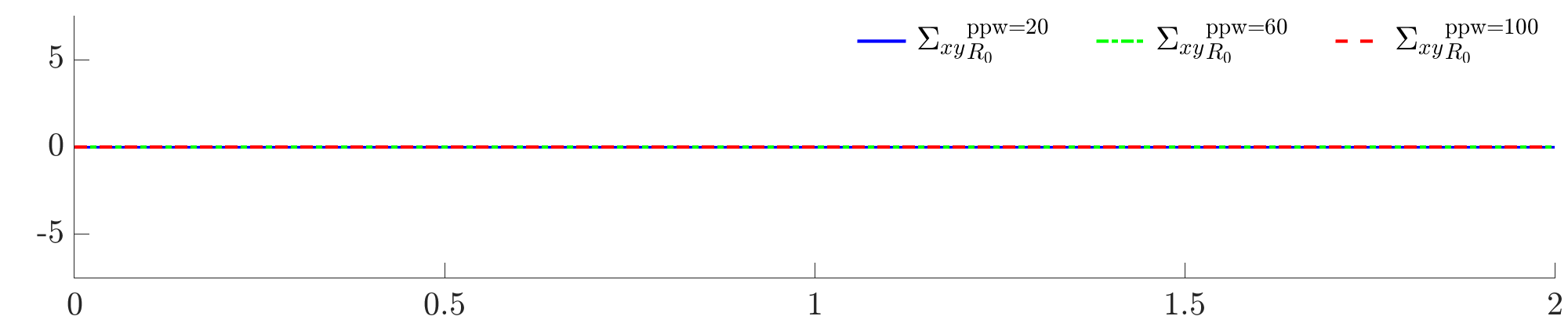}
\caption{Time history of $\Sigma_{xy}$ at the surface.}
\label{stress_test_Sxy_Elastic_full_FSBC_Sxy_0_strong}
\end{subfigure}\hfill
\\[2ex]
\begin{subfigure}[b]{1\textwidth}
\captionsetup{width=1\textwidth, font=footnotesize,labelfont=footnotesize}
\centering\includegraphics[scale=0.2]{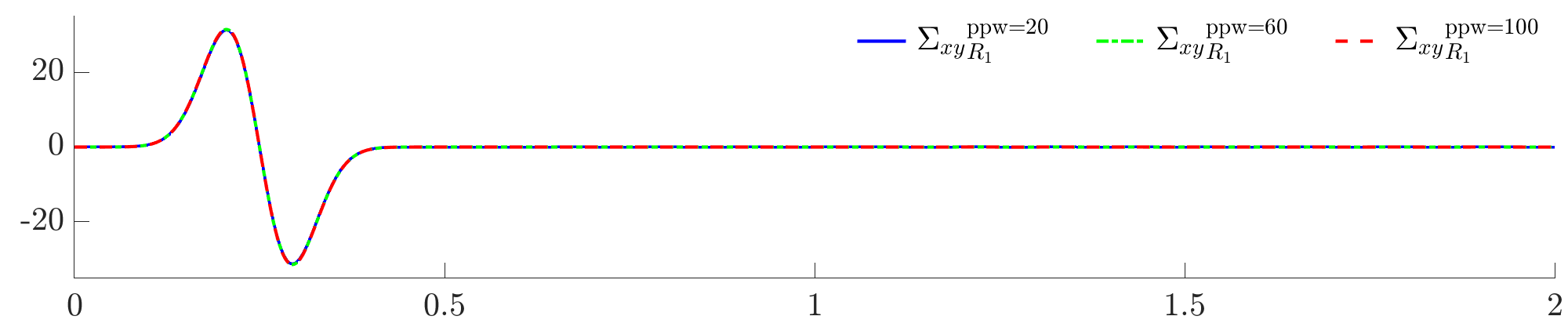}
\caption{Time history of $\Sigma_{xy}$ at 1 grid points below the surface for the case ppw = 20.}
\label{stress_test_Sxy_Elastic_full_FSBC_Sxy_1_strong}
\end{subfigure}\hfill
\\[2ex]
\begin{subfigure}[b]{1\textwidth}
\captionsetup{width=1\textwidth, font=footnotesize,labelfont=footnotesize}
\centering\includegraphics[scale=0.2]{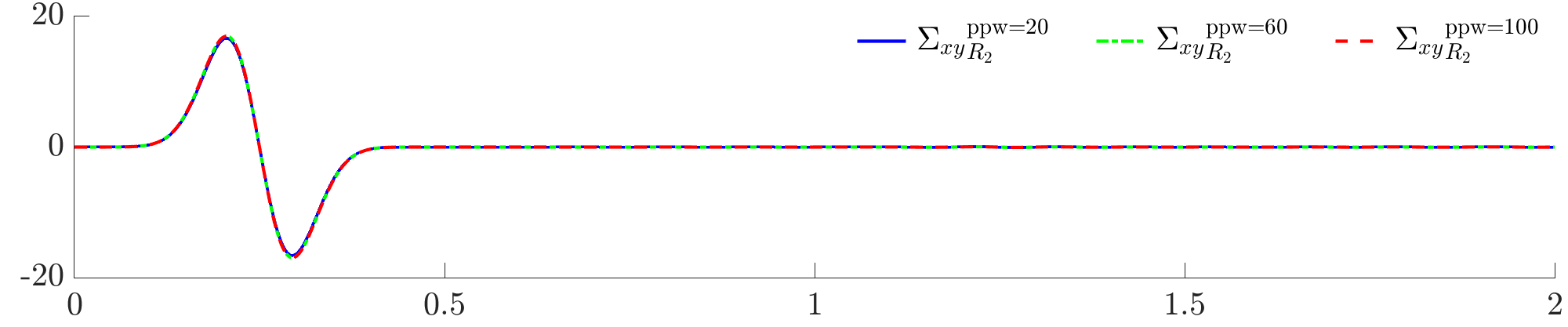}
\caption{Time history of $\Sigma_{xy}$ at 2 grid points below the surface for the case ppw = 20.}
\label{stress_test_Sxy_Elastic_full_FSBC_Sxy_2_strong}
\end{subfigure}\hfill
\caption{Time histories of the 2D elastic experiments. The free surface boundary condition is imposed strongly.
}
\label{stress_test_Sxy_Elastic_full_FSBC_Sxy_strong}
\end{figure}
\end{minipage}
\hspace{0.0075\textwidth}
\begin{minipage}[t]{.5\textwidth}
\begin{figure}[H]
\captionsetup{width=1\textwidth, font=footnotesize,labelfont=footnotesize}
\centering
\begin{subfigure}[b]{1\textwidth}
\captionsetup{width=1\textwidth, font=footnotesize,labelfont=footnotesize}
\centering\includegraphics[scale=0.2]{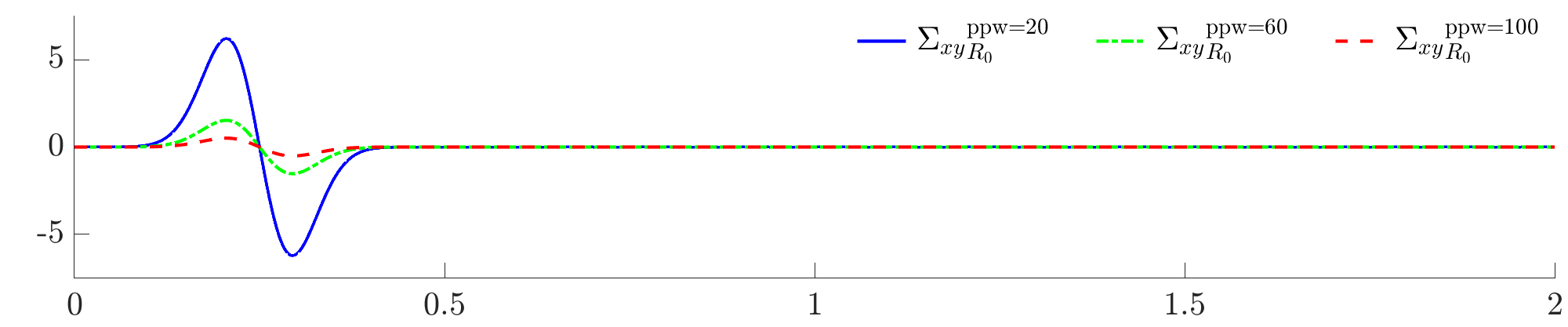}
\caption{Time history of $\Sigma_{xy}$ at the surface.}
\label{stress_test_Sxy_Elastic_full_FSBC_Sxy_0_weak}
\end{subfigure}\hfill
\\[2ex]
\begin{subfigure}[b]{1\textwidth}
\captionsetup{width=1\textwidth, font=footnotesize,labelfont=footnotesize}
\centering\includegraphics[scale=0.2]{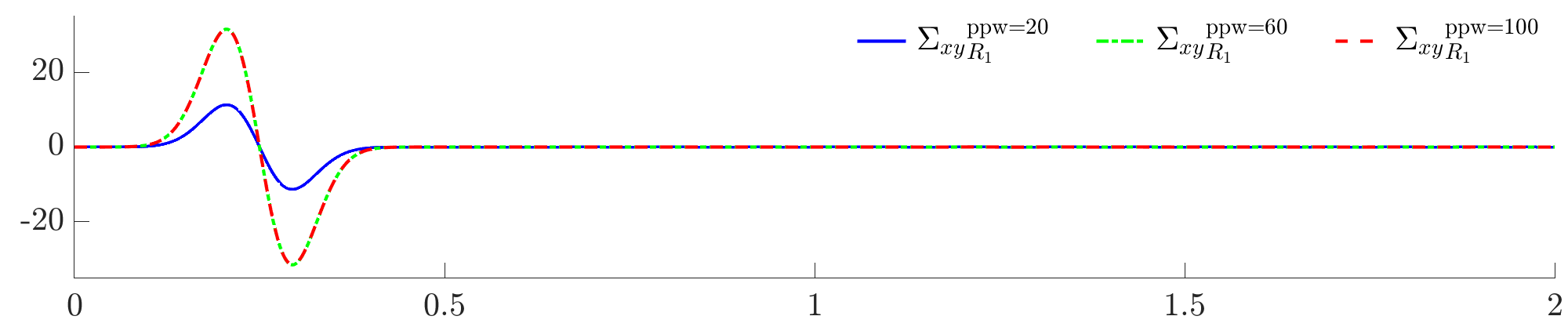}
\caption{Time history of $\Sigma_{xy}$ at 1 grid points below the surface for the case ppw = 20.}
\label{stress_test_Sxy_Elastic_full_FSBC_Sxy_1_weak}
\end{subfigure}\hfill
\\[2ex]
\begin{subfigure}[b]{1\textwidth}
\captionsetup{width=1\textwidth, font=footnotesize,labelfont=footnotesize}
\centering\includegraphics[scale=0.2]{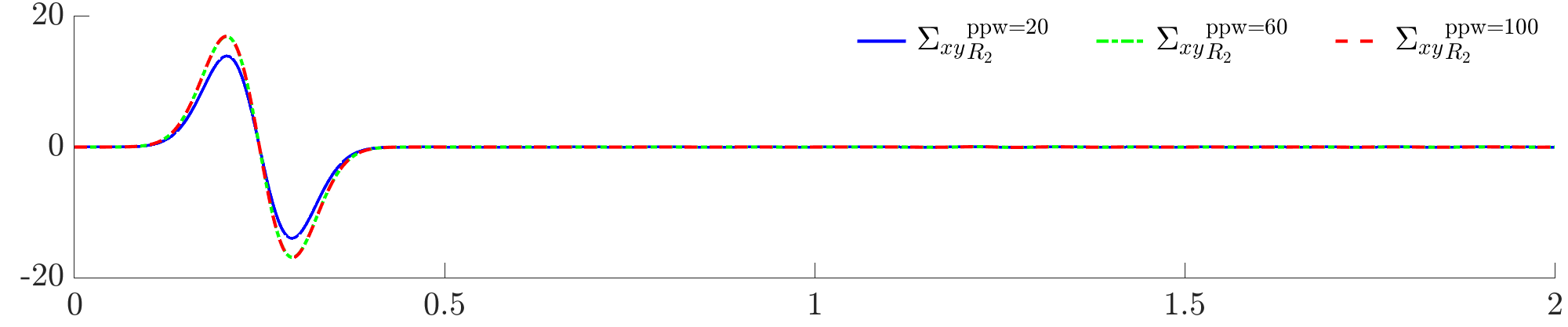}
\caption{Time history of $\Sigma_{xy}$ at 2 grid points below the surface for the case ppw = 20.}
\label{stress_test_Sxy_Elastic_full_FSBC_Sxy_2_weak}
\end{subfigure}\hfill
\caption{Time histories of the 2D elastic experiments. The free surface boundary condition is imposed weakly.
}
\label{stress_test_Sxy_Elastic_full_FSBC_Sxy_weak}
\end{figure}
\end{minipage}
\ \\
\ \\

\hspace{-0.0625\textwidth}
\begin{minipage}[t]{.5\textwidth}
\begin{figure}[H]
\captionsetup{width=1\textwidth, font=footnotesize,labelfont=footnotesize}
\centering
\begin{subfigure}[b]{1\textwidth}
\captionsetup{width=1\textwidth, font=footnotesize,labelfont=footnotesize}
\centering\includegraphics[scale=0.2]{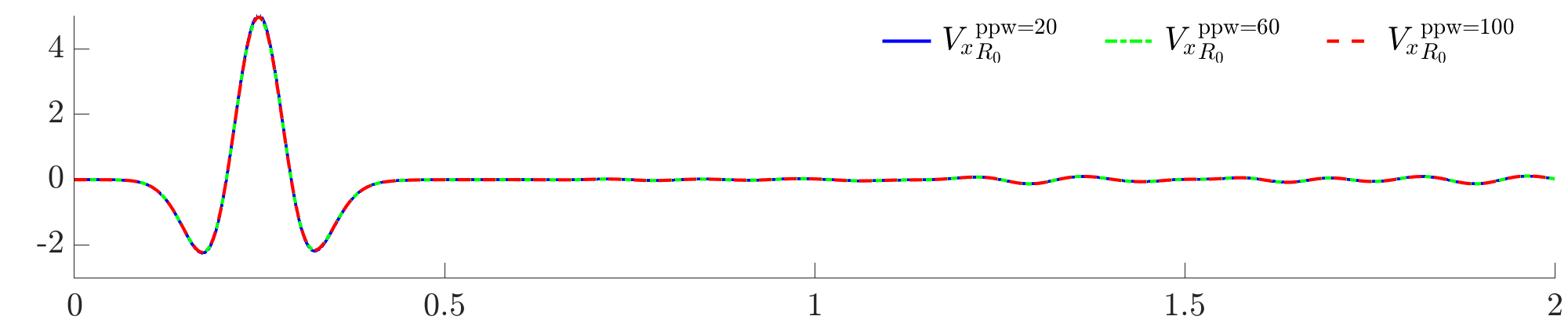}
\caption{Time history of $V_x$ at the surface.}
\label{stress_test_Sxy_Elastic_full_FSBC_Vx_0_strong}
\end{subfigure}\hfill
\\[2ex]
\begin{subfigure}[b]{1\textwidth}
\captionsetup{width=1\textwidth, font=footnotesize,labelfont=footnotesize}
\centering\includegraphics[scale=0.2]{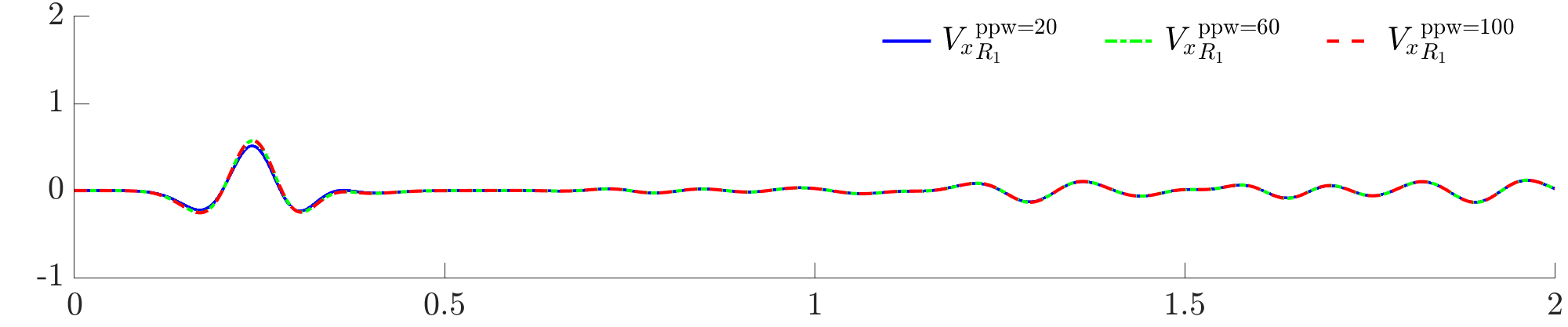}
\caption{Time history of $V_x$ at 1 grid points below the surface for the case ppw = 20.}
\label{stress_test_Sxy_Elastic_full_FSBC_Vx_1_strong}
\end{subfigure}\hfill
\\[2ex]
\begin{subfigure}[b]{1\textwidth}
\captionsetup{width=1\textwidth, font=footnotesize,labelfont=footnotesize}
\centering\includegraphics[scale=0.2]{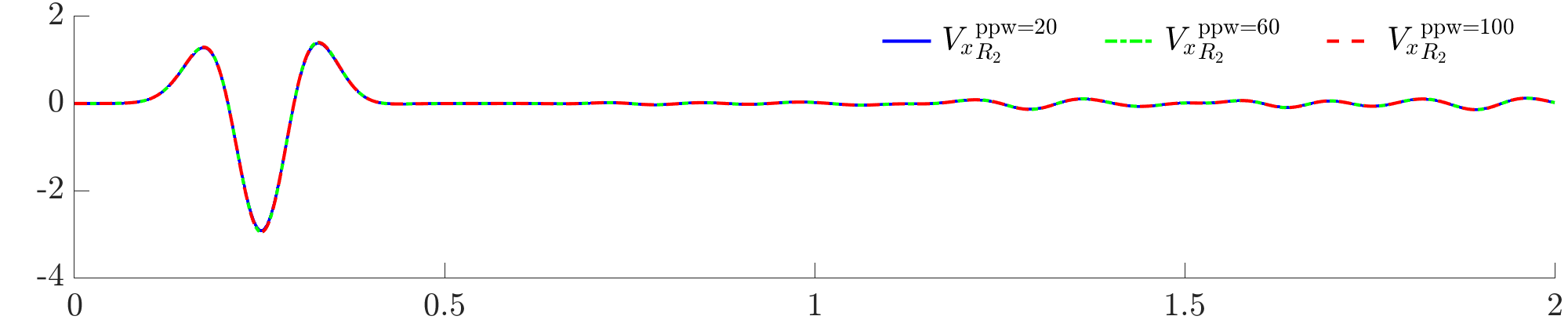}
\caption{Time history of $V_x$ at 2 grid points below the surface for the case ppw = 20.}
\label{stress_test_Sxy_Elastic_full_FSBC_Vx_2_strong}
\end{subfigure}\hfill
\caption{Time histories of the 2D elastic experiments. The free surface boundary condition is imposed strongly.
}
\label{stress_test_Sxy_Elastic_full_FSBC_Vx_strong}
\end{figure}
\end{minipage}
\hspace{0.0075\textwidth}
\begin{minipage}[t]{.5\textwidth}
\begin{figure}[H]
\captionsetup{width=1\textwidth, font=footnotesize,labelfont=footnotesize}
\centering
\begin{subfigure}[b]{1\textwidth}
\captionsetup{width=1\textwidth, font=footnotesize,labelfont=footnotesize}
\centering\includegraphics[scale=0.2]{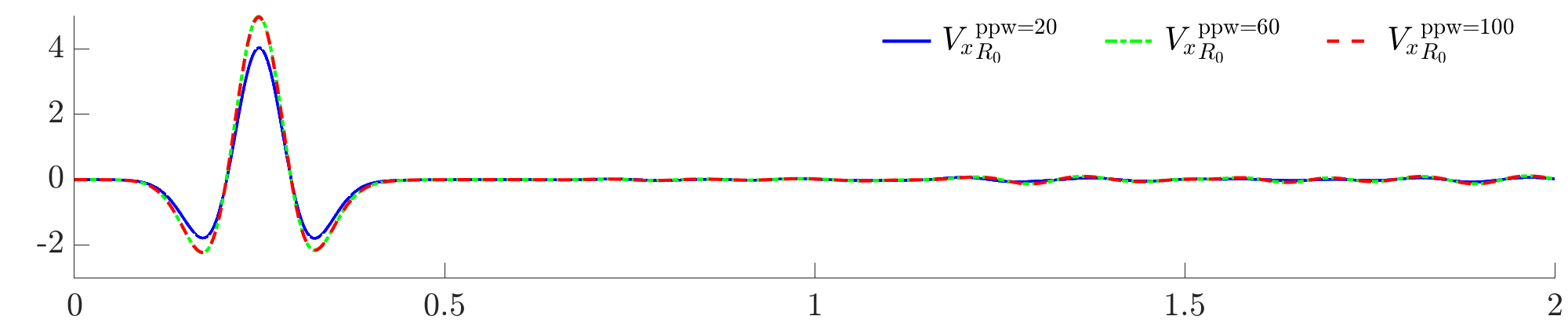}
\caption{Time history of $V_x$ at the surface.}
\label{stress_test_Sxy_Elastic_full_FSBC_Vx_0_weak}
\end{subfigure}\hfill
\\[2ex]
\begin{subfigure}[b]{1\textwidth}
\captionsetup{width=1\textwidth, font=footnotesize,labelfont=footnotesize}
\centering\includegraphics[scale=0.2]{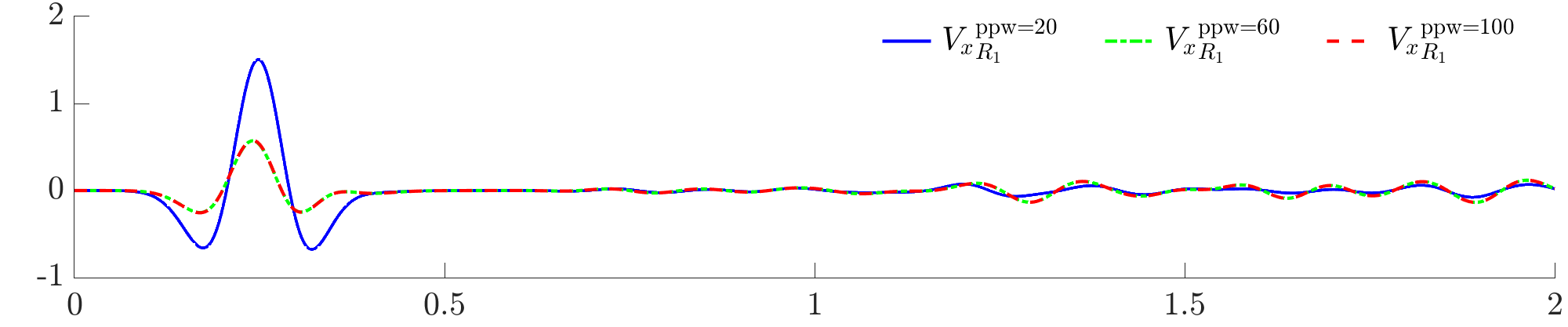}
\caption{Time history of $V_x$ at 1 grid points below the surface for the case ppw = 20.}
\label{stress_test_Sxy_Elastic_full_FSBC_Vx_1_weak}
\end{subfigure}\hfill
\\[2ex]
\begin{subfigure}[b]{1\textwidth}
\captionsetup{width=1\textwidth, font=footnotesize,labelfont=footnotesize}
\centering\includegraphics[scale=0.2]{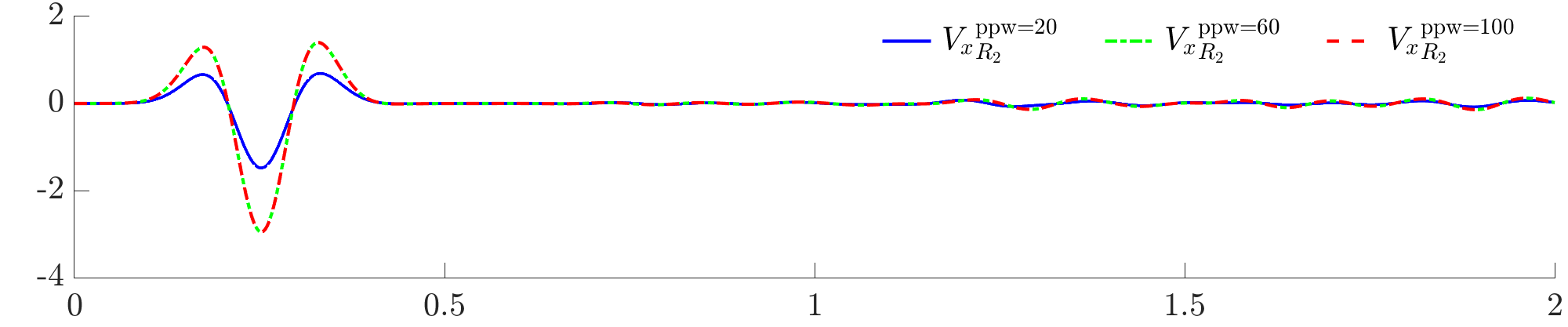}
\caption{Time history of $V_x$ at 2 grid points below the surface for the case ppw = 20.}
\label{stress_test_Sxy_Elastic_full_FSBC_Vx_2_weak}
\end{subfigure}\hfill
\caption{Time histories of the 2D elastic experiments. The free surface boundary condition is imposed weakly.
}
\label{stress_test_Sxy_Elastic_full_FSBC_Vx_weak}
\end{figure}
\end{minipage}
\ \\
\ \\

\hspace{-0.0625\textwidth}
\begin{minipage}[t]{.5\textwidth}
\begin{figure}[H]
\captionsetup{width=1\textwidth, font=footnotesize,labelfont=footnotesize}
\centering
\begin{subfigure}[b]{1\textwidth}
\captionsetup{width=1\textwidth, font=footnotesize,labelfont=footnotesize}
\centering\includegraphics[scale=0.2]{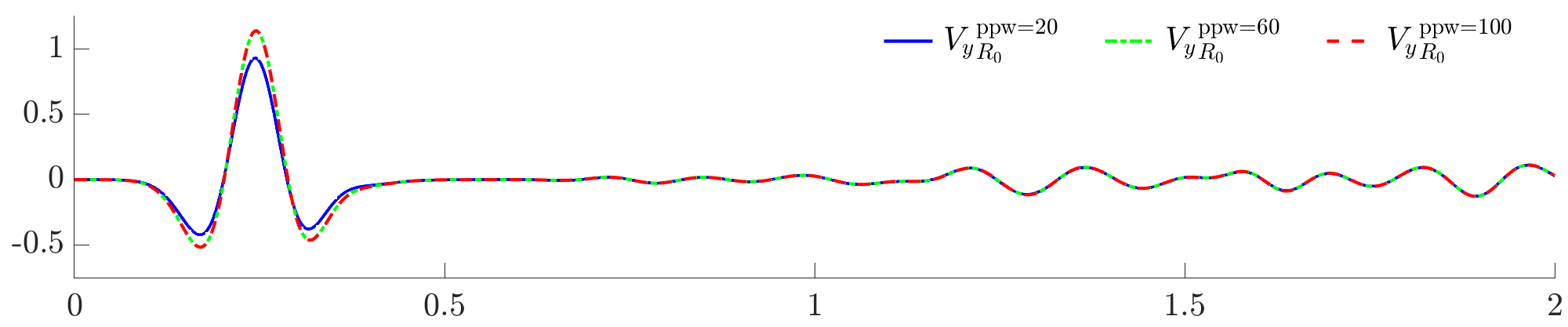}
\caption{Time history of $V_y$ at the surface.}
\label{stress_test_Sxy_Elastic_full_FSBC_Vy_0_strong}
\end{subfigure}\hfill
\\[2ex]
\begin{subfigure}[b]{1\textwidth}
\captionsetup{width=1\textwidth, font=footnotesize,labelfont=footnotesize}
\centering\includegraphics[scale=0.2]{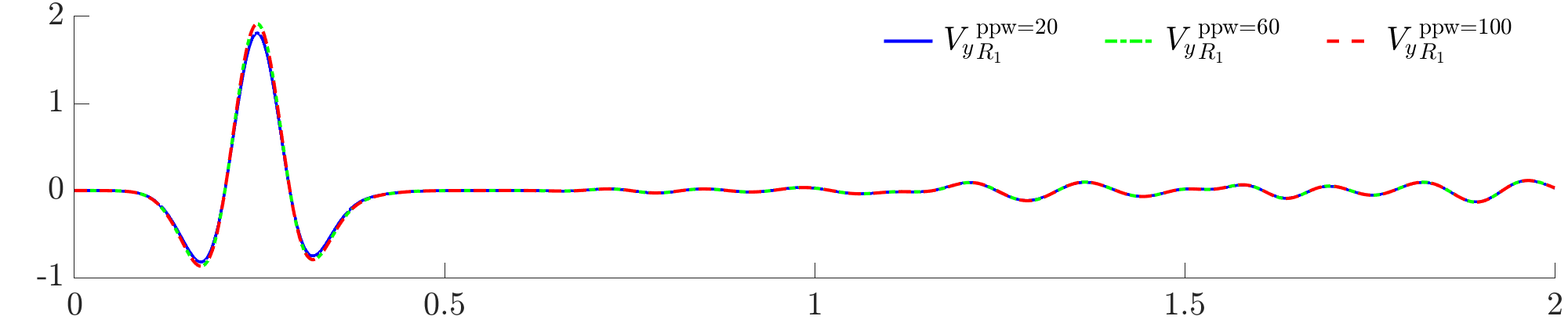}
\caption{Time history of $V_y$ at 1 grid points below the surface for the case ppw = 20.}
\label{stress_test_Sxy_Elastic_full_FSBC_Vy_1_strong}
\end{subfigure}\hfill
\\[2ex]
\begin{subfigure}[b]{1\textwidth}
\captionsetup{width=1\textwidth, font=footnotesize,labelfont=footnotesize}
\centering\includegraphics[scale=0.2]{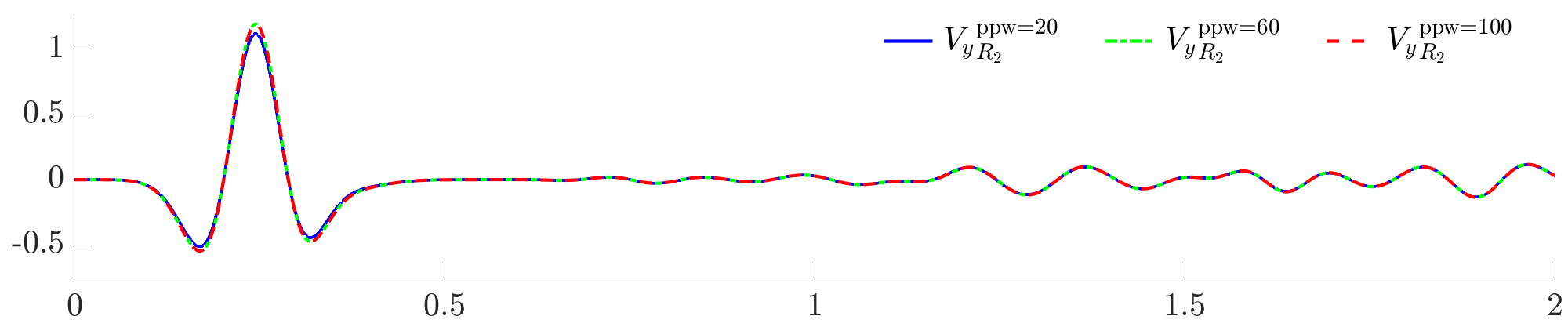}
\caption{Time history of $V_y$ at 2 grid points below the surface for the case ppw = 20.}
\label{stress_test_Sxy_Elastic_full_FSBC_Vy_2_strong}
\end{subfigure}\hfill
\caption{Time histories of the 2D elastic experiments. The free surface boundary condition is imposed strongly.
}
\label{stress_test_Sxy_Elastic_full_FSBC_Vy_strong}
\end{figure}
\end{minipage}
\hspace{0.0075\textwidth}
\begin{minipage}[t]{.5\textwidth}
\begin{figure}[H]
\captionsetup{width=1\textwidth, font=footnotesize,labelfont=footnotesize}
\centering
\begin{subfigure}[b]{1\textwidth}
\captionsetup{width=1\textwidth, font=footnotesize,labelfont=footnotesize}
\centering\includegraphics[scale=0.2]{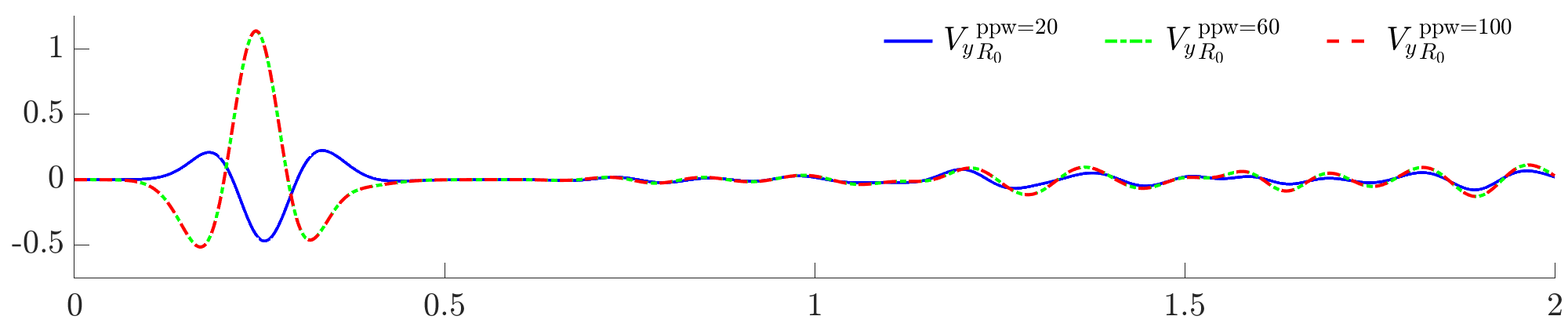}
\caption{Time history of $V_y$ at the surface.}
\label{stress_test_Sxy_Elastic_full_FSBC_Vy_0_weak}
\end{subfigure}\hfill
\\[2ex]
\begin{subfigure}[b]{1\textwidth}
\captionsetup{width=1\textwidth, font=footnotesize,labelfont=footnotesize}
\centering\includegraphics[scale=0.2]{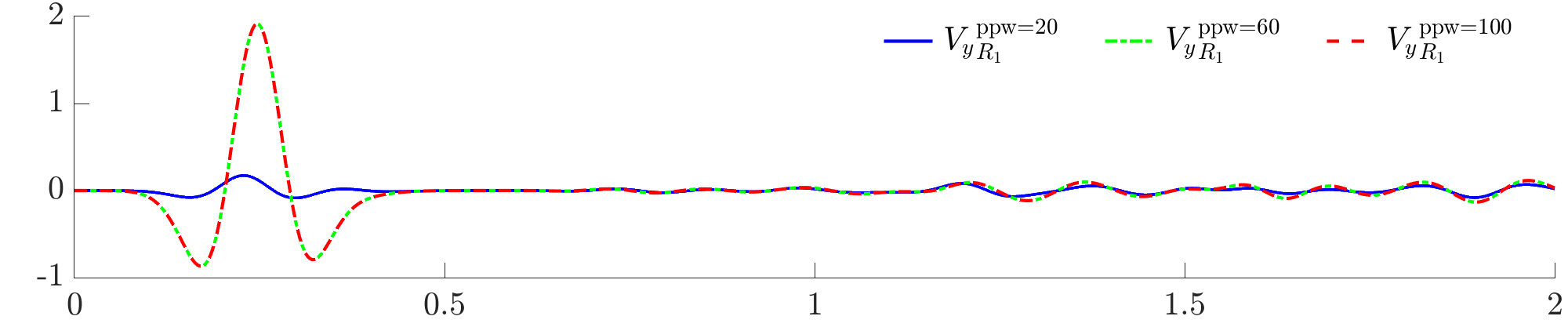}
\caption{Time history of $V_y$ at 1 grid points below the surface for the case ppw = 20.}
\label{stress_test_Sxy_Elastic_full_FSBC_Vy_1_weak}
\end{subfigure}\hfill
\\[2ex]
\begin{subfigure}[b]{1\textwidth}
\captionsetup{width=1\textwidth, font=footnotesize,labelfont=footnotesize}
\centering\includegraphics[scale=0.2]{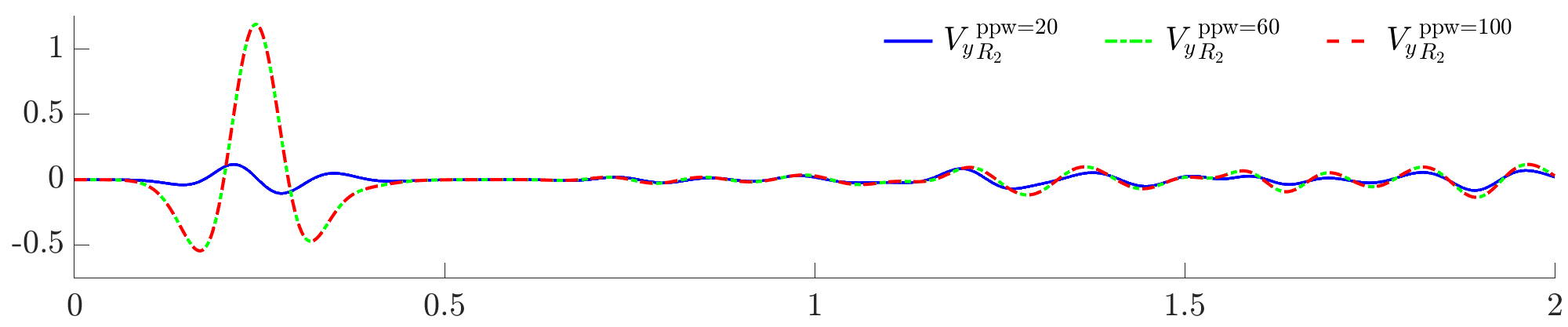}
\caption{Time history of $V_y$ at 2 grid points below the surface for the case ppw = 20.}
\label{stress_test_Sxy_Elastic_full_FSBC_Vy_2_weak}
\end{subfigure}\hfill
\caption{Time histories of the 2D elastic experiments. The free surface boundary condition is imposed weakly.
}
\label{stress_test_Sxy_Elastic_full_FSBC_Vy_weak}
\end{figure}
\end{minipage}